\documentclass[12pt,preprint,nofootinbib,nobibnotes,nonumbib]{article}
\usepackage{amsmath} 
\usepackage{amsxtra} 
\usepackage{amscd}
\usepackage{amsthm} 
\usepackage{amsfonts} 
\usepackage{amssymb} 
\usepackage{eucal} 
\usepackage{bm} 
\usepackage{hyperref}

\newcommand{\C}{{\mathbb C}} 
\newcommand{\Z}{{\mathbb Z}} 
\newcommand{\N}{{\mathbb N}} 
\newcommand{\T}{{\mathbb T}}

\begin{document}

\title{Dense nuclear Fr\'echet ideals in $C^\star$-algebras}

\author{Larry B. Schweitzer}

\date{October 29, 2015}

\maketitle 

\begin{abstract}
We show that a $C^\star$-algebra $B$ contains 
a dense left or right Fr\'echet ideal $A$, 
with $A$ a nuclear locally convex 
space, if and only if the primitive ideal space Prim$(B)$ 
of $B$ is discrete and countable,
and $B/I$ is finite dimensional for each $I \in $ Prim$(B)$.
We show the forward implication holds
for a general Banach algebra $B$,
if the ideal is assumed two-sided. 
For $C^\star$-algebras, we construct 
all two-sided
dense nuclear ideals by defining a set 
of matrix-valued Schwartz functions on 
the countable discrete space Prim$(B)$.

\end{abstract}

\tableofcontents

\pagebreak

\section{Introduction}

Dense subalgebras 
of $C^\star$-algebras are well-known to
be useful in the study of $C^\star$-algebras.
The dense subalgebra can be viewed as $C^\infty$ functions
on a manifold, where instead of a manifold we have an 
underlying \lq\lq noncommutative space\rq\rq.
To be more useful, the dense subalgebra often has
a Fr\'echet topology and is
{\it spectral invariant} in the $C^\star$-algebra
(see Remark 3.2).
While some $C^\star$-algebras have nice dense
subalgebras, others don't seem to.

If we insist that the dense subalgebra be an {\it ideal} in
the $C^\star$-algebra, 
we get a much stronger condition than spectral invariance.
Few $C^\star$-algebras can have a \lq\lq smooth\rq\rq\ dense ideal. 
For example, any compact manifold $M$ without boundary
has a spectral invariant Fr\'echet algebra 
of smooth functions $C^\infty(M)$,
dense in the $C^\star$-algebra of 
continuous functions $C(M)$.   
But $C^\infty(M)$ is an ideal in $C(M)$, if and only
if $M$ is discrete (and therefore finite by compactness).  

In this paper, I will classify which $C^\star$-algebras
have dense nuclear Fr\'echet 
ideals.   The nuclearity property
plays the role of making elements of the subalgebra 
 \lq\lq smooth\rq\rq\ 
or \lq\lq differentiable\rq\rq.  For example, $C^\infty(M)$ 
is a nuclear Fr\'echet space  
[Treves, 1967], Chapters 10 and 51. 

The notion of nuclearity for a locally 
convex space is different than nuclearity 
for a $C^\star$-algebra [Kad Ring II, 1997], Chapter 11.
If a $C^\star$-algebra were nuclear 
as a locally convex space, 
it would be finite dimensional (Proposition 2.1 (b) below).  
However, it is reasonable that 
a dense Fr\'echet subalgebra of an
infinite dimensional
nuclear $C^\star$-algebra be a 
nuclear locally convex space
(see Corollary 4.7 below). 

In \S 2, we recall definitions and properties of nuclear 
locally convex spaces and Fr\'echet spaces from the literature.
We define Schwartz functions on a countable set.
Every nuclear Fr\'echet space with basis is one of these.

In \S 3, we go over definitions and 
basic lemmas on dense Fr\'echet ideals of Banach algebras, and
state the example of $\ell^p(X)$ and Schwartz functions
${\cal S}(X)$ as dense ideals of the commutative
pointwise-multiplication $C^\star$-algebra $c_0(X)$.  
We show the property of
having a dense nuclear ideal
is preserved
by taking quotients and subalgebras formed using
idempotents.

In \S 4, we apply nuclearity and results of \S 3 to see that  
a $C^\star$-algebra
has a dense nuclear left or right ideal only if its
spectrum is discrete and countable.
We do this by showing that 
the \lq\lq finite socle\rq\rq\  is dense.
In the commutative case, a shorter proof is given.

In \S 5, 
we generalize our results on discrete spectrum
to the case of an arbitrary Banach algebra.
We prove that any Banach algebra with a dense nuclear
two-sided ideal is both left and right completely
continuous.  We also show the primitive ideal space
is countable, and that primitive quotients 
are finite dimensional.
Complete continuity on both sides
of a Banach algebra is already known to imply
discrete spectrum and finite dimensional 
primitive quotients [Kaplansky, 1949], [Kaplansky, 1948], Lemma 4.

In \S 6, we construct dense nuclear two-sided ideals
for every $C^\star$-algebra $B$ which 
is the direct sum of simple finite dimensional
$C^\star$-algebras (i.e. full matrix algebras).
We define our dense
ideal to be matrix-valued Schwartz functions on the countable
discrete spectrum of $B$.

In \S 7, we note that both the $C^\star$-algebra and dense 
ideal can be decomposed into direct sums of two ideals,
the matrix algebras which occur with finite
multiplicity and those that occur with infinite multiplicity.
In many cases, the underlying 
Fr\'echet space structure of the dense ideal
is isomorphic to standard Schwartz functions on 
a countably infinite set $X$,
namely 
\begin{equation}
{\cal S}_\gamma (X) = \{ \varphi \colon X 
\rightarrow \C \,|\, \text{for each } n=0,1,2,\dots,  
\,\sup_{x \in X} \gamma(x)^n |\varphi(x)|< \infty  \},
\label{eqn:standSchw}
\end{equation}
where $\gamma$ is an enumeration of $X$.

In \S 8, we give examples of dense nuclear ideals,
including the convolution algebra of $C^\infty$ functions
on a compact Lie group, and
$C^\infty$ functions on the Cantor set.

In Appendix A, we show that the constants $m_n$ and $C_n$
in the dense ideal inequality (\ref{eqn:idealcond}) can be made
to satisfy $m_n = n$ and $C_n = 1$ by choosing an equivalent 
family of increasing norms
$\{\| \cdot \|_n\}_{n=0}^\infty$ for the 
topology of the dense Fr\'echet ideal.

In Appendix B, we give counterexamples to our 
theorems, when various hypotheses are dropped.

Throughout the paper, 
$X$ will denote a countable set (with discrete topology), 
usually infinite. 
All algebras in this paper are over the field 
of complex numbers $\C$.
The set of natural
numbers $\{0, 1, 2, \dots\}$
is denoted by $\N$, and $\N^+ = \{1,2,3,  \dots \}$.

\vskip\baselineskip
\vskip\baselineskip
\section{Nuclearity and Fr\'echet Spaces with Basis}

We state some basics facts about nuclear Fr\'echet spaces
used throughout the paper.
The reader unfamiliar with these concepts can 
consult standard references, such as
[Pietsch, 1972] and [Treves, 1967].

\vskip\baselineskip
\noindent{\bf Proposition 2.1.  Nuclearity Facts.}

\noindent
{\it 
(a) A nuclear Fr\'echet space is separable.

\noindent
(b) A nuclear Banach space must be finite dimensional.

\noindent
(c) If a Fr\'echet space $F$ contains an infinite 
dimensional Banach space $E$, 
and the topology on $E$ inherited from $F$ is
stronger than the Banach space topology, 
then $F$ is not nuclear.

\noindent
(d) A bounded set in a nuclear Fr\'echet space is
relatively compact.
}

\noindent{Proof:} 
For separability, see [Pietsch, 1972], Theorem 4.4.10.
For finite dimensionality of a nuclear normed space,
see [Pietsch, 1972], Theorem 4.4.14, 
or [Treves, 1967], Chapter 50-12, 
Corollary 2.  

Let $F$ be a Fr\'echet space containing an 
infinite dimensional Banach space
$E$.  Since the Banach space topology on $E$ 
is assumed weaker than the Fr\'echet space topology 
inherited from $F$, these two topologies
agree on $E$ [Treves, 1967], Chapter 17-7, Corollary 2.
But a subspace of a nuclear Fr\'echet space is 
nuclear [Pietsch, 1972], Proposition 5.1.5, 
or
[Treves, 1967], Proposition 50.1 (50.3).
$E$ is not nuclear by (b), so therefore 
$F$ cannot be nuclear.  

A bounded set in a nuclear locally convex space is
precompact by [Pietsch, 1972], 4.4.7, or 
[Treves, 1967], Proposition 50.2 (50.12).
\qed

\vskip\baselineskip
\noindent{\bf Remark 2.2.}
The hypotheses of Proposition 2.1 (c) 
imply that $E$ is a closed subspace of $F$, namely ${\overline E}^F = E$.
Here is a counterexample when the Banach space is not a closed
subspace of the Fr\'echet space. 
Let $F=l_u(X)$ be the nuclear Fr\'echet space of
complex-valued sequences on $X$ (allowed to be unbounded), 
topologized by the seminorms
$\| \varphi \|_n = \max_{\gamma(x) \leq n} | \varphi(x)|$ 
(see [Pietsch, 1972], 4.3.4 or [Treves, 1967], Theorem 51.1),
and let $E$ be the Banach space 
$\ell^p(X)$ for any $1\leq p \leq \infty$.

\vskip\baselineskip
We will often cite the following obvious Corollary to Proposition 2.1 (a).

\noindent{\bf Corollary 2.3.} {\it A Banach space
with a dense nuclear Fr\'echet subspace is separable.}

\vskip\baselineskip
For any locally convex space $E$, 
recall that the {\it polar} $M^0$ of a set 
$M \subseteq E$ is defined to be the set
\begin{equation}
\{ \varphi \in E' \,\, |\,\, | \varphi(m) | \leq 1, \quad m \in M \}.  
\end{equation}
Here $E'$ denotes the topological dual of $E$, the set of linear
functionals on $E$ which are continuous.  
We will use the following consequence of the Uniform Boundedness Theorem
in \S 5.  

\noindent{\bf Theorem 2.4.  Uniform Boundedness.}
{\it Let $E$ be a Banach space with unit ball $U$.
If $\{ \varphi_n \}_{n=0}^\infty \subseteq E'$
is a sequence of continuous linear functionals 
converging to zero pointwise,
 then 
$\sup\, \{ \, \| \varphi_n \|_{U^0} \,\, | \,\, n \in \N\, \} < \infty$.}

\noindent{Proof:} See [Rudin, 1973], Theorem 2.6.  \qed

\vskip\baselineskip
\noindent{\bf Definition 2.5.  Fr\'echet Space with Basis.} 
A sequence $\{ e_k \}^\infty_{k=0}$ in a Fr\'echet space $F$ 
is a {\it basis} if every $f \in F$ has a unique series expansion
$f = \sum^\infty_{k=0} f_k e_k$ which converges in $F$, 
where each $f_k$ is a complex number.
By [Husain, 1991], Chapter I, Theorem 4.3,
any basis in a Fr\'echet space is a {\it Schauder basis},
meaning that the coordinate 
functional $f \mapsto f_k$ is a continuous 
linear map from $F$ to $\C$
for each $k \in \N$.
The basis is {\it unconditional} if the series expansion converges
to $f$ in any order.
If  for every $n \in \N$, there is some $C>0$ and $q \in \N$ such that
$|f_k | \| e_k \|_n = 
\| f_k e_k \|_n 
\leq C \|f \|_q$, for all $f \in F$ and $k \in \N$,
then the basis is {\it equicontinuous} [Pietsch, 1972], 10.1.2.
According to [Pietsch, 1972], Theorem 10.1.2, every Schauder basis
in a Fr\'echet space is equicontinuous.
If $F$ is nuclear, an equicontinuous basis is {\it absolute}
[Pietsch, 1972], 10.2.1.  This means that for each $n \in \N$
there exists $C>0$ and $q \in \N$ for which
$\sum_{k=0}^\infty |f_k | \| e_k \|_n  
\leq C \|f \|_q$, for all $f \in F$.
Then the series $\sum_{k=0}^\infty f_k e_k$ converges 
absolutely and unconditionally to $f$.
For each $e_k$, let $e'_k \colon F \rightarrow \C$ denote
the dual coordinate functional $e'_k(f) = f_k$, $f \in F$.

We say two bases $\{ e_k \}_{k=0}^\infty$ and $\{ {\tilde e}_k \}_{k=0}^\infty$
for (possibly different) Fr\'echet spaces
are {\it equivalent} if 
for every sequence of complex numbers $\{ f_k \}_{k=0}^\infty$,
the series $ \sum_{k=0}^\infty f_k e_k$
converges if and only if 
the series $\sum_{k=0}^\infty f_k {\tilde e}_k$ converges,
in their respective Fr\'echet spaces.
The two bases are {\it semi-equivalent} if 
there exists some sequence of non-zero scalars 
$r_k \not= 0$ for which
 $\{ e_k \}_{k=0}^\infty$ and $\{ r_k {\tilde e}_k \}_{k=0}^\infty$
are equivalent [Crone Rob, 1975], \S 1. 
If the bases are absolute, we take $r_k >0$ in the
definition [Dubinsky, 1979], (6.3.4). 

\vskip\baselineskip
\noindent{\bf Definition 2.6.  Scales
and Families of Scales.}
Let $X$ be a countable set.
A function
$\sigma \colon X \rightarrow  [1, \infty)$
is a {\it scale} on $X$.
A scale $\sigma$
{\it dominates} another scale $\tau$ if there
exists $C>0$ such that $\tau(x) \leq C \sigma(x)$ holds
for all $x \in X$.  
We use the notation $\tau \lesssim \sigma$.
If $\tau$ also dominates $\sigma$, we say that $\sigma$
and $\tau$ are {\it equivalent},
and write $\sigma \thicksim \tau$.
The scales $\sigma$ and $\tau$ are {\it semi-equivalent}
if there exists a function $r \colon X \rightarrow (0,\infty)$
for which $\sigma$ and $r\tau$ are equivalent,
and use the notation $\sigma \thicksim_{si} \tau$.

Let $\sigma$ be a scale on $X$ and $\tau$ a scale on
another countable set
$Y$, with $\pi \colon X \cong Y$ a bijection.
Noting that $\tau \circ \pi$ is a scale on $X$,
we say $\sigma \thicksim \tau$ when $\sigma \thicksim
\tau \circ \pi$, and $\sigma \thicksim_{si} \tau$
when $\sigma \thicksim_{si} \tau \circ \pi$, 

If ${\underline \sigma} = \{ \sigma_n \}_{n=0}^\infty$
is a {\it family of scales} on  $X$,
we require that for any pair $\sigma_n, \sigma_m \in \underline \sigma$,
there exists some $\sigma_p \in \underline \sigma$ which dominates 
$\sigma_n$ and $\sigma_m$.
Then we can 
replace $\underline \sigma$ with an equivalent increasing family
of scales, so that 
$ \sigma_0 \leq  \sigma_1
\leq \cdots  \sigma_n \leq \cdots$.
We sometimes adopt the convention
$ \sigma_0\equiv 1$.
A family $\underline \sigma$
dominates another family $\underline \tau$ if
for every $n \in \N$
there exists $m  \in \N$ such that
$ \tau_n \lesssim \sigma_m$.
As for single scales, 
we say that $\underline \sigma$ is equivalent to $\underline \tau$
if both ${\underline \sigma} \lesssim {\underline \tau}$
and ${\underline \tau} \lesssim {\underline \sigma}$.
To define semi-equivalence ${\underline \sigma} 
\thicksim_{si} {\underline \tau}$, we use the same function 
$r \colon X \rightarrow (0,\infty)$
for all scales in the family $\underline \tau$,
and generalize the definitions in the previous paragraphs.
If $\sigma$ is a scale,
an associated family ${\underline \sigma}$ is given
by the powers $\sigma_n = \sigma^n$.

\vskip\baselineskip
\noindent{\bf Definition 2.7.  Sequence Spaces defined using Scales.}
Following the definition of a sequence space
in [Pietsch, 1972], Definition 6.1.1, we define
the linear space
\begin{equation}
{\cal S}^1_{\underline \sigma} (X) =
\biggl\{ \varphi \colon X \rightarrow \C \,
\biggl| \, \| \varphi \|^1_n < \infty\, ,\, n \in \N \biggr\},
\end{equation}
which is a Fr\'echet space for the norms
\begin{equation}
\| \varphi \|^1_n = \sum_{x \in X}
 \sigma_n(x) | \varphi(x) |.
\label{eqn:rapidNorms1}
\end{equation}
Call this Fr\'echet space ${\cal S}^1_{\underline \sigma} (X)$ the
{\it $\ell^1$-norm ${\underline \sigma}$-rapidly vanishing
functions on $X$}, or simply {\it Schwartz functions on $X$}.
Let $\delta_x $ denote the step function
at $x \in X$.
Then span$\{\delta_x \, | \, x \in X \}$
is the dense
subspace of finite support functions in
${\cal S}^1_{\underline \sigma} (X)$,
denoted by
 $c_f(X)$.
The set of functions $\{ \delta_x \}_{x \in X}$
is an absolute basis for
${\cal S}^1_{\underline \sigma} (X)$ since
$\sum_{x \in X} | \varphi(x) | \| \delta_x \|^1_n
= \| \varphi \|_n^1$.
Similarly, define ${\cal S}^\infty_{\underline \sigma}(X)$
using $\sup$-norms in place of $\ell^1$ norms:
\begin{equation}
\| \varphi \|^\infty_n = \sup_{x \in X}
\bigg\{  \sigma_n(x) | \varphi(x)| \biggr\}.
\label{eqn:rapidNormsSup}
\end{equation}
Define the
{\it $\sup$-norm $\underline \sigma$-rapidly vanishing functions
on $X$}, 
${\cal S}^\infty_{\underline \sigma} (X)$, to be the completion
of ${\cal S}^1_{\underline \sigma} (X)$ (or equivalently of
 $c_f(X)$) in the norms (\ref{eqn:rapidNormsSup}).
We may also refer to
${\cal S}^\infty_{\underline \sigma} (X)$ as
 {\it Schwartz functions on $X$}.
If $\varphi = \{ \varphi_k \}_{k=0}^\infty$ is a
Cauchy sequence in this completion,
then $\{ \varphi_k(x) \}_{k=0}^\infty$
is Cauchy in $\C$ for each $x \in X$,
since we are using $\sup$-norms.
The pointwise limit $\varphi(x)$ is the unique coordinate
of $\varphi$ in the unconditional basis $\{ \delta_x \}_{x \in X}$.
Since
$| \varphi(x) | \| \delta_x \|^\infty_n \leq \| \varphi \|_n^\infty$,
the basis is equicontinuous.

The identity map
id$\colon X \rightarrow X$
extends to a continuous inclusion of
Fr\'echet spaces ${\cal S}^1_{\underline \sigma}(X)
\hookrightarrow {\cal S}^1_{\underline \tau}(X)$
or
${\cal S}^\infty_{\underline \sigma}(X)
\hookrightarrow {\cal S}^\infty_{\underline \tau}(X)$
if and only if ${\underline \tau} \lesssim {\underline \sigma}$.
If $\underline \sigma$ is a family of powers of a single
scale $\sigma$ on $X$, we use the notation
${\cal S}^1_{ \sigma}(X)$ and
${\cal S}^\infty_{\sigma}(X)$ to stand for
${\cal S}^1_{\underline \sigma}(X)$
and
${\cal S}^\infty_{\underline \sigma}(X)$, respectively.

\vskip\baselineskip
\noindent{\bf Remark 2.8.  Power Series of Infinite Type.}
Let $\alpha$ be an ordered sequence of
real numbers $0 \leq \alpha_1 \leq \alpha_2 \leq \cdots 
\alpha_k \leq \cdots$.
The set of sequences $ \varphi \colon X \rightarrow \C$
which satisfy
\begin{equation}
\| \varphi \|^1_{\alpha, \rho} =  
\sum_{x \in X} \rho^{\alpha_{\gamma(x)}} | \varphi(x) | < \infty,
\end{equation}
for every real number $\rho$ satisfying $0< \rho < \infty$, is the
{\it $\alpha$-power series space of infinite type on $X$}, where
$\gamma$ is some fixed enumeration of $X$ [Pietsch, 1972], 6.1.5
or [Dubinsky, 1979], Chapter 1, (4.1).
If $\sigma$ is a scale on $X$, we can define $\alpha$ by
$\alpha(x) = \ln(\sigma(x))$.  Then ${\cal S}^1_{\sigma}(X)$
is precisely the $\alpha$-power series space of infinite type, since
the norms
\begin{equation}
\|\varphi \|^1_{\sigma, n} 
= \sum_{x \in X} 
\sigma(x)^n | \varphi(x) |
= \sum_{x \in X} 
e^{n \alpha(x)} | \varphi(x) |
  = \| \varphi \|^1_{\alpha, \rho}
\end{equation}
are equal, taking $\rho = e^n$.

\vskip\baselineskip
\noindent{\bf Theorem 2.9.  Nuclear Fr\'echet Spaces with Basis.}
{\it 
Let $X$ be a countable set, and
$F$ a Fr\'echet space with 
absolute basis $\{e_x\}_{x \in X}$.
Assume there exists a continuous norm $\| \cdot \|_{00}$ on $F$
for which $\| e_x \|_{00} = 1$, $x \in X$.
(This can be arranged by replacing each $e_x$ with $e_x/\|e_x\|_{00}$.) 
Let $\{ \| \cdot\|_n \}_{n=0}^\infty$ be an increasing
family of norms dominating $\| \cdot \|_{00}$ and topologizing $F$.
Then $ \sigma_n(x) = \| e_x \|_n$ 
defines a family of scales ${\underline \sigma}$
on $X$, and the Fr\'echet 
spaces  $F \cong {\cal S}^1_{\underline \sigma}(X)$
are naturally isomorphic.

Moreover, $F$ is nuclear if and only if
$\underline \sigma$ satisfies the  summability condition
\begin{equation} 
(\forall n \in \N)\,(\exists m > n) \quad
\sum_{x \in X}  \,
{ \sigma_n(x) \over  \sigma_m(x)}\, <\, \infty,
\label{eqn:NucSumSpace}
\end{equation}
and if and only if the Fr\'echet spaces  
$F \cong {\cal S}^1_{\underline \sigma}(X) \cong 
{\cal S}^\infty_{\underline \sigma}(X)$
are naturally isomorphic.
}

\noindent{Proof:} 
The isomorphism 
$F \cong {\cal S}^1_{\underline \sigma}(X)$
is [Pietsch, 1972], Theorem 10.1.4.
The second paragraph is [Pietsch, 1972], Theorems 6.1.2 and 6.1.3.
\qed

\vskip\baselineskip
\noindent{\bf Remark 2.10.  Defining
${\cal S}^\infty_{\underline \sigma} (X)$ without completion.}
We say that a map $f \colon X \rightarrow (0, \infty)$
is {\it proper} if
$1/f$ vanishes at infinity.
If the summability condition
(\ref{eqn:NucSumSpace}) holds for $\underline \sigma$,
then $ \sigma_m / \sigma_n$ is a proper
map $X \rightarrow (0, \infty)$, for any pair
$m$ and $n$ for which the summation holds.
We show that when for any $n \in \N$,
there is an $m \geq n$ for which
$ \sigma_m/ \sigma_n$ is proper, then
${\cal S}^\infty_{\underline \sigma}(X)$ is
identical to the set of functions
\begin{equation}
\biggl\{ \varphi \colon X \rightarrow \C \,
\biggl| \, \| \varphi \|^\infty_n < \infty\, ,\, n \in \N \biggr\}.
\end{equation}
For let $\varphi$ be in this set.
Since $| \sigma_n\varphi(x)| \leq
{ \sigma_n \over  \sigma_m}(x)
\| \varphi \|_m^\infty $,
the function $ \sigma_n \varphi$
vanishes at infinity on $X$.
If $S$ is a finite subset of $X$,
then
$\bigl\| \sum_{x \in S} \varphi(x) \delta_x(\cdot)
- \varphi(\cdot) \bigr\|_n^\infty
= \sup_{x \in X - S} | \sigma_n\varphi(x)|$.
The right hand side can be as small as
we like, by
taking larger and larger finite sets $S$.
Since this is true for each $n \in \N$,
we have shown that $\varphi$
is approximated by elements of
$c_f(X)$ in the $\sup$-norms (\ref{eqn:rapidNormsSup}),
and hence $\varphi \in {\cal S}^\infty_{\underline \sigma} (X)$.

\vskip\baselineskip
\noindent{\bf Definition 2.11. Enumerations and Regularity.}
We call a scale $\gamma$ an {\it enumeration} of a
countably infinite set $X$ 
if $\gamma \colon X \cong \N^+$ is a bijection.
If $\varphi$ is a complex-valued function on $X$, we write
$\varphi = (\varphi_1, \varphi_2, \varphi_3, \dots \varphi_k, \dots)$,
where $\varphi_k$ is short for $\varphi(\gamma^{-1}(k))$.
The step functions $\delta_x$ for $x \in X$ 
may be written $\delta_k$,
where $k=\gamma(x)$.
We write
 $(x_1, x_2, \dots x_k, \dots)$
in place of
 $(\gamma^{-1}(1), \gamma^{-1}(2), \dots \gamma^{-1}(k), \dots)$.

%If $\underline \sigma$ is a family of scales on $X$, 
%then $\underline{\gamma \sigma}$ denotes
%a new family of scales on $X$, 
%given by $\{ \gamma^n  \sigma_n\}_{n=0}^\infty$.

If ${\underline \gamma} = \{ \gamma_n \}_{n=1}^\infty$ is a family
of enumerations on $X$, we define
the {\it ${\underline \gamma}$-family of scales on $X$}
by
${\underline \sigma}_{\underline \gamma} = \{ \gamma_1 
\cdots \gamma_n \}_{n=0}^\infty$.
Also define 
${\underline \sigma}_{\underline \gamma}^2$ 
and $\sqrt{{\underline \sigma}_{\underline \gamma}}$,
the {\it ${\underline \gamma}$-family of squared and 
square-root scales on $X$},
in the obvious way.

If $\sigma$ is a proper scale on $X$, it is not hard to see
there exists an enumeration of $X$ for which
$\sigma(x_1) \leq \sigma(x_2) \leq \cdots \sigma(x_k) \leq \cdots$.

An absolute basis $\{ e_x \}_{x \in X}$
for some Fr\'echet space is {\it $\gamma$-regular}
if $\sigma_n(x_k)/\sigma_m(x_k) \geq 
\sigma_n(x_{k+1})/\sigma_m(x_{k+1})$ for $m \geq n$ and
$k \in \N^+$ [Dubinsky, 1979], Chapter 1, (6.3.1),
where the family of scales $\underline \sigma$ is defined by
$\sigma_n(x_k) = \| e_{x_k} \|_n$.  
In this framework, we also say $\underline \sigma$
is $\gamma$-regular.
Note that if $\underline \sigma$
is powers of a single scale, $\sigma_n = \sigma^n$, then
$\gamma$-regular means that $\sigma$ is {\it $\gamma$-ordered}, 
that is $\sigma(x_1) \leq \sigma(x_2) \leq \cdots \sigma(x_k) \leq \cdots$.

\vskip\baselineskip
\noindent{\bf Proposition 2.12. Enumerations and
Summability.}
{\it 
Let ${\underline \gamma}$ be a family of enumerations of $X$.
The families ${\underline \sigma}_{\underline \gamma}$,
${\underline \sigma}_{\underline \gamma}^2$, and 
$\sqrt{{\underline \sigma}_{\underline \gamma}}$ satisfy
summability with $m=n+2$, $m=n+1$, and $m=n+4$, respectively.

A family of scales $\underline \sigma$ is summable
if and only if
for every $n \in \N$ there exists 
an enumeration $\gamma$ of $X$ such
that $\underline \sigma$ dominates the scale
$\sigma_n \sqrt{\gamma}$.
In particular, if $\underline \sigma$ is summable,
then $\underline \sigma$ dominates 
$\sqrt{{\underline \sigma}_{\underline \gamma}}$ for 
some sequence of enumerations ${\underline \gamma}$.

If we make the additional assumption that
$\underline \sigma$ is $\gamma$-regular
for some enumeration $\gamma$ of $X$, then the
summability of $\underline \sigma$
 implies that $\underline \sigma$ dominates
every power of $\gamma$. 
}

\noindent{Proof:} 
First note that
\begin{equation}
{  \bigl(\sigma^2_{\underline \gamma}\bigr)_n \over 
{  \bigl(\sigma^2_{\underline \gamma}\bigr)_{n+1} } }
\,=\,
{ 1 \over{\gamma_{n+1}^2} }.
\end{equation}
Then $\sum_{x \in X} 1/ \gamma_{n+1}(x)^2
= \sum_{k=1}^\infty 1/k^2 = \pi^2/6 < \infty$
shows that ${\underline \sigma}^2_{\underline \gamma}$ satisfies
summability with $m = n+1$.

Next note that
\begin{equation}
{  \bigl(\sigma_{\underline \gamma}\bigr)_n  \over 
{  \bigl(\sigma_{\underline \gamma}\bigr)_{n+2} }  }
\,=\,
{ 1 \over{ \gamma_{n+1} \gamma_{n+2} }}.
\end{equation}
Then 
$\sum_{x \in X} 1/ \gamma_{n+1}(x) \gamma_{n+2}(x)
\, \leq\, \| 1/ \gamma_{n+1} \|_2\, 
\| 1/ \gamma_{n+2} \|_2 
\, =\, \pi^2/6\, <\, \infty$
shows that ${\underline \sigma}_{\underline \gamma}$ satisfies
summability with $m = n+2$, where we used the Cauchy-Schwartz
inequality. 
Similarly $\sqrt{{\underline \sigma}_{\underline \gamma}}$ 
satisfies summability
with $m = n+4$.

Next let $\underline \sigma$ be a family satisfying the
condition of the second paragraph.
Apply the condition four times 
to get four enumerations $\gamma_1$, $\gamma_2$, $\gamma_3$,
$\gamma_4$ of $X$ such that
$\sqrt{\gamma_1} \leq C_1  \sigma_m/ \sigma_n$, 
$\sqrt{\gamma_2} \leq C_2  \sigma_p/ \sigma_m$, 
$\sqrt{\gamma_3} \leq C_3  \sigma_q/ \sigma_p$, 
and
$\sqrt{\gamma_4} \leq C_4  \sigma_r/\sigma_q$. 
Then, using the Cauchy-Schwartz inequality
and setting $C=C_1 C_2 C_3 C_4$, we have
\begin{eqnarray}
\sum_{x \in X} { \sigma_n(x) \over  \sigma_r(x)} 
\quad & = & \quad 
\sum_{x \in X} { \sigma_n \over  \sigma_m} 
{ \sigma_m \over  \sigma_p} 
{ \sigma_p \over  \sigma_q} 
{ \sigma_q \over  \sigma_r} (x)  \nonumber \\
\quad & \leq & \quad  C \sum_{x \in X}   
\biggl( {1 \over {\gamma_1} 
{\gamma_2} {\gamma_3} {\gamma_4} }(x) \biggr)^{1/2}
\nonumber \\
\quad & \leq & \quad  C\bigl(\,\, \| 1/{\gamma_1}\|_2   
\,\, \| 1/{\gamma_2}\|_2   
\,\, \| 1/{\gamma_3}\|_2   
\,\, \| 1/{\gamma_4}\|_2 \,\, \bigr)^{1/2} \nonumber \\ 
\quad &=& \quad 
C \biggl( \sqrt{ \pi^2 \over 6} \sqrt{ \pi^2 \over 6}
\sqrt{ \pi^2 \over 6} \sqrt{ \pi^2 \over 6}
\biggr)^{1/2} \quad = \quad
 {C \pi^2 \over 6} \quad < \quad \infty, 
\qquad
\end{eqnarray}
so $\underline \sigma$ is summable.

Conversely, find an enumeration $\gamma$ of $X$ such that
$( \sigma_n/ \sigma_m) \circ \mu$ is 
non-increasing, where $\mu = \gamma^{-1}$.
Apply the Cauchy Condensation Test [Mars Hoff, 1999], p. 194,
for convergence of the series
$\sum^\infty_{k =1} ( \sigma_n/ \sigma_m) \circ \mu(k)$, 
to see that $\sum^\infty_{i=0}  2^i 
(  \sigma_n / 
 \sigma_m ) \circ \mu(2^i) < \infty$.
Let $C>0$ be such that
$ \sigma_n/ \sigma_m \circ \mu(2^i) \leq C/2^i$, $i \in \N$.   
Since the terms are non-increasing,
$ \sigma_n/ \sigma_m \circ \mu(k) \leq 
C/2^i = C/\sqrt{2^{2i}}$ for $2^i \leq k < 2^{2i}$,
so $ \sigma_n/ \sigma_m \circ 
\mu(k) \leq C/\sqrt{k}$, $k \in \N^+$.
Replacing $k$ with $\gamma(x)$, we see
$ \sigma_n(x) \sqrt{\gamma(x)} \leq C \sigma_m(x)$, $x \in X$.

To finish the second paragraph, let 
${\underline \gamma} = \gamma_1, \gamma_2, \dots, \gamma_n, \dots$
be a family of enumerations, 
and $\{ m_n \}_{n=0}^\infty$ an increasing sequence of natural
numbers, with the property that
$ \sigma_1 \sqrt{ \gamma_1 } \lesssim  \sigma_{m_1}$,
$ \sigma_{m_1} \sqrt{ \gamma_2 } \lesssim  \sigma_{m_2}$,
$\dots$,
$ \sigma_{m_{n-1}} \sqrt{ \gamma_n } \lesssim  \sigma_{m_n}$, 
and so on.
Then $\bigl(\sqrt{ {\underline \sigma}_{\underline \gamma} } 
\bigr)_n
= \sqrt{\gamma_1 \dots \gamma_n}
\leq  \sigma_1 \sqrt{\gamma_1 \dots \gamma_n}
\lesssim \sigma_{m_1} \sqrt{\gamma_2 \dots \gamma_n}
\lesssim \cdots \lesssim \sigma_{m_{n-1}} \sqrt{\gamma_n}
\lesssim \sigma_{m_n}$,
so $\underline \sigma$ dominates 
$\sqrt{{\underline \sigma}_{\underline \gamma}}$.

Now assume $\gamma$-regularity.
Then 
$( \sigma_n/ \sigma_m) \circ \mu$ is 
non-increasing for any $m, n$ satisfying $m \geq n$, 
where $\mu = \gamma^{-1}$.
For those $m > n$ for which the series is summable,
we have just seen 
there exists constants $C_{m,n}>0$ such
that $\sigma_n(x) \sqrt{\gamma(x)} \leq
C_{m,n} \sigma_m(x)$.
For $d \in \N^+$,
find a sequence of natural numbers
$0 < m_1 < m_2<  \dots < m_{2d}$ for which
\begin{eqnarray}
\sigma_{0}(x) \sqrt{ \gamma(x) } & \leq & C_{m_1,0} \sigma_{m_1}(x),\nonumber\\
\sigma_{m_1}(x) \sqrt{ \gamma(x) } 
&\leq & C_{m_2,m_1} \sigma_{m_2}(x),\nonumber\\
& \vdots & \nonumber \\
\sigma_{m_{2d-1}}(x) \sqrt{ \gamma(x) } & \leq  &
C_{m_{2d},m_{2d-1}} \sigma_{m_{2d}}(x). 
\end{eqnarray}
Combining these $2d$ inequalities and canceling 
the telescoping terms gives
\begin{eqnarray}
\gamma(x)^d &=& \sqrt{\gamma(x)}^{2d} 
 \nonumber \\
&\leq & C_{m_1,0} C_{m_2,m_1} \cdots C_{m_{2d}, m_{2d-1}}
\sigma_{m_{2d}}(x)/\sigma_0(x)\nonumber \\ 
&=& C \sigma_{m_{2d}}(x)/\sigma_0(x) 
\qquad  \mbox{$C=$  product of constants} \nonumber \\
&\leq & C \sigma_{m_{2d}}(x)
\qquad\qquad\quad   \mbox{since $\sigma_0 \geq 1$}, 
\end{eqnarray}
so $\underline \sigma$ dominates $\gamma^d$.
\qed

\vskip\baselineskip
\noindent{\bf Corollary 2.13. 
Single Scales and Summability.}
{\it  
The family of scales associated with
a single scale $\sigma$ is summable if and only if 
$\sigma$ dominates 
every enumeration $\gamma$ of $X$ which orders $\sigma$.
}

\noindent{Proof:} 
If the scale $\sigma$ is summable,
it is proper and there is an enumeration $\gamma$ of
$X$ which orders $\sigma$.
Then $\sigma$ is $\gamma$-regular, and by the last statement
of Proposition 2.12, $\sigma$ dominates every power of $\gamma$.
On the other hand, if there exists $C, d$ such that 
$\gamma(x) \leq C \sigma(x)^d$, $x \in X$, then
\begin{equation}
{1\over C^2} \sum_{x \in X} {1 \over \sigma(x)^{2d}}
\leq
\sum_{x \in X} {1 \over \gamma(x)^2} = \sum_{k=1}^\infty {1\over k^2}
= \pi^2/6 < \infty,
\end{equation}
and $\sigma$ is summable.
\qed

\vskip\baselineskip
\noindent{\bf Proposition 2.14. 
Equivalence and Semi-Equivalence.}
{\it  
Let $\underline \sigma$ and $\underline \tau$ be
 families of scales on sets $X$ and $Y$, respectively,
and let $\pi$ be a bijection $\pi \colon X \cong Y$.
Then $\underline \sigma$ is equivalent (semi-equivalent)
to $\underline \tau$ if and only if
the basis $\{ \delta _x \}_{x \in X}$
for
 ${\cal S}^1_{\underline \sigma}(X)$ 
is equivalent (semi-equivalent) to
the basis $\{ \delta _y \}_{y \in Y}$
of ${\cal S}^1_{\underline \sigma}(Y)$. 

The map
$\pi \colon X \cong Y$ 
extends to an isomorphism of
Fr\'echet spaces ${\cal S}^1_{\underline \sigma}(X)
\cong {\cal S}^1_{\underline \tau}(Y)$
if and only if ${\underline \sigma}$ is equivalent
to ${\underline \tau}$.

If ${\underline \sigma}$ is semi-equivalent
to ${\underline \tau}$,
the map
$\pi \colon X \cong Y$
extends to an isomorphism of
Fr\'echet spaces ${\cal S}^1_{\underline \sigma}(X)
\cong {\cal S}^1_{\underline \tau}(Y)$,
when scaled by the semi-equivalence function $r$.

Assume the family $\underline \sigma$ satisfies the summability condition,
and is regular for some enumeration of $X$.
If ${\cal S}^1_{\underline \sigma}(X)$ is isomorphic to a Fr\'echet
space $E$, then $E$ is identified with 
${\cal S}^1_{\underline \tau}(Y)$ for some countable set $Y$ and
summable family $\underline \tau$ which is regular
for some enumeration of $Y$, with $\underline \sigma$
semi-equivalent to $\underline \tau$.
}

\noindent{Proof:} 
{\bf Equivalence.}
Assume that $\underline \sigma$ and $\underline \tau$
are equivalent.  Then
\begin{eqnarray}
\biggl\| \sum_{x \in X}  f_x \delta_x \biggr\|_d &=&
\sum_{x \in X} |f_x| \sigma_d(x)
\nonumber \\
& \leq & C_d \sum_{x \in X} |f_x| \tau_{m_d}(\pi(x)) 
\qquad \mbox{since $\tau \circ \pi$ dominates $\sigma$}
\nonumber \\
& = & C_d \sum_{y \in Y} |f_{\pi^{-1}(y)}| \tau_{m_d}(y) \nonumber \\
& = & C_d 
\biggl\| \sum_{y \in Y}  f_{\pi^{-1}(y)} \delta_y \biggr\|_{m_d}. 
\end{eqnarray}
A similar inequality can be derived in the reverse direction.
Therefore the tail of a series $\sum_{x \in X} f_x \delta_x$
converges to zero 
in ${\cal S}^1_{\underline \sigma}(X)$
if and only if the corresponding 
tail of the series $\sum_{y \in Y} f_{\pi^{-1}(y)} \delta_y$ 
converges to zero
in ${\cal S}^1_{\underline \tau}(Y)$, 
so the two bases are equivalent.

To prove that 
the map $\pi \colon X \cong Y$ extends to an isomorphism
if the two bases 
$\{ \delta_x \}_{x \in X}$
and
$\{ \delta_y \}_{y \in Y}$
are equivalent, we give a general argument that
equivalent bases imply isomorphic Fr\'echet spaces.
This is an application of 
the Uniform Boundedness Theorem [Rudin, 1973], Theorem 2.6.
Let  $\{ e_k \}_{k=0}^\infty$ and $\{ {\tilde e}_k \}_{k=0}^\infty$
be equivalent bases for two possibly 
different Fr\'echet spaces $E$ and $\tilde E$ (see
Definition 2.5).  Define a linear map
$T \colon E \rightarrow {\tilde E}$  by $T(e_k) = {\tilde e}_k$.
Let $P_n \colon E 
\rightarrow \mbox{span}\{e_0, \dots e_{n-1} \}\subseteq E$
be the projection $P_n\bigl(\sum_{k=0}^\infty f_k e_k\bigr)  = 
\bigl(\sum_{k=0}^{n-1} f_k e_k\bigr)$.  
Each $P_n$ is a continuous linear 
map since $\{ e_k \}_{k=0}^\infty$ is a basis for $E$.
Since any linear map of finite dimensional spaces
is continuous, the composition $T \circ P_n$ is
continuous for each $n \in \N^+$.
Fix $e \in E$, and let $e = \sum_{k=0}^\infty f_k e_k$
be the expansion for $e$.  By assumption of equivalent bases,
$\sum_{k=0}^\infty f_k T(e_k) = \sum_{k=0}^\infty f_k {\tilde e}_k$ 
is a convergent series in $\tilde E$.
Then the partial sums
\begin{equation}
T \circ P_n (e)
= 
T \biggl( \sum_{k=0}^{n-1} f_k e_k \biggr) 
= 
\sum_{k=0}^{n-1} f_k T(e_k)
=
\sum_{k=0}^{n-1} f_k {\tilde e}_k 
\end{equation}
are bounded in $\tilde E$, as $n$ varies.
Applying uniform boundedness, we get
for each seminorm $\| \cdot \|_d$ on $\tilde E$, a
seminorm $\| \cdot \|_{m_d}$ on $E$ such that
\begin{equation}
\| T \circ P_n(e) \|_d \leq C_d \| e \|_{m_d}, 
\end{equation}
for {\it all} $e \in E$, where $m_d$ and $C_d>0$ do not
depend on $n$.
Let $x_n \rightarrow 0$ in $E$.  Choose $M$ large enough
so that $\| x_m \|_{m_d} < \epsilon$ for $m > M$.
Then for fixed $m > M$,
\begin{eqnarray}
\| T(x_m) \|_d \quad & = & \quad \biggl\| T \circ P_n(x_m)  
 \quad  +  \quad \sum_{k=n}^\infty f_{m,k} {\tilde e}_k \biggr\|_d
\qquad \mbox{for any $n\in \N^+$} \nonumber \\
& \leq &  \quad C_d \epsilon \quad  +   \quad
\biggl\| \sum_{k=n}^\infty f_{m,k} {\tilde e}_k \biggr\|_d,
\label{eq:convSeries}
\end{eqnarray}
where $\sum_{k=0}^\infty f_{m,k} e_k$ is the series expansion
for $x_m$ in $E$.  
Since $m$ is fixed, the tail of the series in (\ref{eq:convSeries}) 
tends to zero for large $n$.
Thus we have $\|T(x_m) \|_d \leq C_d \epsilon + \epsilon$ for
$m > M$, and $T$ is a continuous map of Fr\'echet spaces.
The same argument in reverse shows that $T^{-1}$ is continuous,
so $T$ is an isomorphism of Fr\'echet spaces.

Now assume the map
$\pi \colon X \cong Y$ 
extends to an isomorphism of
Fr\'echet spaces ${\cal S}^1_{\underline \sigma}(X)
\cong {\cal S}^1_{\underline \tau}(Y)$.
Then for every $d \in \N$,
there are $C_d> 0$ and $m_d \in \N$ for which
$\| \varphi \circ \pi^{-1} \|_d \leq C_d \| \varphi \|_{m_d}$,
$\varphi \in {\cal S}^1_{\underline \sigma}(X)$.
Taking $\varphi = \delta_x$ gives
\begin{equation}
\tau_d(\pi(x)) = \| \delta_{\pi(x)} \|_d =  
\| \delta_x \circ \pi^{-1} \|_d \leq C_d \| \delta_x \|_{m_d}
= C_d \sigma_{m_d}(x),
\end{equation}
so $\underline \sigma$ dominates $\underline \tau$.
(This spells out the \lq\lq only if\rq\rq\  direction of the remark
at the end of Definition 2.7.)
Similarly $\underline \tau$ dominates $\underline \sigma$,
so the two families of scales are equivalent.
This completes the proof of statements about equivalence.

{\bf Semi-Equivalence.}
Let $\underline \sigma$ and $\underline \tau$ be semi-equivalent
families of scales, with $r \colon X \rightarrow (0, \infty)$
and $\underline \sigma \thicksim r \cdot \underline \tau \circ \pi$.
Then
\begin{eqnarray}
\biggl\| \sum_{x \in X}  f_x \delta_x \biggr\|_d &=&
\sum_{x \in X} |f_x| \sigma_d(x)
\nonumber \\
& \leq & C_d \sum_{x \in X} |f_x| r(x) \tau_{m_d}(\pi(x)) 
\qquad \mbox{since $r \cdot \tau \circ \pi$ dominates $\sigma$}
\nonumber \\
& = & C_d 
\sum_{y \in Y} |f_{\pi^{-1}(y)}|r(\pi^{-1}(y)) \tau_{m_d}(y) \nonumber \\
& = & C_d 
\biggl\| \sum_{y \in Y}  f_{\pi^{-1}(y)} r(\pi^{-1}(y))
\delta_y \biggr\|_{m_d}. 
\end{eqnarray}
A similar inequality can be derived in the reverse direction.
Therefore the tail of a series $\sum_{x \in X} f_x \delta_x$
converges to zero 
in ${\cal S}^1_{\underline \sigma}(X)$
if and only if the corresponding 
tail of the series  $\sum_{y \in Y} f_{\pi^{-1}(y)} r(\pi^{-1}(y)) \delta_y$ 
converges to zero
in ${\cal S}^1_{\underline \tau}(Y)$, 
so the two bases $\{ \delta_x \}_{x \in X}$ 
and $\{ r\circ \pi^{-1}  \cdot \delta_y \}_{y \in Y}$ are equivalent.
The two bases $\{ \delta_x \}_{x \in X}$
and $\{ \delta_y \}_{y \in Y}$ are then semi-equivalent,
by definition.

Next assume the conclusion we just reached, that
the two bases $\{ \delta_x \}_{x \in X}$ 
and $\{ r\circ \pi^{-1}  \cdot \delta_y \}_{y \in Y}$ are equivalent.
The proof we gave of the second paragraph of the proposition, 
gives an isomorphism of Fr\'echet spaces 
${\cal S}^1_{\underline \sigma}(X)$
and ${\cal S}^1_{\underline \tau}(Y)$,
where the isomorphism is given by 
\begin{eqnarray}
\theta(\delta_x) &=& r \circ \pi^{-1} \cdot \delta_y \nonumber \\
& \mbox{or equivalently} \nonumber \\
\theta(\varphi )(y) & = & (r\varphi) \circ \pi^{-1}(y),
\end{eqnarray}
for $\varphi \in {\cal S}^1_{\underline \sigma}(X)$.
Then for every $d \in \N$,
there are $C_d> 0$ and $m_d \in \N$ for which
$\| \theta(\varphi) \|_d \leq C_d \| \varphi \|_{m_d}$.
Taking $\varphi = \delta_x$ for fixed $x$ gives
\begin{equation}
r(x) \cdot \tau_d(\pi(x)) = \| r(x) \delta_{\pi(x)} \|_d =  
\| r\delta_x \circ \pi^{-1} \|_d = 
\| \theta(\varphi)\|_d \leq C_d \| \varphi \|_{m_d}
= C_d \sigma_{m_d}(x).
\end{equation}
Since this is true for every $x \in X$,
$\underline \sigma$ dominates $r \underline \tau$.
Similarly $r \underline \tau$ dominates $\underline \sigma$,
so $\underline \sigma$ and $\underline \tau$ are
semi-equivalent.

Assume $\underline \sigma$ is summable, and $\gamma$-regular
for some enumeration $\gamma$ of $X$.  Then if
${\cal S}^1_{\underline \sigma}(X)$ is isomorphic to some
Fr\'echet space $E$,
[Crone Rob, 1975], Theorem, says there exists
a regular basis $\{ e_k \}_{k=1}^\infty$ 
for $E$ which is semi-equivalent
to the basis $\{ \delta_x \}_{x \in X}$ for
${\cal S}^1_{\underline \sigma}(X)$
(see also [Dubinsky, 1979], (6.3.5)).
By the nuclearity of $E$, we make the identification
$E \cong {\cal S}^1_{\underline \tau}(Y)$
for some countable set $Y$, by Theorem 2.9, where
$\{ e_k \}_{k=1}^\infty$ is identified with the basis
$\{ \delta_y \}_{y \in Y}$.
Let $\zeta$ be the enumeration
of $Y$ for which $e_{\zeta(y)} = \delta_y$.
Then $\tau_n(y) = \| \delta_y \|_n$ with
$\tau$ $\zeta$-regular.
The semi-equivalence of 
 $\{ e_k \}_{k=1}^\infty$ 
with $\{ \delta_x \}_{x \in X}$ says there exists
some 
function $r \colon \N^+ \rightarrow (0, \infty)$ 
such that 
 $\{ \delta_{\gamma^{-1}(k)} \}_{k=1}^\infty$ 
and $\{ r_k \delta_{\zeta^{-1}(k)} \}_{k=1}^\infty$
are equivalent bases for
${\cal S}^1_{\underline \sigma}(X)$ and
${\cal S}^1_{\underline \tau}(Y)$, respectively.
Since
 $\{ \delta_{\gamma^{-1}(k)} \}_{k=1}^\infty$ 
is the same as $\{ \delta_x \}_{x \in X}$ and
$\{ r_k \delta_{\zeta^{-1}(k)} \}_{k=1}^\infty$
is the same as
 $\{ r\circ \pi^{-1} \cdot \delta_y \}_{y \in Y}$, where
$\pi = \zeta^{-1} \circ \gamma \colon X \rightarrow Y$ 
and $r(x) = r_{\gamma(x)}$,
the semi-equivalence of $\underline \sigma$ with
 $\underline \tau$ is what was shown
in the previous paragraph.
\qed

\vskip\baselineskip
\noindent{\bf Example 2.15.  Semi-Equivalence Classes of Scales.}
Let $\gamma$ be an enumeration of $X$.
We show that powers of the two scales $\gamma$ and $e^\gamma$
are not semi-equivalent.
Assume for a contradiction there is some 
$r \colon X \rightarrow (0,\infty)$
for which the families 
$\underline \gamma$ and $r \underline e^\gamma$
are equivalent.
Then for each $d \in \N^+$, we have numbers $n_d, m_d \in \N$ 
and constants $C_d, D_d > 0$ such that
\begin{eqnarray}
\gamma(x)^d & \leq &  C_d\, r(x)\,  e^{n_d\gamma(x)}
\nonumber \\
r(x)\, e^{d\gamma(x)} & \leq &  D_d\, \gamma(x)^{m_d},
\end{eqnarray}
for $x \in X$.
By the first equation with $d=1$,
the reciprocal $1/r(x)$ is bounded by $C \gamma(x)^{-1} e^{n\gamma(x)}$.
Then the second equation says
for each $d \in \N$,
\begin{equation}
e^{d\gamma(x)} \leq {D_d \gamma(x)^{m_d} \over r(x)} \leq
D_d \gamma(x)^{m_d} \bigl( C \gamma(x)^{-1} e^{n\gamma(x)} \bigr)
= C D_d \gamma(x)^{m_d-1} e^{n \gamma(x)},
\end{equation}
for $x \in X$.
For fixed $d > n$, we have
$e^{(d-n)\gamma(x)} \leq C D_d \gamma(x)^{m_d -1}$, which
is a contradiction since  $e^{(d-n) k}$ is not bounded
by $k^{m_d -1}$ for $k \in \N^+$.

By Proposition 2.14, there is no Fr\'echet space isomorphism
between ${\cal S}^1_\gamma(X)$ and ${\cal S}^1_{e^\gamma}(X)$.

\vskip\baselineskip
\noindent{\bf Proposition 2.16. 
Single Scales and Semi-Equivalence.}
{\it  
The families of scales associated with
single scales $\sigma$, $\tau$ are semi-equivalent 
if and only if $\sigma$ is equivalent to $\tau$.
}

\noindent{Proof:} 
Recall that semi-equivalence means
that for all $d \in \N$ we have numbers $n_d, m_d \in \N$ 
and constants $C_d, D_d > 0$ such that
\begin{eqnarray}
 \sigma(x)^d & \leq &  C_d\, r(x) \, \tau(x)^{n_d}
 \nonumber \\
r(x) \, \tau(x)^d & \leq &  D_d\, \sigma(x)^{m_d},
\label{eq:gammaTauSemi}
\end{eqnarray}
for $x \in X$, where $r\colon X \rightarrow (0, \infty)$.
By the first equation with $d=1$,
the reciprocal $1/r(x)$ is bounded by $C_1 \sigma(x)^{-1} \tau(x)^{n_1}$.
Then the second equation says
for each $d \in \N$,
\begin{equation}
\tau(x)^d \leq {D_d \sigma(x)^{m_d} \over r(x)}
\leq D_d \sigma(x)^{m_d}\bigl( C_1 \sigma(x)^{-1} \tau(x)^{n_1} \bigr)
= C_1 D_d \sigma(x)^{m_d -1} \tau(x)^{n_1}.
\end{equation}
Plugging $d=n_1 +1$ into this equation, and dividing both sides by
$\tau(x)^{n_1}$, gives
\begin{equation}
\tau(x) 
\leq C_1 D_{n_1+1} \sigma(x)^{m_{n_1+1} -1},
\end{equation}
so $\tau \lesssim \sigma$.

By the second equation 
of (\ref{eq:gammaTauSemi})
with $d=1$,
 $r(x)$ is bounded by $D_1 \tau(x)^{-1} \sigma(x)^{m_1}$.
Then the first equation 
of (\ref{eq:gammaTauSemi})
says
for each $d \in \N$,
\begin{equation}
\sigma(x)^d \leq C_d r(x) \tau(x)^{n_d} 
\leq C_d \tau(x)^{n_d}\bigl( D_1 \tau(x)^{-1} \sigma(x)^{m_1} \bigr)
= C_d D_1 \tau(x)^{n_d -1} \sigma(x)^{m_1}.
\end{equation}
Plugging $d=m_1 +1$ into this equation, and dividing both sides by
$\sigma(x)^{m_1}$, gives
\begin{equation}
\sigma(x) 
\leq C_{m_1+1} D_1 \sigma(x)^{n_{m_1+1} -1},
\end{equation}
so $\sigma \lesssim \tau$.
\qed

\vskip\baselineskip
\vskip\baselineskip
\section{Fr\'echet Ideals}

\noindent{\bf Definition 3.1. Fr\'echet Ideals and Continuous Inclusion.}
Let $B$ be a Banach algebra, with norm $\|\cdot \|_B$,
and $A$ a subalgebra of $B$. 
The algebra $A$ is called a 
{\it Fr\'echet algebra} when endowed with a 
locally convex Fr\'echet space 
topology for which multiplication is 
jointly continuous. 
Let $\{ \| \cdot \|_n \}_{n=0}^\infty$ be an increasing
family of seminorms giving the topology for $A$.  
When we say that $A$ is a {\it Fr\'echet subalgebra} of $B$, 
we require the inclusion map $A \hookrightarrow B$ be continuous.
In terms of seminorms, this means that if $n \in \N$ 
is sufficiently large, 
there exists a constant $C > 0$ such that
\begin{equation}
\| a \|_B \leq C \|a \|_n, 
\label{eqn:contincl}
\end{equation}
for $a \in A$. 
Define new norms $\{\| \cdot \|'_n \}_{n=0}^\infty$
on $A$ by $\| a \|'_0 = \| a \|_B$ and 
$\| a \|'_{n+1}= \max \{  \| a \|_n, \| a \|_B \}$.
By continuous inclusion (\ref{eqn:contincl}), 
the \lq\lq primed\rq\rq\  norms 
also topologize $A$,
and we will use them in place of our original family,
so that we may work with a countable family of norms for $A$,
with zeroth norm equal to $\| \cdot \|_B$.

We say that $A$ is a {\it right Fr\'echet ideal} in $B$ 
if $a \in A$ and $b \in B$
implies $ab \in A$, and this multiplication operation 
is continuous for the respective topologies.  
In terms of norms, this means that for each 
$n\in \N$ there exists an 
integer $m \geq n$ and constant $C_n > 0$ such that
the inequality
\begin{equation}
\| ab \|_n \leq C_n \| a\|_m \| b \|_B
\label{eqn:idealcond}
\end{equation}
holds for all $a \in A$ and $b \in B$.
Similarly, a {\it left Fr\'echet ideal} 
and {\it two-sided Fr\'echet ideal} 
is defined.

If $\{ \| \cdot \|'_n \}_{n=0}^\infty$ is an equivalent increasing 
family of seminorms for $A$, then the inequality $(\ref{eqn:idealcond})$
is still satisfied but with adjusted constants $C'_n$ and integers
$m'_n$.  If $\| \cdot \|'_B$ is an equivalent Banach algebra
norm on $B$, then the constants $C_n$ will scale uniformly in $n$,
with each $m_n$ staying the same for a given $n$ 
in $(\ref{eqn:idealcond})$.

\vskip\baselineskip
\noindent{\bf Remark 3.2.  Spectral Invariance.}
A subalgebra $A$ of an algebra $B$ is {\it spectral invariant}
in $B$ if every element $a \in A$ is quasi invertible in $B$ if and
only if it is quasi invertible in $A$.
An element $x \in B$ is a  {\it quasi-inverse} for $y \in B$ if
$x \circ y = y \circ x = 0$, where
$x \circ y$ is defined as $x + y + xy$ for any $x, y \in B$.
If $A$ is a left or right ideal in an
algebra $B$, then $A$ is spectral invariant in $B$.
For let $x \in A$ have quasi-inverse $y \in B$.
If $A$ is a right ideal, then $xy \in A$.  So $0=x \circ y =
x + y + xy$ and $y = -x - xy \in A$.
For left ideals, apply the same argument 
with $yx$ in place of $xy$.

\vskip\baselineskip
We say that a scale $\sigma$ on $X$ 
is {\it proper} if 
the inverse map $\sigma^{-1}$ takes
bounded subsets of $[1, \infty)$ to finite subsets of $X$.

\noindent{\bf Example 3.3.}
Let $B=c_0(X)$ be the commutative $C^\star$-algebra of 
complex-valued sequences which vanish at infinity, 
with pointwise multiplication and
sup-norm $\| \cdot \|_B = \| \cdot \|_\infty$.

\noindent
(a) Let $A$ be the dense Banach subalgebra $\ell^p(X)$ 
for some $1 \leq p < \infty$, with pointwise multiplication.
The inequality $\| \varphi \psi \|_p \leq 
\|\varphi \|_p \| \psi \|_\infty$
is satisfied for all $\varphi, \psi \in A$, so $A$ is a dense 
Banach ideal in $B$.  $A$ is not nuclear by Proposition 2.1 (b).

\noindent
(b) Let $\sigma$ be a proper scale on
$X$, and
let $A$ be the Fr\'echet space of 
$\sigma$-rapidly vanishing 
sequences ${\cal S}_\sigma^\infty(X)$,
topologized by the sup-norms 
$\| \varphi \|_n = \| \sigma^n \varphi \|_\infty$
(see Definition 2.7).
The inequalities $\| \varphi \psi \|_n 
\leq \|\varphi \|_n \| \psi \|_\infty$
are satisfied, so $A$ is a dense Fr\'echet ideal
in $B$.  
$A$ is nuclear if and only if there exists a $p \in \N^+$ for which
$\sum_{x \in X} {1\over {\sigma(x)^p}} < \infty$ 
(Theorem 2.9).
For example, this sum is bounded with $p=2$
when $\sigma$ is an enumeration of $X$ (see Corollary 2.13).

\vskip\baselineskip
\noindent{\bf Proposition 3.4. Unital Banach Algebras.} 
{\it 
Let $A$ be a subalgebra of a Banach algebra $B$.  
If $A$ is dense in $B$, then a (left or right) unit for $A$ is also 
a (left or right) unit for $B$.

If $B$ has a two-sided unit or left (right) unit,
and $A$ is a dense right (left) ideal, then $A$ 
contains the same unit.   
If $B$ has a left (right) unit $1_B$,
and $A$ is a dense left (right) ideal,
then $A$ contains a left (right) unit, 
possibly different than $1_B$.

If $A$ and $B$ have the same left (right) unit,
and $A$ is a right (left) ideal in $B$,
then $A=B$.
If $A$ is a right (left) Fr\'echet ideal, the equality
$A=B$ is topological, and $A$ is a Banach algebra exactly 
equal to $B$.}

\vskip\baselineskip
\noindent
In other words, dense ideals can only be proper when the 
Banach algebra $B$ is non-unital.  

\vskip\baselineskip
\noindent{\bf Remark 3.5. One-sided Units Not Unique.}  
If $B$ has a left unit $1_L$, and $A$ is a dense left ideal,
then $A$ will contain a left unit, but maybe different
than $1_L$.  (A similar statement holds for right units.)
Let $B= \ell^2(\N)$ with multiplication 
$\chi * \eta = \chi_0 \eta$.  Then $1_L = (1, 0, \dots, 0, \dots)$
is a natural left unit to take for $B$.
Let $\xi = (1,\, 1/2,\, 1/3, \dots 1/k,\, \dots) \in \ell^2(\N)$,
and take $A = \C \xi \oplus (\xi^\perp \cap c_f(\N))$.
Note that $A$ is a left ideal in $B$, but not a right ideal, 
and $\xi \in A$ is a left unit for $B$, but $1_L \notin A$.
To see that $A$ is dense in $B$, let 
$\chi \in \xi^\perp$ and $\epsilon > 0$.
Let
$\eta \in c_f(\N)$ satisfy $\| \eta - \chi \|_2 < \epsilon$.
Then $\eta' = \eta - < \eta, \xi > \delta_0$ is 
in $\xi^\perp\cap c_f(\N)$,
and $\| \eta'  - \chi \|_2 \leq  
\epsilon + |<\eta, \xi>|
= \epsilon + |< \eta - \chi, \xi>| 
\leq  \epsilon + \|\eta -\chi \|_2 \| \xi \|_2 \leq 
 \epsilon + \epsilon \| \xi \|_2$, using the 
Cauchy-Schwartz inequality.
To make $A$ a nuclear left Fr\'echet ideal, 
replace $\xi^\perp \cap c_f(\N)$ with
$\xi^\perp \cap {\cal S}(\N)$.

\vskip\baselineskip
\noindent{Proof of Proposition 3.4:} 
Let $A$ be a dense subalgebra of $B$, and $1_L$ a left unit in $A$,
satisfying $1_L a = a$ for every $a \in A$. 
Let $b \in B$,
and $a_\epsilon \in A$ satisfy $\| a_\epsilon - b \|_B < \epsilon$.
Then $\|1_L b - b \|_B \leq 
\| 1_L b - 1_L a_\epsilon \|_B + \| 1_L a_\epsilon - b \|_B
\leq \| 1_L \|_B \epsilon + \epsilon$.  Letting $\epsilon \rightarrow 0$,
we see that $1_L b = b$ and $1_L$ is a left unit for $B$.
The same argument works for a right unit.

Let $1_B$ be a two-sided unit in $B$, and let $a \in A$ be close to $1_B$
in $\| \cdot \|_B$, so that $b = a^{-1}$ exists in $B$.  
If $A$ is a right ideal, then $1_B = ab \in A$, 
so $A$ is unital with the same unit as $B$. 
The same argument holds if $A$ is a left ideal.

Let $1_L\in B$ be a left unit for $B$.  
If $A$ is a right ideal, 
then $A1_L$ is a dense right ideal of the closed
subalgebra $B1_L$ of $B$ since
$ a 1_L \bigl( b1_L \bigr) = \bigl( a_1 1_L b  \bigr) 1_L
\in AB 1_L \subseteq A 1_L$.
If $A$ is a left ideal, 
then $A1_L$ is a dense left ideal 
of $B1_L$ since
$\bigl( b1_L \bigr) a1_L = \bigl( b 1_L a \bigr) 1_L
\in BA 1_L \subseteq A 1_L$.
Since $1_L$ is a two-sided
unit for $B1_L$, we know $1_L \in A1_L$
by the previous paragraph.
If $A$ is a right ideal, $1_L \in A1_L \subseteq A$.
If $A$ is a left ideal,
let $a_0 \in A$ be such that $1_L = a_0 1_L$.
Then $a_0$ is a left unit for $B$ contained in $A$.
A similar argument applies to right units.

Next assume that $A$ is a right
ideal and $1_L \in A$ is a left unit for $A$ and $B$.  
Then for every $b \in B$, $b = 1_L b \in A$,
so $A = B$.   If $A$ is a right Fr\'echet ideal, 
then
$\| a  \|_n  = \| 1_L a \|_n  \leq C_n \| 1_L \|_m \| a\|_B  = 
C'_n \| a \|_B$ for all $a \in A$,
so each seminorm for $A$ is bounded by $\| \cdot \|_B$.  
By the continuity of inclusion $A\hookrightarrow B$,
$A$ has a topology equivalent to $B$'s.
\qed

\vskip\baselineskip
\noindent{\bf Proposition 3.6. Quotients and Idempotent Subalgebras.}  
{\it Let $A$ be a Fr\'echet subalgebra of a Banach algebra $B$.  

\noindent
(a) If $J$ is a closed two-sided ideal of $B$, then
the image of $A$ in the quotient Banach algebra $B/J$
is a Fr\'echet subalgebra.  
The image of $A$ in the quotient $B/J$ is
a left (right) Fr\'echet ideal if $A$ is a 
left (right) Fr\'echet ideal in $B$.
The image is nuclear if $A$ is nuclear, 
and dense if $A$ is dense in $B$.

\noindent
(b) Let $e$ be an idempotent in $B$. 
Then $Be$ and $eB$ are closed subalgebras of $B$.
If $A$ is a left (right) Fr\'echet ideal in $B$,
then the subalgebra $Ae$ ($eA$) inherits a natural quotient 
Fr\'echet topology from $A$, and is a left (right)
Fr\'echet ideal in $B$.
Both $Ae$ and $eA$ are nuclear if $A$ is, and $Ae$ ($eA$) is
dense in $Be$ ($eB$) if $A$ is dense in $B$.
If $e \in A$, then 
$Ae$ and $eA$ are closed subalgebras of $A$.
}

\vskip\baselineskip
\noindent
Note that Remark 3.5 gives an example 
when $e \notin A$, and $Ae$ is not contained in $A$. 

\vskip\baselineskip
\noindent{Proof of (a):}
Let $\pi \colon B \longrightarrow B/J$ be the canonical
quotient map.  The Fr\'echet algebra $A$ maps into
$B/J$ by a composition of continuous 
maps $\pi \circ \iota \colon A \rightarrow B/J$, where
$\iota \colon A \hookrightarrow B$ is the
inclusion map.  So the kernel $I = J \cap A$ is a
closed two-sided ideal in $A$, and we identify 
$A/I$ with a Fr\'echet subalgebra of $B/J$.

Assume that $A$ is a right Fr\'echet ideal in $B$,
and let $a \in A$ and $b \in B$.  Then $ab \in A$ so
$\pi(a) \pi(b) \in A/I$.  Let $\{ \| \cdot \|_n \}_{n=0}^\infty$
be increasing seminorms for $A$, and $\| \cdot \|_B$ the
norm on $B$.  The $n$th quotient seminorm for $A/I$ is
\begin{displaymath}
\| \pi(a) \|_n = \inf_{i \in I} \| a + i \|_n,
\end{displaymath}
and similarly $\inf_{j \in J}\|b + j\|_B$ is
the norm on $B/J$.
Using the seminorm inequality (\ref{eqn:idealcond}) 
for right ideals,
\begin{eqnarray}
\| \pi(a) \pi(b) \|_n & = & \inf_{i \in I} \| ab + i \|_n \nonumber \\
& = & \inf_{i \in I} \inf_{j \in J} \| a(b + j) + i \|_n
\qquad \text{since $AJ \subseteq A\cap J = I$} \nonumber \\
& \leq & \inf_{i \in I} \inf_{j \in J} \| (a + i)(b + j)\|_n 
\qquad \text{since $IB \subseteq A\cap J = I$} \nonumber \\
& \leq  & \inf_{i \in I} \inf_{j \in J} 
C_n \| (a + i)\|_m \| (b + j) \|_B 
\qquad \text{by inequality (\ref{eqn:idealcond})} \nonumber \\
& = & C_n \| \pi(a)\|_m \| \pi(b) \|_B,
\end{eqnarray}
so $A/I$ is a right Fr\'echet ideal in $B/J$.
The \lq\lq left\rq\rq\  case is handled in the same way.

If $A$ is nuclear, then $A/I$ is nuclear,
since quotients by closed linear subspaces preserve 
nuclearity [Pietsch, 1972], Proposition 5.1.3,
or
[Treves, 1967], Proposition 50.1 (50.4).
If $a_k \in A$ is a sequence
approaching $b$ in $B$, then by continuity 
$\pi(a_k)$ approaches $\pi(b)$ in $B/J$.
So $A/I$ is dense in $B/J$ whenever $A$ is dense in $B$.

\vskip\baselineskip
\noindent{Proof of (b):}
Let $e$ be an idempotent in $B$, 
and consider the linear subspace 
$Be = \{ b e \,\, | \,\, b \in B \}$ of $B$.  
If $\{ b_k \}_{k \in \N}$ is a sequence in
the Banach algebra $B$, such that
$b_k e$ converges in $B$ to a limit $b_0$, then $b_0 e = b_0$,
since multiplication is continuous.  
So $Be$ is a closed linear subspace of $B$, 
and clearly a subalgebra.  Similarly for $eB$.  If $e \in A$,
the same argument shows $Ae$ and $eA$ are closed 
subalgebras of the Fr\'echet algebra $A$.

Assume that $A$ is a left Fr\'echet ideal in $B$.  Then $Ae$ 
is a left ideal in $B$ since
$b (a e) \in B (A e) \subseteq (BA) e \subseteq A e$, 
for $a \in A$ and $b\in B$.
Note that $I = \{ \, a \in A \, | \, ae=0 \, \}$ is a closed
left ideal in $A$, 
with quotient $A/I = Ae$.
The inherited topology on $Ae$ is given by the quotient seminorms
$\|  ae  \|'_n = \inf_{ i \in I} \| a + i \|_n$, for $ae \in Ae$.
For $b \in B$, $ae \in Ae$, we have
\begin{eqnarray}
\| b(ae) \|'_n & = &  \inf_{i \in I}\, \|  ba + i \|_n \nonumber \\
& \leq & \inf_{i \in I}\, \| b(a + i) \|_n \nonumber \\
& \leq & C_n\, \| b \|_B \,\, \inf_{i \in I} \| a + i \|_m \nonumber \\
& = & C_n\, \| b \|_B\, \| ae \|'_m,
\end{eqnarray}
so $Ae$ is a left Fr\'echet ideal in $B$.

If $A$ is nuclear, then $Ae$ is nuclear,
since quotients by closed linear subspaces preserve 
nuclearity [Pietsch, 1972], Proposition 5.1.3,
or
[Treves, 1967], Proposition 50.1 (50.4).
If $b e \in B$ and $a_k \rightarrow b$ in $B$,
then $a_k e \rightarrow be$ in $Be$, by continuity
of multiplication, so $Ae$ is dense in $Be$ if $A$ is
dense in $B$.
\qed

\vskip\baselineskip
\noindent{\bf Remark 3.7.}  Note that the quotient 
Fr\'echet ideal $A/I \subseteq B/J$ in Proposition 3.6 (a)
satisfies inequality (\ref{eqn:idealcond}) with the 
same integers $m_n$ and constants $C_n$ as the
original Fr\'echet ideal $A \subseteq B$.
The same is true for the ideals $Ae \subseteq B$
(or $eA \subseteq B$).

\vskip\baselineskip
\noindent{\bf Proposition 3.8. Algebraic Ideals.}  
{\it Let $A$ be a Fr\'echet subspace of a Banach algebra $B$,
with continuous inclusion $\iota \colon A \hookrightarrow B$.
If $A$ is a left (right) ideal in the purely algebraic sense,
then $A$ is a left (right) Fr\'echet ideal in $B$.}

\noindent{Proof:}
Let $L_b \colon A \rightarrow A$ denote left multiplication 
by some $b \in B$.
Assume $a_\alpha \rightarrow 0$ and $L_b(a_\alpha) 
\rightarrow a_0$ in $A$.
Then $a_\alpha \rightarrow 0$ and $L_b(a_\alpha) 
\rightarrow a_0$ in $B$, by continuous inclusion.
But since $B$ is a Banach algebra, the multiplication
is continuous, and $a_0=0$.  Apply the closed graph theorem to see that
$L_b$ is a continuous linear map from $A$ to $A$.

Let $R_a \colon B \rightarrow A$ denote right 
multiplication by some $a \in A$.
Let $b_\alpha \rightarrow 0$ in $B$ and 
$b_\alpha  a \rightarrow a_0$ in $A$.
Then $b_\alpha a \rightarrow a_0$ in $B$, by continuous inclusion,
and so $a_0=0$.  Apply the closed graph theorem again to see that $R_a$
is continuous from $B$ to $A$.

We have shown that the bilinear map of multiplication 
$M \colon B \times A \rightarrow A$
is separately continuous.
By [Rudin, 1973], Theorem 2.17, $M$ is jointly continuous.
\qed

\vskip\baselineskip
\vskip\baselineskip
\section{Dense Nuclear Ideals in $C^\star$-Algebras}

In this section and the next, we apply the dense ideal
condition together with nuclearity, to get 
results about the structure of the Banach algebra.
Here we obtain complete results for $C^\star$-algebras,
but wait until \S 5 for the general Banach
algebra case.

\vskip\baselineskip
An algebra is {\it semiprime} if it 
has no non-zero nilpotent ideal.  (One can use two-sided,
left, or right ideals to define semiprime; the resulting definitions
are equivalent [Palmer, 1994], Proposition 4.4.2 (d).)

\noindent{\bf Proposition 4.1.}  
{\it Let $A$ be a dense nuclear left (right) Fr\'echet 
ideal of a Banach algebra $B$.  
Let $e \in B$ be an idempotent.  
Then $Ae$ ($eA$)  
is a finite dimensional Banach algebra equal to $Be$ ($eB$).

Assume further that $B$ is semiprime.
Then $Be$, $eB$, and $BeB$ 
are all finite dimensional Banach algebras equal to $Ae$, $eA$,
and $AeA$ respectively. 
The Fr\'echet ideal $A$ can be either left or right for this to work.}

\vskip\baselineskip
\noindent{\bf Remark 4.2.}  
If $B$ is not semiprime, it may happen that $A$ is a dense nuclear
left Fr\'echet ideal,  but $eB$ and $BeB$ are not finite dimensional,
for some idempotent $e \in A$.
Let $B$ be the Hilbert space $\ell^2(X)$ 
and $A$ be Schwartz functions ${\cal S}_\gamma(X)$, with the natural
inclusion map $A \hookrightarrow B$.
For $\chi, \eta \in B$, define multiplication by 
$\chi * \eta = \chi_1 \eta $.
Then $B$ is a Banach algebra and $A$ is a Fr\'echet algebra
for this multiplication.  
$B$ is not semiprime
since the (two-sided) ideal $I = \{ \chi \in B \,\, | \,\, \chi_1 = 0 \}$ 
satisfies $I^2=0$.
The Fr\'echet algebra $A$ is a 
dense nuclear left Fr\'echet ideal
in $B$, but not a right ideal since $AB = B$.  Let 
$e = (1, 0, 0,  \dots) \in A$.  Note that $e^2 = e$,
$Be = \C e$, $eB = B$ and $BeB = B$.

\vskip\baselineskip
\noindent{Proof of Proposition 4.1:}
By Proposition 3.6 (b), $Ae$ is a dense nuclear left Fr\'echet 
ideal in $Be$.
Since $e$ is a right unit for $Be$,
Proposition 3.4 tells us $e \in Ae$.   So
by the last part of Proposition 3.4, 
$Ae$ is a Banach algebra exactly
equal to $Be$.  Being a nuclear Banach space, $Ae=Be$ is 
finite dimensional (Proposition 2.1 (b)).

Next, I will imitate the proof of the first Lemma of [Smyth, 1980], 
to show that $eB$ is also finite dimensional, with
the added assumption that $B$ is semiprime.
We just proved that $Be$ is finite dimensional, so $eBe$
must also be finite dimensional.   
Assume the dimension is some positive integer $N$,
and let $\theta \colon eBe \rightarrow \C^N$ be a linear
bijective map.
For each $y \in eB$, define a 
linear map $\varphi_y \colon Be \rightarrow \C^N$
by $\varphi_y(x) = \theta(yx)$, 
$x \in eB$.   The map  $y \mapsto \varphi_y$ from $eB$
to the (finite dimensional) space 
of linear maps ${\cal L}(Be, \C^N)$ is linear.
If $\varphi_{y_0} = 0$ for some $y_0 \in eB$, then
$ y_0 Be = 0$,  and $(By_0 )^2 = (B e y_0)^2$ (since $y_0 \in eB$)
$= Be (y_0 Be) y_0 =  0$.  So the left ideal $B y_0$
is nilpotent with order 2, 
contradicting our assumption that $B$ is semiprime.
Therefore the mapping $ y \in eB \mapsto   
\varphi_y \in {\cal L}(Be, \C^N)$ 
is one-to-one, and hence dim$(eB) <$
dim$({\cal L}(Be, \C^N)) = $ dim$(Be) \cdot N  < \infty$. 

Now we know that both $eB$ and $Be$ are finite dimensional.
Let $b_1, \dots b_K \in B$ and $c_1, \dots c_L \in B$
satisfy $eB = \C-$span$ \{ e b_i  \,\, | \,\, i = 1, \dots K\}$
and $Be = \C-$span$ \{ c_i e \,\, | \,\, i = 1, \dots L\}$.
Then $BeB$ has dimension at most $K  L $ since
$BeB = \C-$span$ \{ b_i e c_j \,\, | \,\, i = 1, 
\dots K, \quad j = 1, \dots L \}$.
\qed

\vskip\baselineskip
\noindent{\bf Theorem 4.3.}  {\it  Let $B$ be a commutative 
$C^\star$-algebra, with maximal ideal space $M$.
Then $B$ has a dense nuclear Fr\'echet ideal if and only if
$M$ is discrete and countable.}

\noindent{Proof:}
The maximal ideal space $M$ is 
locally compact, and $B$ is isomorphic to the 
$C^\star$-algebra $C_0(M)$ of continuous 
functions vanishing
at infinity on $M$ [Dixmier, 1982], \S 1.4.1. 
Let $m \in M$,
and let $U\subseteq M$ be an open, relatively compact
set about $m$.  The set of functions $f \in B$ which
vanish on  the closure $\overline U$, is a 
closed ideal $I_{\overline U}$ in $B$. 
Let $A$ be a dense nuclear Fr\'echet ideal in $B$.  
Applying Proposition 3.6 (a),
we find the image $\pi(A)$ of $A$ in the quotient 
$B/I_{\overline U}$
is again a dense nuclear Fr\'echet ideal.  
Since $B/I_{\overline U} \cong C({\overline U})$ is unital
by the compactness of ${\overline U}$,
Proposition 4.1 (with $e= 1_{C({\overline U})}$) 
tells us that $\pi(A)$ is equal to $B/I_{\overline U}$ 
and finite dimensional.
But $C({\overline U})$ can only be finite dimensional
if ${\overline U}$ is a finite set of  points.
Thus every point $m \in M$ has a finite neighborhood.
This proves that $M$ is discrete.  

To see that $M$ is countable, let
$\delta_m$ denote the unit step function at $m \in M$.
An element $b \in B$ can satisfy $\| b - \delta_m \|_B < 1/4$
for at most one $m \in M$, since
$\|\delta_{m_1}  -\delta_{m_2}\|_B= 1$
for distinct $m_1, m_2 \in M$.
By Corollary 2.3, $B$ has a countable dense set $S$.  
The correspondence
$m \mapsto$ \lq\lq choose $b \in S \text{ within distance 1/4 of } 
\delta_m$\rq\rq\ 
gives an injective map from $M$ into $S$, so $M$ is countable.

Since $M$ is discrete and countable, it must be either finite
or countably infinite.
For the infinite case, $B \cong c_0(X)$,
and $A={\cal S}_\gamma(X)$, Schwartz functions on $X$,
is a dense nuclear Fr\'echet ideal 
(see Equation (\ref{eqn:standSchw})
or Example 3.3 (b)).
\qed

\vskip\baselineskip
\noindent{\bf Lemma 4.4.  Existence of Projections.}  
{\it 
Let $B$ be a $C^\star$-algebra with a dense
nuclear left or right Fr\'echet ideal $A$.  
Let $I$ be any proper closed two-sided ideal in $B$.
Then $A$ contains a nontrivial projection which does
not lie in $I$.
}

\noindent{Proof:}
Let $ a$ be an element of $A - I$, with image
$[a]$ in the quotient $B/I$.  
Let
$\{ e_\lambda \}_{\lambda \in \Lambda}$ be an approximate
identity for the ideal $I$.  
Assume $A$ is a right Fr\'echet ideal,
and note that $\|a - a e_\lambda \|_B$
approaches $\|[a]\|_{B/I}$.  
By the Fr\'echet ideal condition,
$\| a e_\lambda \|_n \leq C_n \| a \|_m 
\| e_\lambda \|_B = C_n \| a \|_m$,  $n \in \N$, $\lambda \in \Lambda$.
Since a bounded set in a nuclear locally convex space is
relatively compact (Proposition 2.1 (d)),
there is a cluster point $i_0 \in I \cap A$ 
for the net 
$\{ a e_\lambda \}_{\lambda \in \Lambda}$,
converging in the Fr\'echet topology of $A$.
Replace $a$ with $a - i_0$ so now we have 
$\| a \|_B = \|[a]\|_{B/I} \not= 0$.
Replace $a$ with $a a^\star/\| a\|_B^2$. 
This new $a$ remains in $A$, is positive, and
satisfies $\| a \|_B = \|[a]\|_{B/I} = 1$.

For each $k \in \N$,  
$a^{2^k}$ has $C^\star$-norm equal to one, by applying the
$C^\star$-identity $\|a^{2^{k+1}}\|_B = \|a^{2^k}\|^2_B$
repeatedly.
By the Fr\'echet ideal condition,
$\| a^{2^k} \|_n \leq C_n \| a^{2^k -1} \|_B 
\| a \|_m \leq C_n \| a \|_m$,  $k, n \in \N$.
Since a bounded set in a nuclear locally convex space is
relatively compact (Proposition 2.1 (d)),
there is a subsequence $\{ a^{k_i}\}_{i=0}^\infty$,
converging in the Fr\'echet topology of $A$.
The limit point, $a_0 \in A$, must also have unit norm in $B$.
Moreover, the image $[a_0]$ of $a_0$ in $B/I$
also has norm $1$, since $\| [a^{k_i}] \|_{B/I} = 1$ for
eack $k_i$.

Consider the commutative $C^\star$-subalgebra $C^\star(a)$ of $B$
generated by $a$.  This algebra must be isomorphic to
continuous functions on a locally compact space $M$, vanishing at 
infinity.  We may think of $a$ and $a_0$ as real-valued
functions on $M$, with range in $[0, 1]$, both taking on the
value $1$ for at least one point of $M$.  Since $a^{k_i} 
\rightarrow a_0$ in the sup-norm, $a^{k_i}(m) \rightarrow a_0(m)$
for each $m \in M$.
It follows that $a_0(m) \in \{ 0, 1 \}$, and $a_0$ is
a projection. 
\qed

\vskip\baselineskip
\noindent{\bf Definition 4.5.  The Finite Socle.}  
Let $A$ be any algebra. 
The {\it left (right) finite socle of $A$}
is the sum of all minimal left (right) finite 
dimensional ideals of $A$.
If the left and right finite socles are equal,
their common value is the {\it finite socle of $A$}, denoted
by $A_{\rm{fin}}$.

For $a \in A$ and minimal left ideal $L$ of $A$, 
$La$ is a minimal left ideal 
or $\{ 0 \}$ [Palmer, 1994], Proposition 8.2.8.
If $L$ is finite dimensional, so is $La$.  Hence
the left finite socle is a right, 
and therefore two-sided, ideal of $A$.
Similarly the right finite socle is a two-sided ideal of $A$.

Let $S$ be the
set $\{\, e \, | \, e \textrm{ is a minimal idempotent in } A \, \}$,
and assume $A$ is semiprime.  Then
$\{\, Ae \, | \, e  \in S \, \}$
($\{\, eA \, | \, e  \in S \, \}$)
gives all the minimal left (right) ideals of $A$
[Palmer, 1994], Corollary 8.2.3.
By the first Lemma of [Smyth, 1980], 
we know that
dim$(Ae)$ is finite if and only if
dim$(eA)$ is finite.
%, and dim$(Ae)$=dim$(eA)$ when finite.
Let $S_{\rm{fin}}$ denote those
$e\in S$ for which $Ae$ (or $eA$) is finite dimensional.
Then the left (right) finite socle is
equal to $AS_{\rm{fin}}$
($S_{\rm{fin}}A$).  By the previous paragraph,
these are both two-sided ideals in $A$.
But then they must 
both equal $AS_{\rm{fin}}A$, so the
left and right finite socles agree, 
and $A_{\rm{fin}} =AS_{\rm{fin}} = S_{\rm{fin}}A = AS_{\rm{fin}}A$.
Note that $C^\star$-algebras are semiprime.

\vskip\baselineskip
\noindent{\bf Theorem 4.6.}  
{\it 
Let $B$ be a $C^\star$-algebra containing a dense 
nuclear left or right Fr\'echet ideal $A$.
Then the finite socle $B_{\rm{fin}}$ is dense in both $B$ 
and $A$, 
the primitive ideal space $\rm{Prim}(B)$ is discrete
and countable, and $B/I$ is finite dimensional
for any $I \in \rm{Prim}(B)$.  $B$ is the countable 
direct sum, or restricted product, of finite 
dimensional matrix algebras.
}

\noindent{Proof:}
We do the proof for a left ideal $A$.
The closure 
$I = {\overline B}_{\rm{fin}}^B$ 
is a two-sided
ideal in $B$.  Assume for a contradiction that
$I \not= B$.  
Apply Lemma 4.4 to
get a nontrivial projection $p \in A  - I$.
By Proposition 4.1, $Ap = Bp$ is finite dimensional.
If  $I \cap Bp \not= \{0\}$,
find only finitely
many $e_1, \dots e_k \in S_{\rm{fin}}$ such that
$e_i p \not= 0$.  
%For any other $e \in S_{\rm{fin}}$,
%$pe=ep=0$.  
A simple calculation shows $e' = p - e_1 - \dots - e_k$
is an idempotent in $B$ which is orthogonal to $I$.
Since $p \notin I$, $e' \not= 0$.  
By Proposition 4.1, $Ae'=Be'$ is finite dimensional.
This is a left ideal of $B$ whose 
intersection with $I$ is $\{0\}$.
Let $e_{\rm{min}}$ be a minimal idempotent contained therein.
Then $Be_{\rm{min}}$ is a minimal left finite dimensional 
ideal in $B$, which is not in $B_{\rm{fin}}$, a contradiction.

Let $e$ be any idempotent in $B$.
The second part of Proposition 4.1 tells us $eA=eB$. 
Since $e \in eB$, 
we know $e \in eA$.   But $eA \subseteq A$ since $A$ is a left ideal.
So $A$ contains any idempotent of $B$, 
$S_{\rm{fin}} \subset A$, and $B_{\rm{fin}} \subseteq A$.

The distance in $B$ between two 
idempotents $e, f \in S_{\rm{fin}}$ 
is at least one, since 
$\| e - f \|_B \| e \|_B \geq \| e - ef \|_B = \| e \|_B \geq 1$.
By Corollary 2.3, $B$ is separable, so $S_{\rm{fin}}$ is at most
countable.

For $e \in S_{\rm{fin}}$, $BeB$ is a finite dimensional
matrix algebra with dimension $($dim$Be)^2$.  
So as a $C^\star$-algebra, $B$ is the direct
sum of finite dimensional matrix algebras, and Prim$(B)$
is well-known to be discrete [Fell Dor, 1988], Proposition 5.21.
Also $B/I$ is one of the finite dimensional 
direct summands for each $I \in \rm{Prim}(B)$.

To see that $B_{\rm{fin}}$ is dense in $A$,
let $P_K = \sum_{k\leq K} 1_k$ be the sum of 
of the first $K$ identity matrices, from the matrix
algebras in the direct sum for $B$.  Then
$\{ P_K\}_{K=0}^\infty$ is a bounded approximate unit for $B$.
For any $a \in A$, the left ideal condition tells us
that the set $\{ P_K a \, | \, K \in \N\}$  is bounded
in $A$.  By Proposition 2.1 (d), some $a_0 \in A$
is a cluster point.  But $a_0$ is also a cluster point
in the topology of $B$, so since $P_K a \rightarrow a$ in
$B$, we must have $a_0=a$.
Since $P_K a \in B_{\rm{fin}}$, this shows that
$B_{\rm{fin}}$ is dense in $A$.
\qed

\vskip\baselineskip
There is a notion of nuclearity for $C^\star$-algebras 
[Kad Ring II, 1997], Chapter 11,
which is different from the notion of nuclearity 
for locally convex spaces.

\noindent{\bf Corollary 4.7.  Nuclearity of the $C^\star$-algebra.}  
{\it Let $B$ be a $C^\star$-algebra with a dense left 
or right Fr\'echet ideal $A$.  If $A$ is nuclear as a locally
convex space, then $B$ is nuclear as a $C^\star$-algebra.}

\noindent{Proof:}
By Theorem 4.6, every irreducible 
representation of $B$ is finite dimensional, and
$B$ is separable by Corollary 2.3.
Applying [Dixmier, 1982],  \S 9.1 Theorem (iv) $\Rightarrow$ (i), 
we find that $B$ is a Type I $C^\star$-algebra.
But every Type I $C^\star$-algebra
is nuclear in the $C^\star$-algebraic sense [Paterson, 1988], (1.31).
\qed

\vskip\baselineskip
\vskip\baselineskip
\section{Complete Continuity and the General Banach Case}

In this section we generalize
Theorem 4.6 for an arbitrary Banach algebra.

For any Banach algebra $B$, we give Prim$(B)$ 
the Jacobson topology [Dixmier, 1982], \S  3.1.1. 
Our main goal in this section
is to show that Prim$(B)$ has a discrete topology,
whenever $B$ has a dense nuclear Fr\'echet ideal.
The proof of Theorem 4.3 does not easily generalize
to the noncommutative case, and the proofs of
Theorem 4.6 and Lemma 4.4 rely on the theory of $C^\star$-algebras.
In this section, we find a different approach to the problem,
with the concept of complete continuity.

\vskip\baselineskip
If $E$ and $F$ are Banach spaces, a continuous linear 
map $T \colon E \rightarrow F$ is {\it completely continuous}  
if it maps weakly converging sequences in $E$
to norm converging sequences in $F$.  
A Banach algebra $B$
is said to be {\it left completely continuous} if for every $b \in B$, 
the left  multiplication 
operator $L_b \colon B \rightarrow B$, given by
$x \in B \mapsto bx$,
is a completely
continuous map from $B$ to itself. 
Similarly, $B$ is {\it right completely continuous} if 
right multiplication $R_b$ is
completely continuous for every $b \in B$, and $B$ is
{\it  completely continuous} if it is both left and
right completely continuous [Kaplansky, 1949].

\noindent{\bf Theorem 5.1.}  {\it Let $B$ be a Banach algebra
containing a dense nuclear left (right) Fr\'echet ideal $A$.
Then $B$ is right (left) completely continuous.}

\noindent{Proof:}  
We do the proof for right Fr\'echet ideals.
Let $\{ x_k \}_{k=0}^\infty $ be a sequence in $B$, and $a$  
a fixed element of  $A$.
It is not hard to show that if
$x_k $ converges to a limit $x$ in the norm of $B$, then
$a x_k$  converges to $a x$ in the Fr\'echet topology of $A$, using
the right ideal inequality (\ref{eqn:idealcond}).
But we must begin by assuming $x_k$ converges weakly,
{\it not} in norm.

Let $\varphi$ be a continuous linear functional on
$A$.  
Then  $\varphi \circ L_a$ is  continuous
on $B$,  since
\begin{eqnarray}
| \varphi \circ L_a (x)| = |\varphi(a x) | & \leq &    C \| a x \|_n  \qquad  \text{$C, n$ exist since $\varphi$ is continuous on $A$} \nonumber \\
& \leq &  C C_n \| a \|_m \| x \|_B  \qquad  \text{by the ideal inequality (\ref{eqn:idealcond}),} 
\label{eqn:leftweakbound}
\end{eqnarray}
for any $ x \in B$.
Thus if our original 
sequence $\{ x_k \}_{k=0}^\infty $ converges weakly to
zero in $B$, then $a  x_k$ converges weakly to zero in $A$.  
To simplify notation, define $y_k = ax_k$.  
We wish to prove that $y_k \rightarrow 0$ in
the norm of $B$, by using the 
fact that $y_k \rightarrow 0$ weakly in $A$.

Let $U$ be the unit ball of $B$.  Then $U\cap A$ is absolutely convex, 
and by continuous inclusion $A \hookrightarrow B$, 
$U\cap A$ is a zero neighborhood of $A$. 
Apply nuclearity 
[Pietsch, 1972], Proposition 4.1.4,
to find an absolutely convex zero neighborhood
$V$ of $A$, and sequence of continuous 
linear functionals $\varphi_n \in A'$
satisfying 
\begin{equation}
\sum_{n=0}^\infty \| \varphi_n \|_{V^0} < \infty, 
\label{eqn:boundedvarphi}
\end{equation}
such that the inequality 
$
\| x \|_B \leq \sum_{n=0}^\infty | \varphi_n(x)|
$
holds for all $x \in A$.
Taking $x = y_k$, we get
\begin{equation}
\| y_k \|_B \leq \sum_{n=0}^\infty | \varphi_n (y_k) |
\label{eqn:rhsgoal}
\end{equation}
holds for every $k \in \N$.

\vskip\baselineskip

\noindent{\bf Lemma 5.2.}  {\it  The sequence $\{ y_k \}_{k=0}^\infty$
is bounded in the seminorm $\| \cdot \|_V$ on $A$.}

\noindent{Proof:}  
Since $y_k \rightarrow 0 $ weakly in the Fr\'echet space $A$,
the set $\{ y_k \}_{k=0}^\infty$ is weakly bounded in $A$.
In the terminology of [Rudin, 1973], Theorem 3.18, $V$ is an
original neighborhood of $A$.  So $\{ y_k \}_{k=0}^\infty$ 
is contained in $tV$ for some $t>0$.
\qed

\vskip\baselineskip

We use Lemma 5.2 to show that the right hand 
side of inequality (\ref{eqn:rhsgoal})
tends to zero as $k \rightarrow \infty$.  
Let $\epsilon > 0$.   
Let $N$ be large enough so that 
$\sum_{n=N+1}^\infty \| \varphi_n \|_{V^0} < \epsilon /2M$, where
$M$ is the bound on $\| y_k\|_V$ from Lemma 5.2.  This is possible
since the full series in (\ref{eqn:boundedvarphi}) converges.
Using weak convergence of $\{y_k\}_{k=0}^\infty$ to zero, find $K$
large enough so that 
$\sum_{n=0}^N | \varphi_n(y_k) | < \epsilon /2$ for $k \geq K$.  
So we have
\begin{eqnarray}
\sum_{n=0}^\infty |\varphi_n(y_k) | & \leq & 
\sum_{n=0}^N | \varphi_n (y_k) |
\quad + \quad
\sum_{n=N+1}^\infty| \varphi_n (y_k) |
 \nonumber \\
 &\leq & \epsilon/2 \qquad + \qquad  
\sum_{n=N+1}^\infty 
| \varphi_n (y_k) | 
\quad \qquad \text{since $k \geq K$} \nonumber \\
&\leq &  \epsilon/2\qquad  +  \qquad
\sum_{n=N+1}^\infty 
\| \varphi \|_{V^0} M
\quad \qquad \text{since $\| y_k \|_V \leq M$ by Lemma 5.2} \nonumber \\
& \leq & \epsilon/2 \qquad + \qquad \epsilon/2 \qquad = \qquad \epsilon. 
\end{eqnarray}
Thus the right hand side of inequality (\ref{eqn:rhsgoal}) is less than $\epsilon$
for $k \geq K$, so the left hand side $\| y_k\|_B$ is also less 
than $\epsilon$ for $k \geq K$.  Thus $y_k \rightarrow 0$
in the norm of $B$.   

We have proved that if $\{ x_k \}_{k=0}^\infty$
is a sequence in $B$ converging weakly to zero, and $a \in A$,
then $y_k = a x_k$ converges to zero in the norm of $B$.
To finish Theorem 5.1 and prove the complete 
continuity of $B$, we would like to replace $a \in A$ with 
an arbitrary element $b$ of $B$.
Since $x_k \rightarrow 0 $ weakly in $B$,
the Uniform Boundedness Principle (Theorem 2.4), 
%$\{x_k\}_{k=0}^\infty$ is a bounded sequence in $B$, so let
gives us an $M  < \infty$ 
%satisfy 
for which $\| x_k \|_B \leq M$ for all $k \in \N$.
Let $\epsilon > 0 $   and 
%using the density of $A$ in $B$, 
pick
$a \in A$ such that $\| a - b \|_B <\epsilon/2M$. 
Let $K$ be big enough so that $\|a x_k \|_B < \epsilon /2$
for $k \geq K$.  
Then
\begin{eqnarray}
\| b x_k \|_B  & \leq & 
\| (a - b) x_k \|_B + \| a x_k \|_B
 \nonumber \\
 &\leq & \| a - b \|_B \| x_k \|_B  \quad + \quad  
\epsilon /2 
\quad \qquad \text{since $k \geq K$} \\
 &\leq & (\epsilon/2M) M  \qquad\quad + \quad  
\epsilon /2 \qquad  = \qquad \epsilon,
\nonumber
\end{eqnarray}
so $bx_k \rightarrow 0$ in the norm of $B$.
This completes the proof of Theorem 5.1.
\qed

\vskip\baselineskip
A {\it primitive ideal} of a Banach algebra $B$
is the kernel of a non-zero algebraically
irreducible, continuous left Banach-space 
representation of $B$.  (This is equivalent to
the purely algebraic definition 
of primitive ideal [Palmer, 1994], Corollary 4.2.9.)
The {\it primitive ideal space} of $B$, or  Prim($B$), 
is the set of all primitive ideals in $B$.  We
also call Prim$(B)$ the {\it spectrum} of $B$.

\noindent{\bf Corollary 5.3.}  {\it Let $B$ be a Banach algebra
containing a dense nuclear two-sided Fr\'echet ideal $A$.
Then $B$ is completely continuous and
the primitive ideal space of $B$ is discrete.}

\noindent{Proof:}  
By Theorem 5.1 just proved, 
containment of a dense nuclear two-sided Fr\'echet 
ideal implies 
$B$ is completely continuous on both sides.
By [Kaplansky, 1949], Theorem 5.1, a completely 
continuous Banach algebra
has discrete primitive ideal space. 
\qed

\vskip\baselineskip
\noindent{\bf Proposition 5.4.}  {\it  Let $B$ be a Banach
algebra with a dense nuclear left or right Fr\'echet ideal.
Then $B/I$ is finite dimensional 
for each $I\in$\rm{Prim}$(B)$.}

\noindent{Proof:}
Let $V$ be an algebraically irreducible
left Banach $B$-module.
Then $I = \{ b \in B \,\, | \,\, bV=0 \} \in $Prim$(B)$ is the 
primitive ideal corresponding to $V$.
Let $v$ be a nonzero element of $V$, and define
$M = \{ b \in B \,\, | \,\, bv=0 \}$.  Then $M$
is a maximal modular closed left ideal of $B$,
and $V= B/M$ as left $B$-modules.

Let $A$ be a dense nuclear left Fr\'echet ideal of $B$.
Then $Av$ is dense in $V$ 
(since $Bv = V$ and $A$ is dense in $B$), so there must be
an $a_0 \in A$ such that $a_0 v \not= 0$.
Since $A$ is a left ideal in $B$, 
we have $Av \supseteq (Ba_0)v = B(a_0v) = V$.
Let $N = M \cap A = \{ a \in A \,\,|\,\, av=0\}$.  
Then $N$ is a closed left ideal in $A$,  $V=A/N$ as 
left $B$-modules, and
$A/N$ is nuclear 
[Pietsch, 1972], Proposition 5.1.3,
or
[Treves, 1967], Proposition 50.1 (50.4).

The map $a \in A \mapsto av \in V$ is continuous,
so the Fr\'echet topology on $A/N$ is at least as strong
as the original Banach space topology on $V$.\footnote{The same is
true if we give $V$ the Banach space topology of $B/M$.}
By [Treves, 1967], Chapter 17-7, Corollary 2,
the topologies must agree.
So $V$ is a nuclear Banach space and therefore finite 
dimensional (Proposition 2.1 (b)).   Since $B/I$ is
represented faithfully on $V$, it must also be 
finite dimensional.
\qed

\vskip\baselineskip
\noindent{\bf Remark 5.5.}  
An alternative proof of Proposition 5.4 
(assuming a dense nuclear two-sided ideal)
is given by applying the 
complete continuity of $B$ from Theorem 5.1, and
[Kaplansky, 1948], Lemma 4.

\vskip\baselineskip
The {\it Jacobson radical ${\cal A}_J$} of an algebra $\cal A$ is
the intersection of all primitive ideals of $\cal A$. 
$\cal A$ is {\it semisimple} if ${\cal A}_J=0$, and {\it radical} if
${\cal A}_J=\cal A$ [Palmer, 1994], \S 4.3.1.
A semisimple algebra is also 
semiprime [Palmer, 1994], Theorems 4.4.6 and 4.5.9.

\noindent{\bf Theorem 5.6.}  {\it 
Let $B$ be a Banach algebra
containing a dense nuclear two-sided Fr\'echet ideal $A$.
Then {\rm Prim}$(B)$ is discrete and countable, and
$B/I$ is finite dimensional for each $I \in \rm{Prim}(B)$.}

\noindent{Proof:}  Proposition 3.6 (a) shows the quotient of $B$ by 
its Jacobson radical still 
has a dense nuclear Fr\'echet ideal.  Since the Jacobson radical $B_J$
is the intersection of all kernels of primitive ideals, Prim$(B)$ = 
Prim$(B/B_J)$.  So without loss of generality, I will assume 
that $B$ is semisimple.

I will construct an idempotent for each primitive ideal.  
Let $I \in $ Prim$(B)$ and let $\{I\}^c$ denote the complement of 
the singleton $\{I\}$, namely $\{ J \in $Prim$(B)| J \not= I \}$.  
Discreteness (Corollary 5.3) tells us the intersection $L$
of all elements of $\{ I \}^c$ is not contained in $I$.  
(The closure of  
$\{ I \}^c$ in the Jacobson topology is by definition 
the set of primitive ideals
containing $L$.  By discreteness this closure must 
only be $\{ I \}^c$ and nothing more.)
Note that $L$, being the intersection of closed two-sided ideals, 
is itself a closed two-sided ideal in $B$.
Intersecting $L$ with $I$ gives the zero ideal,
since $L \cap I = \bigcap\{ J | J \in $Prim$(B)\}$
and $B$ is semisimple.
So $LI=IL \subseteq L \cap I=0$, and $B$ is the
direct sum of ideals $B = L \oplus I$. 
Since Prim$(B)$ is discrete, the singleton $\{I\}$ is a closed
set, so $I$ is not contained in any other primitive ideal.
This implies $L$ must be simple.

We need a unit for $L$.  By
Proposition 5.4, $L$ is finite dimensional, and $L$ is
semiprime since $B$ is [Palmer, 1994], Proposition 4.4.2 (e).
By the Wedderburn Theorem [Palmer, 1994], Theorem 8.1.1,
a finite dimensional semiprime algebra which is simple
is isomorphic to a full matrix algebra, and therefore unital.
Let $e_I$ be the unit of $L$. 

Since $LI=IL=0$, we know $e_I I= I e_I = 0$.  
Let $J \in $Prim$(B)$, $J \not= I$, with corresponding
idempotent $e_J$, satisfying $e_J J = J e_J = 0$.  
Since the singleton $\{ J \}$ is closed
in Prim$(B)$, $J$ cannot be contained in $I$, so $L \cap J \not=0$.
Since $L \cap J$ is a two-sided ideal, and $L$ is simple, we
have $L \subseteq J$, so $e_I \in J$.  
It follows that $e_I e_J = e_J e_I = 0$.
The distance between these two orthogonal idempotents is at least one, because
$\| e_I - e_J \|_B   \| e_I \|_B \geq \| e_I - e_J e_I \|_B = \| e_I \|_B$.

An element $b \in B$ can satisfy $\| b - e_I \|_B < 1/4$
for at most one $I \in$Prim$(B)$, since
$\|e_I  -e_J\|_B= 1$
for distinct $I, J \in$Prim$(B)$.
By Corollary 2.3, $B$ has a countable dense set $S$.  
The correspondence
$I \mapsto$ \lq\lq choose $b \in S \text{ within distance 1/4 of } 
e_I$\rq\rq\ 
gives an injective map from Prim$(B)$ into $S$, so Prim$(B)$ is countable.

The last statement follows from Proposition 5.4.
\qed

\vskip\baselineskip
\vskip\baselineskip

\section{Construction of Dense Nuclear Ideals 
for $C^\star$-Algebras}

Assume that a $C^\star$-algebra $B$ 
is a countably infinite direct sum of 
finite dimensional matrix algebras.
Let $\frak p = \{ p_z \}_{z \in Z}$ be
the dimensions, with $Z$ a countably
infinite set,
so that 
$B = \bigoplus_{z \in Z} M_{p_z}(\C)$.
Note that $\frak p$ is a scale on $Z$,
and repeated values of the same 
dimension are allowed.
We think of elements of $B$ as 
matrix-valued functions $f$ on $Z$, where 
$f(z) \in M_{p_z}(\C)$ for each
$z\in Z$.
The $C^\star$-norm on $B$ is
\begin{equation}
\| f \|_B = \sup_{z \in Z}
\| f(z) \|_{\text{op}},
\label{eqn:beesnorm}
\end{equation}
and $B$ consists of those functions $f$ which vanish at $\infty$
[Dixmier, 1982], \S 1.9.14.  
In this section, we construct dense nuclear ideals in $B$.

Let $Y$ be the disjoint union
of finite sets 
\begin{equation}
Y\quad  =\quad \bigcup_{z\in Z} \{z \}\times \{1, \dots p_z \} 
\times \{1, \dots p_z \},
\label{eqn:Ydef}
\end{equation}
and
let $e_y = e_{z, ij}$ be matrix elements for 
the $p_z\times p_z$ matrices $M_{p_z}(\C)$, for each
tuple $y = \{z, i, j\} \in Y$.
These are partial isometries which form a Schauder basis
for the $C^\star$-algebra $B$ (Definition 2.5).
Any $b \in B$ has
a coordinate functional
$b_y = < e_{z, ii} b e_{z, jj}, e_y> \in \C$,
and unique series expansion 
$b = \sum_{y \in Y} b_y e_y$ 
which converges unconditionally 
in $B$. 

Let $c_f(Y)$ be the linear span of the
$e_y$'s.  The finite socle $B_{\textrm{fin}}$
of $B$ (Definition 4.5) is identified with $c_f(Y)$, and
equals the algebraic direct sum 
of the matrix algebras $\bigoplus_{z \in Z} M_{p_z}(\C)$.

\vskip\baselineskip
\noindent{\bf Definition 6.1. Socle-Specific Schwartz Spaces.}
Let $\underline \ell$ be any family of scales on $Z$, with countably
infinite sets $Y$ and $Z$ defined above.
The Fr\'echet space
${\cal S}^{\infty, \textrm{op}}_{\underline \ell} (Y)$ is defined to be 
the completion of
$c_f(Y)$ in the norms
\begin{equation}
\| \varphi \|_n^{\infty, \textrm{op}} = 
\sup_{z \in Z} \ell_n(z) \|\varphi(z) \|_{\textrm{op}},
\label{eqn:rapidNormsOp}
\end{equation}
where $\| \cdot \|_{\textrm{op}}$ is the operator norm on $M_{p_z}(\C)$.

\vskip\baselineskip
\noindent{\bf Theorem 6.2. Existence of Dense Nuclear Fr\'echet Ideals.} 
{\it 
Let a $C^\star$-algebra $B$ be the
countably infinite direct sum of finite dimensional
$C^\star$-algebras, with
dimensions $\frak p = \{ p_z \}_{z \in Z}$.
If $\underline \ell$
is any family of scales on $Z$, 
then  } 
${\cal S}^{\infty, \textrm{op}}_{\underline \ell} (Y)$
{\it is a dense two-sided 
Fr\'echet ideal in $B$, in which 
$\{ e_y \}_{y \in Y}$ is 
an equicontinuous, unconditional basis.

The Fr\'echet ideal}
${\cal S}^{\infty, \textrm{op}}_{\underline \ell} (Y)$
{\it is nuclear if and only if
$\underline \ell$ satisfies the $\frak p$-summability condition
\begin{equation} 
(\forall n \in \N)\,(\exists m > n) \quad
\sum_{z \in Z}  \,
{p_z^2 \ell_n(z) \over \ell_m(z)}\, <\, \infty,
\label{eqn:NucSum}
\end{equation}
and if and only if}
${\cal S}^{\infty, \textrm{op}}_{\underline \ell} (Y) 
\cong {\cal S}^1_{\underline \sigma} (Y) 
\cong {\cal S}^\infty_{\underline \sigma}(Y)$, 
{\it where $\underline \sigma$ is the family of scales on $Y$ defined by 
$\sigma_n(z, i, j) = \ell_n(z)$
for $\{z, i, j\} \in Y$, $n \in \N$.}

\vskip\baselineskip
\noindent{\bf Remark 6.3.  Existence of Dense Nuclear Fr\'echet Ideals.}
It follows from Theorem 6.2 that a dense 
nuclear Fr\'echet ideal always exists, 
for any $C^\star$-algebra $B$ satisfying the hypotheses of Theorem 6.2,
since if $\zeta$ is any enumeration of $Z$,
then $\ell  = \zeta \frak p$ is a scale whose associated family of 
scales satisfies (\ref{eqn:NucSum}).
Since the ideal produced by this construction is
isomorphic to the Fr\'echet space ${\cal S}^1_\sigma(Y)$, 
the functions which vanish rapidly with 
respect to powers of the single scale $\sigma$,
it has the additional property of being a 
power series space of infinite type, as described in Remark 2.8.

\vskip\baselineskip
\noindent{Proof of Theorem 6.2:}
Since ${\underline \ell}  \geq 1$, 
$\| \cdot \|^{\infty,\textrm{op}}_0 \geq \| \cdot \|_B$, 
and the inclusion 
map ${\cal S}^{\infty, \textrm{op}}_{\underline \ell} (Y) \hookrightarrow B$
is continuous.
For $f \in B$ and $\varphi \in 
{\cal S}^{\infty, \textrm{op}}_{\underline \sigma}(Y)$,
we have
\begin{eqnarray}
\| f \varphi  \|^{\infty, \textrm{op}}_n 
&=& \sup_{z \in Z}\, 
\ell_n(z)  \,
\bigl\|f(z) * \varphi(z) \bigr\|_{\textrm{op}}
\qquad \qquad
\text{definition of $\| \cdot \|^{\infty, \textrm{op}}_n$} \nonumber\\
 &\leq& 
\biggl( \sup_{z\in Z} \| f(z) \|_{\text{op}}\biggr)
\,\, 
\biggl( \sup_{z\in Z} \ell_n(z) \| \varphi(z)\|_{\textrm{op}} \biggr)
 \nonumber\\
 & = & 
 \| f \|_B 
\| \varphi\|^{\infty, \textrm{op}}_n,
\end{eqnarray}
so ${\cal S}^{\infty, \textrm{op}}_{\underline\ell}(Y)$ 
is a left Fr\'echet ideal in $B$.  
Similarly it is 
a right, and therefore two-sided, Fr\'echet ideal in $B$.

The basis is equicontinuous since  for
$\varphi \in {\cal S}^{\infty, \textrm{op}}_{\underline\ell}(Y)$, \ 
$y = (z, i, j) \in Y$,  and $n \in \N$,
$$| \varphi(z)_{ij} | \| e_y \|^{\infty, \textrm{op}}_n
= \|  \varphi(z)_{ij} e_y \|^{\infty, \textrm{op}}_n
= \ell_n(z) \| e_{z, ii} \varphi(z) e_{z, jj} \|_{\textrm{op}}
\leq \ell_n(z) \| \varphi(z) \|_{\textrm{op}}  \leq
\| \varphi \|_n^{\infty,\textrm{op}}.$$ 
Use the fact that $\| \varphi(z) \|_{\text{op}}$
becomes smaller if matrix entries are set to zero,
to prove the basis is unconditional.

Since the operator norm on $M_{p_z}(\C)$ is bounded
by the sum of the matrix entries, and is greater
than or equal to any single matrix entry,
we have $\| \cdot \|^\infty_n \leq 
\| \cdot \|^{\infty, \textrm{op}}_n
\leq \| \cdot \|^1_n$.
This proves the continuity of inclusion maps
${\cal S}^1_{\underline\sigma} (Y)  \hookrightarrow
{\cal S}^{\infty, \textrm{op}}_{\underline\ell} (Y)  \hookrightarrow
{\cal S}^\infty_{\underline\sigma}(Y)$.
Since $\sigma$ is constant along each matrix algebra, 
the $\frak p$-summability condition (\ref{eqn:NucSum}) is
equivalent to the summability condition (\ref{eqn:NucSumSpace}), so
Theorem 2.9 tells us the three spaces are
isomorphic and nuclear if the summability condition is satisfied.

Conversely, assume 
${\cal S}^{\infty, \textrm{op}}_{\underline \ell} (Y)$ is a nuclear
Fr\'echet space.
We need to find a 
sequence space, as defined in Definition 2.7, which is nuclear,
to deduce the $\frak p$-summability 
condition (\ref{eqn:NucSum}) for $\underline \ell$.
The subspace of diagonal matrices is an obvious candidate.
Let $Y_D$ be the countable disjoint union
of  diagonal sets  $Y_D = \bigcup_{z \in Z} \{z \}
\times \{\{1, 1\}, \dots \{p_z, p_z\} \}$.
Define a family of scales $\underline \tau$ on 
$Y_D$ by $\tau_n(z, \{i, i\}) = \ell_n(z)$, and 
embed $\theta \colon {\cal S}^\infty_{\underline \tau} (Y_D) \hookrightarrow 
{\cal S}^{\infty, \textrm{op}}_{\underline \ell} (Y)$, via the
diagonal map 
$\theta(\varphi )(y) = \varphi(z, \{i, i\})  \delta(i-j)$.
For $\varphi \in {\cal S}^\infty_{\underline\tau} (Y_D)$,
$$\| \theta(\varphi) \|_n^{\infty, \textrm{op}} =
\sup_{z \in Z} \ell_n(z) \| \theta(\varphi)(z) \|_{\textrm{op}}
= \sup_{z \in Z}  \biggl\{ \ell_n(z) 
\sup_{i\leq p_z} \bigl| \varphi(z,\{i, i\}) \bigr| \biggr\} 
= \| \varphi \|^\infty_n,$$ 
since each $\theta(\varphi)(z)$ is a diagonal matrix.
So $\theta$ is isometric in all the norms.  Since
a subspace of a nuclear Fr\'echet space, with inherited 
topology, is nuclear
[Pietsch, 1972], Proposition 5.1.5, or
[Treves, 1967], Proposition 50.1 (50.3), 
${\cal S}^\infty_{\underline \tau} (Y_D)$ is nuclear.

Notice that for each $z \in Z$ and $i\leq p_z$, 
the linear functional
$x_{z,i} = \ell_n(z) e'_{zi}$
on ${\cal S}^\infty_\tau (Y_D)$ is continuous, since
$| x_{z,i} (\varphi) | = \ell_n(z) | \varphi(z, i) | 
\leq \| \varphi \|^\infty_n$.   
For each $n \in \N$, the countable set 
$S_n = \{ x_{z,i} \}_{z \in Z, i \leq p_z}$ 
converges weakly to zero in 
${\cal S}^\infty_\tau (Y_D)'$.
{\it Proof:}  For $\varphi \in c_f(Y_D)$, only finitely
many $x_{z,i} \in S_n$ have $x(\varphi)$ nonzero,
those whose pair
${z,i}$ lie in the support of $\varphi$. 
Since by definition ${\cal S}^\infty_\tau (Y_D)$ is 
the completion of $c_f(Y_D)$, for any 
 $\varphi \in {\cal S}^\infty_\tau (Y_D)$ 
find $\varphi_\epsilon \in c_f (Y_D)$
 $\epsilon$-close to $\varphi$.  Then
\begin{eqnarray}
\bigl| x_{z,i} (\varphi) \bigr| &\leq& 
\bigl| x_{z,i}\bigl( \varphi - \varphi_\epsilon \bigr) \bigr|
\quad + \quad 
\bigl| x_{z,i}\bigl( \varphi_\epsilon \bigr) \bigr|
\nonumber \\
& \leq &
 \bigl\| \varphi - \varphi_\epsilon  \bigr\|^\infty_n
\quad + \quad 
\bigl| x_{z,i}\bigl( \varphi_\epsilon \bigr) \bigr|
\nonumber \\
& \leq & \epsilon
\quad + \quad 
\bigl| x_{z,i}\bigl( \varphi_\epsilon \bigr) \bigr|.
\end{eqnarray}
For $z,i$ outside the finite set supp$(\varphi_\epsilon)$,
this shows $|x_{z,i}(\varphi)|  \leq \epsilon$. 

The set $S_n$ is an {\it essential}
subset of the polar of the unit ball for $\| \cdot \|^\infty_n$,
in the sense of [Pietsch, 1972], \S 2.3.1.   
By nuclearity and [Pietsch, 1972], Theorem 2.3.3, 
for any $n \in \N$, there exists an $m \in \N$ and a summing sequence
$c_{zi}$  of positive numbers so that for all $\varphi \in 
{\cal S}^\infty_\tau (Y_D)$, 
$$ \| \varphi \|^\infty_n \leq 
\sum_{z \in Z, i \leq p_z} c_{zi} | \varphi(z, i) | \ell_m(z).$$
For each $z_0$ and $i_0$, 
plug in $\varphi(z, i) = \delta(z-z_0, i-i_0)/\ell_m(z)$.
The result is $\ell_n(z_0) / \ell_m(z_0) \leq c_{z_0i_0}$.
Hence we have
$$\sum_{z \in Z}\, {p_z \ell_n(z) \over \ell_m(z)}
 \, = \, \sum_{z \in Z, i \leq p_z}\, {\ell_n(z) \over \ell_m(z)}
\,  \leq  \,
\sum_{z\in Z, i \leq p_z} c_{zi} \, < \infty.$$
Repeat the same argument to find a $p \in \N$ such that
$$\sum_{z \in Z}\, {p_z \ell_m(z) \over \ell_p(z)}\, < \infty.$$
Let $C_1, C_2>0$ be constants bounding these respective sums.
Then
\begin{eqnarray}
\sum_{z \in Z}\, {p_z^2 \ell_n(z) \over \ell_p(z)} & = &
\sum_{z \in Z}\, {p_z \ell_n(z) \over \ell_m(z)} 
\, {p_z \ell_m(z) \over \ell_p(z)} 
 \nonumber \\
& \leq & 
C_1 \,\,
\sup_{z \in Z} {p_z \ell_m(z) \over \ell_p(z)} 
\qquad < \qquad C_1 C_2, 
\end{eqnarray}
which is the $\frak p$-summability condition (\ref{eqn:NucSum}).
\qed

\vskip\baselineskip
We show that every dense nuclear two-sided Fr\'echet ideal
is given by Theorem 6.2. 

\noindent{\bf Theorem 6.4. Classification of 
Dense Nuclear Two-Sided Ideals.} {\it 
Let a $C^\star$-algebra $B$ be the
countable infinite direct sum of finite dimensional
$C^\star$-algebras, with sequence of
dimensions $\frak p = \{ p_z \}_{z \in Z}$.
Let $A$ be any dense nuclear two-sided Fr\'echet
ideal in $B$.
Then there exists a family of scales 
$\underline \ell$
on $Z$ satisfying the $\frak p$-summability condition}
(\ref{eqn:NucSum}) {\it such that the map
$a \mapsto \{ y \mapsto a_y \}$
gives an isomorphism of Fr\'echet ideals }
$A \cong {\cal S}^{\infty, \textrm{op}}_{\underline \ell} (Y)$.

\noindent{Proof:}
First we show that 
$\{ e_y \}_{x \in Y}$ is an absolute basis
for $A$. 
Since $A \subset B$, every element $a$ of $A$ also
has a unique expansion in
$\{ e_y \}_{y \in Y}$.
We will show the
series for $a$ converges in the Fr\'echet topology.
Pick some enumeration $\zeta$ of $Z$.
For $K\in \N^+$, let $P_K = \sum_{\zeta(z) \leq K} 1_z$ be the
sum of the units of the first $K$ matrix algebras
which make up $B$.  
By Theorem 4.6, the socle $c_f(Y)$ is dense in $A$.
Let $\varphi \in c_f(Y)$ be $\epsilon$
close to $a$ in $\| \cdot \|_m$.
For large enough $K$, $P_K \varphi = \varphi$, and we have
\begin{eqnarray}
\biggl\| \sum_{\zeta(z) \leq K} a_x e_x - a \biggr\|_n 
& \leq  &
 \biggl\| \sum_{\zeta(z) \leq K} a_x e_x - 
\varphi  \biggr\|_n
+ \| \varphi - a \|_n \nonumber  \\
& = & 
 \| P_K (a  - 
\varphi ) \|_n
+ \| \varphi - a \|_n \nonumber  \\
& \leq & 
 C_n \|a  - 
\varphi  \|_m
+ \| \varphi - a \|_n \quad  \leq \quad C_n \epsilon + \epsilon, \nonumber
\end{eqnarray} 
where we used the left ideal condition,
$\| P_K \|_B = 1$, 
and $\| \cdot \|_n \leq \| \cdot \|_m$ 
in the last step.
This shows that $\{ e_y \}_{y \in Y}$ is a basis
for the Fr\'echet algebra $A$. 
By the discussion in Definition 2.5,  
this basis is Schauder and equicontinuous,
and by the nuclearity of $A$, it is also 
absolute.\footnote{Note we can apply the dense ideal condition
on both sides 
$|a_x | \| e_x \|_n = 
\| a_x e_x \|_n = 
\| e_{z, ii} a e_{z, jj} \|_n
\leq C \|a \|_q$, 
to see directly that the basis 
is equicontinuous.}

By Appendix A, we can find an equivalent
family $\{ \| \cdot \|_n \}_{n=0}^\infty$  of norms for the
topology of $A$ which are increasing, and
satisfy $\| \cdot \|_0 = \| \cdot \|_B$ 
and $\| ab \|_n \leq \| a \|_n \| b \|_0$, 
$\| ba \|_n \leq \| b \|_0 \| a \|_n$ for all $b \in B$, $a \in A$.
Use these norms to 
define a family of 
scales $\sigma$ 
on $Y$ by $\sigma_n(y) = \| e_y \|_n$.
Since $\| e_y \|_B = 1$ for each $y \in Y$,
$\sigma_n \geq 1$.   Also $\sigma_0 = 1$ and 
$\sigma_n \leq \sigma_{n+1}$.
By the ideal condition, $\sigma_n(y) =\|e_y \|_n
= \| e_{z, i1} e_{z, 11} e_{z, 1j} \|_n \leq \| e_{z, 11}\|_n
= \sigma_n(z, 1, 1)$.  Similarly $\sigma_n(z, 1, 1) \leq
\sigma_n(y)$ for any $y \in Y$ with first component
$z$.  Therefore the $\sigma_n$'s are constant on each 
matrix algebra.
For each $z \in Z$, let $\ell_n(z)$
be the common value of $\sigma_n(z, i, j)$, $i,j\leq p_k$.

Apply Theorem 2.9 to see that 
$A \cong {\cal S}_\sigma^1(Y) \cong {\cal S}_\sigma^\infty(Y)$,
where, by nuclearity, $\sigma$ satisfies the 
summability condition (\ref{eqn:NucSumSpace}), and
$\ell$ satisfies the $\frak p$-summability condition 
(\ref{eqn:NucSum}).
By Theorem 6.2, 
${\cal S}^{\infty,\textrm{op}}_\ell(Y)$
is isomorphic to 
${\cal S}_\sigma^1(Y) \cong {\cal S}_\sigma^\infty(Y)$,
and hence also isomorphic to $A$.
\qed

\vskip\baselineskip
\noindent{\bf Remark 6.5. Automatically Involutive.} 
Note that the ideals of Theorems 6.4 and 6.2 are involutive.

\vskip\baselineskip
\vskip\baselineskip
\section{The Structure of Dense Nuclear Fr\'echet Ideals in $C^\star$-Algebras}

In this section, we decompose the $C^\star$-algebra 
and dense Fr\'echet ideal constructed in \S 6
into a direct sum of subideals, one with
finite multiplicities, the other with infinite multiplicities.
We investigate the properties of the subideals 
in each case.
In particular, we ask when a maximal
dense nuclear Fr\'echet ideal is
 isomorphic to the Fr\'echet space of 
standard Schwartz functions.

\vskip\baselineskip
\noindent{\bf Definition 7.1.  Decomposition Into Direct Sums of Subideals.}
Let $A$ be a dense two-sided Fr\'echet 
ideal of a $C^\star$-algebra $B$.
If $B$ is a direct sum of two-sided
$C^\star$-subideals $B_1$ and
$B_2$,  and $A$ is a direct sum of Fr\'echet subideals $A_1$
and $A_2$ such that $A_1$ is a dense Fr\'echet ideal in $B_1$ and $A_2$
is a dense Fr\'echet ideal in $B_2$, we call 
the pair of direct sums $B=B_1 \oplus B_2$,
$A=A_1 \oplus A_2$ a {\it decomposition into direct sums of subideals}.

\vskip\baselineskip
\noindent
For the dense nuclear ideals constructed in \S 6, 
a decomposition 
into direct sums of subideals occurs
when the countable set $Z$ is partitioned 
into two disjoint subsets $Z_1$ and $Z_2$.
{\it Proof:}
Since $B$ in \S 6 is already a direct sum of matrix algebras 
$\bigoplus_{z \in Z} M_{p_z} (\C)$, we can partition the
sum into two summands
\begin{equation}
B\quad =\quad \biggl( \bigoplus_{z \in Z_1} M_{p_z} (\C) \biggr)
\quad \oplus \quad \biggl( \bigoplus_{z \in Z_2} M_{p_z} (\C) \biggr)
= \quad B_1  \quad \oplus \quad B_2.
\end{equation}
Since multiplication is pointwise along $z \in Z$,
that is $f*g(z) = f(z)* g(z)$, the subalgebras
$B_1$ and $B_2$ are naturally
ideals in $B$.

In \S 6, we constructed dense nuclear Fr\'echet ideals $A$ in $B$,
and showed that every such ideal is of the
form
$A \cong {\cal S}^{\infty, \textrm{op}}_{\underline \ell}(Y)$,
% where $Y$ is defined by (\ref{eqn:Ydef}), and
where $\underline \ell$ is a family of scales on $Z$ satisfying
the $\frak p$-summability condition.
Using our partition of $Z$ into $Z_1$ and $Z_2$, we
can partition $Y$ accordingly:
\begin{eqnarray}
Y_1\quad  &=& \quad \bigcup_{z\in Z_1} \{z \}\times \{1, \dots p_z \} 
\times \{1, \dots p_z \},
\nonumber \\
Y_2\quad  &=& \quad \bigcup_{z\in Z_2} \{z \}\times \{1, \dots p_z \} 
\times \{1, \dots p_z \}.
\label{eq:Ydefs}
\end{eqnarray}
Recall that
${\cal S}^{\infty, \textrm{op}}_{\underline \ell}(Y)$ is defined
as the completion of $c_f(Y)$, and
$A_1=
{\cal S}^{\infty, \textrm{op}}_{\underline \ell \restriction_{Z_1}}(Y_1)$ and
$A_2=
{\cal S}^{\infty, \textrm{op}}_{\underline \ell \restriction_{Z_2}}(Y_2)$ 
are defined
as completions of $c_f(Y_1)$ and $c_f(Y_2)$, respectively.
Since the norms on the two ideals $A_1$ and $A_2$ are restrictions of norms
on 
${\cal S}^{\infty, \textrm{op}}_{\underline \ell}(Y)$,
and $Y$ has discrete topology, it is not hard to show that
$\{ \varphi \restriction_{Y_i}  \, | \,
\varphi \in {\cal S}^{\infty, \textrm{op}}_{\underline \ell}(Y)  \}$
is equal to 
${\cal S}^{\infty, \textrm{op}}_{\underline \ell 
\restriction_{Z_i}}(Y_i)$ for $i=1,2$.
So we have $A \cong A_1 \oplus A_2$. \qed

\vskip\baselineskip
\noindent
We partition $Z$ as follows:
\begin{eqnarray}
Z_{\textrm{fin}} &=& \biggl\{ z \in Z \, \biggl| \,
p_z \mbox{ occurs finitely many times in $\frak p$ } 
\biggr\}, \nonumber \\
\qquad
Z_{\textrm{inf}} &=& \biggl\{ z \in Z \, \biggl| \,
p_z \mbox{ occurs infinitely many times in $\frak p$ } 
\biggr\}.
\end{eqnarray}
The set $Z_{\textrm{fin}}$ could be empty, finite, or infinite.
We are mostly interested in the infinite case.
Fix some enumeration
$\zeta_{\textrm{fin}}$ of $Z_{\textrm{fin}}$, with the constraint
$p_{k-1} \leq p_k$.  
In the direct sum of subideals decomposition,
$B = B_{\textrm{fin}} \oplus B_{\textrm{inf}}$ where
\begin{eqnarray}
B_{\textrm{fin}}\quad & = &
\quad \bigoplus_{z \in Z_{\textrm{fin}}} M_{p_z}(\C) \nonumber \\
& = &
\quad \bigoplus_{k=1}^{d_{\textrm{fin}}}  
  M_{p_k}(\C).
\end{eqnarray}
The set $Z_{\textrm{inf}}$ is either empty or infinite.
Since any $p_z \in Z_{\textrm{inf}}$ occurs 
infinitely many times, we can determine
 $Z_{\textrm{inf}}$ 
by a sequence of dimensions
\begin{equation}
{\tilde {\frak p}}_{\textrm{inf}} 
= \{ {\tilde p}_k \}_{k=1}^{d_{\textrm{inf}}}\, , 
 \quad \mbox{ where ${\tilde p}_k=p_z$ for some $z \in Z_{\textrm{inf}}$,}
\end{equation}
and the ${\tilde p}_k$'s are ordered
so that ${\tilde p}_k < {\tilde p}_{k+1}$.
Then
 ${\tilde p}_1 < {\tilde p}_2 
< \cdots < {\tilde p}_{d_{\textrm{inf}}}$
 if $d_{\textrm{inf}}< \infty$
and
 ${\tilde p}_1 < {\tilde p}_2 < \cdots 
< {\tilde p}_k  < \cdots$ if $d_{\textrm{inf}}=\infty$.
We have
\begin{eqnarray}
B_{\textrm{inf}}\quad &=& \quad 
\bigoplus_{z \in Z_{\textrm{inf}}} M_{p_z}(\C) 
\nonumber \\
& = &
\quad \bigoplus_{k=1}^{d_{\textrm{inf}}}  
\,\underbrace{
 M_{{\tilde p}_k}(\C) \oplus M_{{\tilde p}_k} \oplus \cdots 
\oplus M_{{\tilde p}_k}(\C)
\oplus \cdots}_{\textrm{infinitely many times}} \nonumber \\
& = &
\quad \bigoplus_{k=1}^{d_{\textrm{inf}}}  
 c_0 \bigl( Z_k, M_{{\tilde p}_k}(\C) \bigr),
\end{eqnarray}
where  $c_0 \bigl( Z_k, M_{{\tilde p}_k}(\C) \bigr)$
is matrix-valued functions vanishing at infinity, with
 component-wise multiplication, and $Z_k=\{\, z \in Z \,|\,
p_z = {\tilde p}_k \,\}$ is a countably
infinite set for each $k$.

Similarly, we decompose the Fr\'echet ideal,
$A = A_{\textrm{fin}} \oplus A_{\textrm{inf}}$ where
\begin{equation}
A_{\textrm{fin}}\quad  = 
\quad 
\begin{cases}
B_{\textrm{fin}} & d_{\textrm{fin}} < \infty \\
{\cal S}^{\infty, 
\textrm{op}}_{{\underline \ell}\restriction_{Z_{\textrm{fin}}} }
(Y_{\textrm{fin}})
 & d_{\textrm{fin}} = \infty,
\end{cases}
\end{equation}
where $B_{\textrm{fin}}$ is finite dimensional when 
$d_{\textrm{fin}} < \infty$ and
 $Y_{\textrm{fin}}$ is
defined as $Y_1$ in (\ref{eq:Ydefs})
with $Z_1=Z_{\textrm{fin}}$.
For infinite multiplicities, we have
\begin{eqnarray}
A_{\textrm{inf}}\quad &=& \quad 
{\cal S}^{\infty, \textrm{op}}_{\underline \ell}(Y_{\textrm{inf}})
\qquad \mbox{$Y_{\textrm{inf}}$ 
defined as $Y_2$ in (\ref{eq:Ydefs})
with $Z_2=Z_{\textrm{inf}}$}
\nonumber \\
 &=& \quad 
\begin{cases}
 \bigoplus_{k=1}^{d_{\textrm{inf}}}  
 {\cal S}^{\infty, 
\textrm{op}}_{\underline \ell\restriction_{Z_k}} 
\bigl( Z_k, M_{{\tilde p}_k}(\C) \bigr)
& d_{\textrm{inf}} < \infty \\
{\cal S}^{\infty} 
\biggl(
\N^+,
 {\cal S}^{\infty, 
\textrm{op}}_{\underline \ell\restriction_{Z_k}} 
\bigl( Z_k, M_{{\tilde p}_k}(\C) \bigr) \biggr)
& d_{\textrm{inf}} = \infty, 
\end{cases}
\label{eq:Ainf}
\end{eqnarray}
where the norm of $\varphi \in A_{\textrm{inf}}$ is given by
\begin{equation}
\| \varphi \|_n = \sup_{k =1}^{d_{\textrm{inf}}}
\biggl(
 \sup_{ z \in Z_k}
\ell_n\restriction_{Z_k}(z) 
\bigl\| \varphi(k, z) \bigr\|_{\textrm{op}}
\biggr).
\end{equation}

Next we derive the appropriate summability formulas 
for the restricted scales.
\begin{eqnarray}
\sum_{z \in Z}  \,
{p_z^2 \ell_n(z) \over \ell_m(z)}
\quad &=& \quad 
\sum_{z \in Z_{\textrm{fin}}}  \,
{p_z^2 \ell_n(z) \over \ell_m(z)}
\quad + \quad
\sum_{z \in Z_{\textrm{inf}}}  \,
{p_z^2 \ell_n(z) \over \ell_m(z)}
\nonumber \\
\quad &=& \quad 
\sum_{z \in Z_{\textrm{fin}}}  \,
{p_z^2 \ell_n \restriction_{Z_{\textrm{fin}}}(z) 
\over \ell_m \restriction_{Z_{\textrm{fin}}}(z)}
\quad + \quad
\sum_{z \in Z_{\textrm{inf}}}  \,
{p_z^2 \ell_n \restriction_{Z_{\textrm{inf}}}(z) 
\over \ell_m \restriction_{Z_{\textrm{inf}}}(z)}.
\end{eqnarray}
The left summand, for finite multiplicities, becomes
\begin{equation}
\sum_{z \in Z_{\textrm{fin}}}  \,
{p_z^2 \ell_n \restriction_{Z_{\textrm{fin}}}(z) 
\over \ell_m \restriction_{Z_{\textrm{fin}}}(z)}
\quad = \quad
\sum_{k=1}^{d_{\textrm{fin}}}  \,
{ p_k^2 
\ell_n \restriction_{Z_{\textrm{fin}}}(k) 
\over \ell_m \restriction_{Z_{\textrm{fin}}}(k)}.
\end{equation}
The right summand, for infinite multiplicities, becomes
\begin{eqnarray}
\sum_{z \in Z_{\textrm{inf}}}  \,
{p_z^2 \ell_n \restriction_{Z_{\textrm{inf}}}(z) 
\over \ell_m \restriction_{Z_{\textrm{inf}}}(z)}
\quad &=& \quad
\sum_{k=1}^{d_{\textrm{inf}}}  \,
{\tilde p}_k^2 \biggl(
\sum_{z \in Z_k}  \,
{\ell_n \restriction_{Z_{\textrm{inf}}}(z) 
\over \ell_m \restriction_{Z_{\textrm{inf}}}(z)}
\biggr) \qquad \mbox{note $|Z_k|=\infty$ for each $k$} 
\nonumber \\
\quad &=& \quad
\sum_{k=1}^{d_{\textrm{inf}}}  \,
{\tilde p}_k^2 \biggl(
\sum_{z \in Z_k}  \,
{\ell_n \restriction_{Z_k}(z) 
\over \ell_m \restriction_{Z_k}(z)}
\biggr),
\end{eqnarray}
where the inner sum has infinitely many terms for every $k$.

\vskip\baselineskip \noindent
{\bf Finding a Maximal $A_{\textrm{fin}}$, a Minimal 
${\frak p}_{\textrm{fin}}$-Summable Scale on $Z_{\textrm{fin}}$.}
We restrict to the infinite dimensional case, $d_{\textrm{fin}}=\infty$.

\vskip\baselineskip
\noindent{\bf Definition 7.2.  Scales on $Y_{\textrm{fin}}$.} 
Let $\zeta_{\textrm{fin}}$ be an enumeration
of $Z_{\textrm{fin}}$, which satisfies
$p_{k-1} \leq p_k$.
Define a scale $\ell_{\textrm{fin}} 
\colon Z_{\textrm{fin}} \rightarrow \N^+$ by 
\begin{equation}
\ell_{\textrm{fin}}(k) \, =  \,
 p_1^2 \,  +\,  \cdots\,  +\, p_{k-1}^2 \, + \,
p_k^2,
\qquad \qquad
  k=1,2, \dots,
\end{equation}
and corresponding scale 
$\sigma_{\textrm{fin}} \colon Y_{\textrm{fin}} \rightarrow \N^+$ by 
\begin{equation}
\sigma_{\textrm{fin}}(k,i,j) =  \ell_{\textrm{fin}}(k),
    \qquad 
  k=1,2, \dots,
   \quad i,j=1,2, \dots p_k. 
\end{equation}
Define an enumeration $\gamma_{\textrm{fin}} \colon 
Y_{\textrm{fin}} \rightarrow \N^+$
by
\begin{eqnarray}
\gamma_{\textrm{fin}}(k,i,j)\, 
=& p_1^2 & +\,  \cdots\,  +\, p_{k-1}^2 \, + \,
 (i-1)\, +\, (j-1) p_k\, +\, 1, \nonumber \\
&&  k=1,2,\dots,
\quad i, j = 1, \dots p_k. 
\end{eqnarray}
Then $\gamma_{\textrm{fin}} \colon Y_{\textrm{fin}} 
\rightarrow \N^+$ is a
bijection of sets, since for fixed $k$, 
$\gamma_{\textrm{fin}}$ increases,
in steps of $1$ if we move down rows, one column at a time,
from $p_1^2 + \cdots p_{k-1}^2 + 1$ (at $i=j=1$) up to
$p_1^2 + \cdots p_{k-1}^2 + p_k^2$  (at $i = j = p_k$).
Note that $\gamma_{\textrm{fin}}$ orders
$\sigma_{\textrm{fin}}$, and
 $\gamma_{\textrm{fin}} \leq
\sigma_{\textrm{fin}}$.

Let ${\frak p}_{\textrm{fin}} = {\frak p} \restriction_{Z_{\textrm{fin}}}$
be the sequence of dimensions restricted to $Z_{\textrm{fin}}$.
Then ${\frak p}_{\textrm{fin}}$ is a scale on $Z_{\textrm{fin}}$
which is ordered by $\zeta_{\textrm{fin}}$, and which satisfies
${\frak p}_{\textrm{fin}}^2 \leq 
\ell_{\textrm{fin}}$.
We say that
 ${\frak p}_{\textrm{fin}}$  satisfies the {\it growth condition }
if $p_{k+1} \leq  C\bigl( k p_k\bigr)^d$, 
$k \in \N^+$, for some $C> 0$ and $d \in \N^+$.

\vskip\baselineskip
\noindent{\bf Proposition 7.3. $\ell_{\textrm{fin}}$ is a Minimal 
${\frak p}_{\textrm{fin}}$-Summable
Scale on $Z_{\textrm{fin}}$.}  
{\it 
The scale} $\ell_{\textrm{fin}}$ {\it on} $Z_{\textrm{fin}}$ {\it
satisfies the} ${\frak p}_{\textrm{fin}}${\it-
summability condition, so}
$A_{\textrm{fin}}={\cal S}^{\infty,\textrm{op}}_{\ell_{\textrm{fin}}}
(Y_{\textrm{fin}})$ {\it is a dense nuclear two-sided Fr\'echet ideal
in} $B_{\textrm{fin}}$.

{\it The scale} $\ell_{\textrm{fin}}$ {\it 
is equivalent to the scale } $\zeta_{\textrm{fin}} 
\cdot {\frak p}_{\textrm{fin}}$ {\it on} $Z_{\textrm{fin}}$.
{\it It is minimal in the sense that
if $\kappa$ is any other } $\zeta_{\textrm{fin}}$-{\it ordered scale on}
 $Z_{\textrm{fin}}$ {\it which is} 
 ${\frak p}_{\textrm{fin}}${\it-summable,
 and } $\kappa \lesssim \ell_{\textrm{fin}}$,
{\it then} $\kappa \thicksim \ell_{\textrm{fin}}$.

{\it The scale} $\sigma_{\textrm{fin}}$
 {\it on }  
 $Y_{\textrm{fin}}$ 
{\it is equivalent to an
enumeration of} $Y_{\textrm{fin}}$ {\it
if and only if the sequence of dimensions}
 ${\frak p}_{\textrm{fin}}$  {\it satisfies the growth condition.}

\vskip\baselineskip
\noindent{Proof:} 
Since $\zeta_{\textrm{fin}}(k)=k \leq \ell_{\textrm{fin}}(k)$ and
 ${\frak p}_{\textrm{fin}}
\leq
 {\frak p}_{\textrm{fin}}^2
 \leq \ell_{\textrm{fin}}$,
the product
 $\zeta_{\textrm{fin}}
\cdot {\frak p}_{\textrm{fin}}$
is always less than or equal to
$\ell_{\textrm{fin}}^2$,
so
 $\ell_{\textrm{fin}}$ dominates 
the scale $\zeta_{\textrm{fin}}
\cdot {\frak p}_{\textrm{fin}}$,
which we noted in Remark 6.3 was 
 ${\frak p}_{\textrm{fin}}$-summable.
Therefore 
 $\ell_{\textrm{fin}}$ is
${\frak p}_{\textrm{fin}}$-summable.
By Theorem 6.2, 
$A_{\textrm{fin}}={\cal S}^{\infty,\textrm{op}}_{\ell_{\textrm{fin}}}
(Y_{\textrm{fin}})$  is a dense nuclear two-sided Fr\'echet ideal
in $B_{\textrm{fin}}$.

By definition,
\begin{eqnarray}
 \ell_{\textrm{fin}}(k) &=& p_1^2 + \cdots + p_{k-1}^2 +  p_k^2
\nonumber \\
&\leq &
 \underbrace{p_k^2 + \cdots + p_k^2 + p_k^2}_{k \textrm{ times}} 
\qquad \mbox{since the $p_k$'s are ordered}
\nonumber \\
& = &
 k p_k^2 
\nonumber \\
& = &
 \zeta_{\textrm{fin}}(k) {\frak p}_{\textrm{fin}}(k)^2  
\nonumber \\
& \leq &
 \zeta_{\textrm{fin}}(k)^2 {\frak p}_{\textrm{fin}}(k)^2,  
\end{eqnarray}
so $\zeta_{\textrm{fin}}
\cdot {\frak p}_{\textrm{fin}}$ dominates
 $\ell_{\textrm{fin}}$.  Hence
 $\zeta_{\textrm{fin}}
\cdot {\frak p}_{\textrm{fin}} \thicksim
 \ell_{\textrm{fin}}$.

Let $\kappa$ be a $\zeta_{\textrm{fin}}$-ordered scale on
$Z_{\textrm{fin}}$ which is 
 ${\frak p}_{\textrm{fin}}$-summable.
Then $\zeta_{\textrm{fin}} \leq \kappa$.
Since
\begin{equation}
\sum_{k=1}^\infty {p_k \over \kappa(k)^p} < \infty,
\end{equation}
for some $p \in \N^+$, the terms in the series are
bounded so there is a $C>0$ such that $p_k \leq C \kappa(k)^p$
for $k \in \N^+$.
Then we can bound the product $ \zeta_{\textrm{fin}}(k) p_k 
\leq C \kappa(k)^{p+1}$, so $\zeta_{\textrm{fin}} \cdot 
{\frak p}_{\textrm{fin}} \lesssim \kappa$.
With the additional assumption that $\kappa \lesssim \ell_{\textrm{fin}}$,
we have
$\zeta_{\textrm{fin}} \cdot 
{\frak p}_{\textrm{fin}} \lesssim \kappa
 \lesssim \ell_{\textrm{fin}}
\lesssim
\zeta_{\textrm{fin}} \cdot 
{\frak p}_{\textrm{fin}}$, and so $\kappa \thicksim \ell_{\textrm{fin}}$.

For the last paragraph, assume the growth condition holds.
We show that 
 $\sigma_{\textrm{fin}}$ is
equivalent to 
 $\gamma_{\textrm{fin}}$.
Define a scale $\beta_{\textrm{fin}} 
\colon Y_{\textrm{fin}} \rightarrow \N^+$ by 
\begin{equation}
\beta_{\textrm{fin}}(k,i,j) \, = 
\begin{cases} 
 1 
 & 
 k=1, \quad \quad i,j=1,2,\dots p_1, \\
 p_1^2 \,  +\,  \cdots\,  +\, p_{k-1}^2 
 &
  k=2,3, \dots,
  \quad i,j=1,2,\dots p_k.
\end{cases} 
\end{equation}
By the growth condition,
\begin{eqnarray}
\sigma_{\textrm{fin}}(k) &=& p_1^2 
+ \biggl( p_2^2 + \cdots + p_k^2 \biggr)
\nonumber \\
& \leq & p_1^2 + \biggl( (C\cdot (1\cdot p_1)^d )^2 + 
\cdots + (C \cdot (k-1 \cdot p_{k-1})^d)^2 \biggr)
\nonumber \\
& \leq & p_1^2 + \biggl( (C\cdot (k-1\cdot p_1)^d )^2 + 
\cdots + (C \cdot (k-1 \cdot p_{k-1})^d)^2 \biggr)
\nonumber \\
& = & p_1^2 + C^2 (k-1)^{2d} \biggl( p_1^{2d} + \cdots + p_{k-1}^{2d} \biggr)
\nonumber \\
& \leq &  p_1^2 + C^2 \beta_{\textrm{fin}}(k)^{2d} \beta_{\textrm{fin}}(k)^d
\nonumber \\
& \leq &  (p_1^2 + C^2)  \beta_{\textrm{fin}}(k)^{3d},
\end{eqnarray}
which holds for all $k \in \N^+$.  
So $\sigma_{\textrm{fin}} \lesssim \beta_{\textrm{fin}}$.
Since $\beta_{\textrm{fin}} \leq \gamma_{\textrm{fin}} 
\leq \sigma_{\textrm{fin}}$, we have 
$\beta_{\textrm{fin}} 
\lesssim \gamma_{\textrm{fin}} \lesssim 
\sigma_{\textrm{fin}} \lesssim \beta_{\textrm{fin}}$ 
and so $\sigma_{\textrm{fin}} \thicksim \gamma_{\textrm{fin}}$.

Conversely, let $\gamma: Y_{\textrm{fin}} \cong \N^+$ be an enumeration
which is equivalent to $\sigma_{\textrm{fin}}$.
We prove that the growth condition holds.
Define 
a scale
$\ell(z) = \min_{ij} \gamma(z, i, j)$, for $z \in Z_{\textrm{fin}}$.
Then $\ell \thicksim \ell_{\textrm{fin}}$
since
\begin{eqnarray}
\ell(z) & \leq & \gamma(z,i,j) \qquad \mbox{since $\ell$ is the min}
\nonumber \\ 
& \leq &  
C \sigma_{\textrm{fin}}(z,i,j)^d \qquad \mbox{since $\sigma_{\textrm{fin}}
\thicksim \gamma$ by assumption}
\nonumber \\
& = &  
C \ell_{\textrm{fin}}(z)^d 
\\
\ell_{\textrm{fin}}(z) 
& = & 
\sigma_{\textrm{fin}}(z,i,j) \qquad \mbox{for any $i,j=1, \dots p_z$}
\nonumber \\ 
&\leq & 
C \gamma(z,i,j)^d. \qquad 
 \mbox{since $\sigma_{\textrm{fin}}
\thicksim \gamma$ by assumption}
\label{eq:halfEllsSimple}
\end{eqnarray}
Since the left hand side of  (\ref{eq:halfEllsSimple})
is independent of $i$ and $j$, we can take the min
over $i$ and $j$ to get $\ell_{\textrm{fin}}(z) \leq C\ell(z)^d $.

Let $\zeta$ be an enumeration of $Z_{\textrm{fin}}$,
with inverse $\rho$, which satisfies
$\ell(\rho(1)) <  \ell(\rho(2)) <  \dots$.
Since $\gamma$ maps $Y_{\textrm{fin}}$ onto $\N^+$, 
$\ell(\rho(1))=1$.
The smallest the set of values 
$\{ \gamma(\rho(1),i, j)\}_{i,j\,=\,1,\dots p_{\rho(1)}}$
could be is $\{1, 2, \dots p_{\rho(1)}^2\}$.  So $\ell(\rho(2))$
can be no bigger than $p_{\rho(1)}^2 + 1$.
Similarly $\ell(\rho(k)) \leq p_{\rho(1)}^2 
+ p_{\rho(2)}^2 + \cdots   p_{\rho(k-1)}^2 + 1$.

Since $\sigma_{\textrm{fin}} \thicksim \gamma$, 
there is some $d \in \N^+$  and $C_1>0$ such
that $\gamma(z, i, j) \leq C_1 \ell_{\textrm{fin}}(z)^d$, 
$z \in Z_{\textrm{fin}}$,
$i,j =1,\dots p_z$.
Since $\gamma$ is one to one, for any $z \in Z_{\textrm{fin}}$
the set of values $\{ \gamma(z, i, j) \}_{i,j=1,\dots p_z}$ 
must contain a number as big as $p_z^2$.
Hence $p_z^2 \leq C_1 \ell_{\textrm{fin}}(z)^d$.
So we have 
\begin{eqnarray}
p_{\rho(k)}^2 & \leq &  C_1 \ell_{\textrm{fin}}(\rho(k))^d 
 \nonumber \\
& \leq & C_1 C_d  \ell(\rho(k))^m \qquad 
\qquad \qquad \text{since $\ell \thicksim {\ell_{\textrm{fin}}}$} 
\nonumber \\
& \leq &
C_1 C_d (
 p_{\rho(1)}^2 
+ p_{\rho(2)}^2 + \cdots   p_{\rho(k-1)}^2 
 + 1)^m.
\label{eq:rhoGrowthSimple}
\end{eqnarray}

Let $\rho_{\textrm{fin}} = \zeta_{\textrm{fin}}^{-1}$.
We would like the growth condition
(\ref{eq:rhoGrowthSimple}) to have $\rho_{\textrm{fin}}$'s
instead of $\rho$'s.
If $p_{\rho(k-1)} > p_{\rho(k)}$ for some $k \in 2,3,\dots$,
let $\rho'$ be $\rho$ with the $k$ and $k-1$th values
switched.  Then the growth condition
(\ref{eq:rhoGrowthSimple}) still holds true for $\rho'$:
\begin{eqnarray}
p_{\rho'(k-1)}^2 < p_{\rho(k-1)}^2 & \leq  &
C_1 C_d (
 p_{\rho(1)}^2 
+ p_{\rho(2)}^2 + \cdots   p_{\rho(k-2)}^2 
 + 1)^m \nonumber \\
&=&
C_1 C_d (
 p_{\rho'(1)}^2 
+ p_{\rho'(2)}^2 + \cdots p_{\rho'(k-2)}^2 
 + 1)^m, 
\\
p_{\rho'(k)}^2 = p_{\rho(k-1)}^2 & \leq  &
C_1 C_d (
 p_{\rho(1)}^2 
+ p_{\rho(2)}^2 + \cdots   p_{\rho(k-2)}^2 
 + 1)^m \nonumber \\
&=&
C_1 C_d (
 p_{\rho'(1)}^2 
+ p_{\rho'(2)}^2 + \cdots p_{\rho'(k-2)}
 + 1)^m 
\nonumber \\
&< &
C_1 C_d (
 p_{\rho'(1)}^2 
+ p_{\rho'(2)}^2 + \cdots p_{\rho'(k-2)}
  + p_{\rho'(k-1)}
 + 1)^m.
\end{eqnarray}
Since we didn't change $C_1$, $C_d$, or $m$, we can make 
as many of these switches as we need, and the same
growth condition (\ref{eq:rhoGrowthSimple}) continues
to hold.
Continuing all the way up the sequence gives a non-decreasing
ordering $\rho'(k-1) \leq \rho'(k)$ for all $k=2,3,\dots$.
Next note that by the same argument, 
if $\rho'(k-1) = \rho'(k)$, we can swap
these two equal values and the growth condition will hold
as well.  Since $\rho_{\textrm{fin}}$ can be obtained 
from $\rho'$ in this way, we've shown that 
 the growth condition (\ref{eq:rhoGrowthSimple}) holds for
 $\rho_{\textrm{fin}}$.
We have
\begin{eqnarray}
p_k^2 
& \leq &  
C_1 C_d (
 p_1^2 
+ p_2^2 + \cdots   p_{k-1}^2 
 + 1)^m
\nonumber \\
& \leq &  
C_1 C_d \bigl(
 2 ( p_1^2 
+ p_2^2 + \cdots   p_{k-1}^2 
 )\bigr)^m \qquad \mbox{since $1 \leq p_i$ for each $i$}
\nonumber \\
& \leq &  
C_1 C_d 2^m (
 k-1 \cdot  
  p_{k-1}^2 
 )^m. \qquad \mbox{since $p_i \leq p_{k-1}$ for $i < k-1$}
\end{eqnarray}
It follows that the growth condition 
$p_k \leq C (k-1 \cdot p_{k-1})^m$ holds for each $k=2,3,\dots$.
\qed

\vskip\baselineskip \noindent
{\bf Making $A_{\textrm{inf}}$ Standard Schwartz.}
Assume that 
$A_{\textrm{inf}}$ is nonempty, so 
$d_{\textrm{inf}} \geq 1$.
Recall from (\ref{eq:Ainf}) that 
 $A_{\textrm{inf}}$ is a direct sum of Schwartz spaces.
We can factor out the first direct summand:
\begin{eqnarray}
A_{\textrm{inf}} & = & 
\begin{cases}
 {\cal S}^{\infty, 
\textrm{op}}_{\underline \ell\restriction_{Z_1}} 
\bigl( Z_1, M_{{\tilde p}_1}(\C) \bigr)
\oplus
 \bigoplus_{k=2}^{d_{\textrm{inf}}}  
 {\cal S}^{\infty, 
\textrm{op}}_{\underline \ell\restriction_{Z_k}} 
\bigl( Z_k, M_{{\tilde p}_k}(\C) \bigr)
& d_{\textrm{inf}} < \infty \\
 {\cal S}^{\infty, 
\textrm{op}}_{\underline \ell\restriction_{Z_1}} 
\bigl( Z_1, M_{{\tilde p}_1}(\C) \bigr)
\oplus
{\cal S}^{\infty} 
\biggl(
\N^+ \setminus \{1\},
 {\cal S}^{\infty, 
\textrm{op}}_{\underline \ell\restriction_{Z_k}} 
\bigl( Z_k, M_{{\tilde p}_k}(\C) \bigr) \biggr)
& d_{\textrm{inf}} = \infty
\end{cases} \nonumber \\
&=& {\cal S}^{\infty, 
\textrm{op}}_{\underline \ell\restriction_{Z_1}} 
\bigl( Z_1, M_{{\tilde p}_1}(\C) \bigr)
\oplus
 {\cal S}^{\infty, 
\textrm{op}}_{\underline \ell\restriction_{Z_{\textrm{inf}}\setminus Z_1}} 
\bigl( Y_{\textrm{inf}}\setminus Y_1 \bigr).
\label{eq:AinfStand}
\end{eqnarray}
By taking ${\underline \ell} \restriction_{Z_1}$
to be an enumeration of the infinite set $Z_1$,
the first direct summand becomes standard Schwartz
functions on an infinite set (\ref{eqn:standSchw}),
and the second direct summand 
is Schwartz functions on $Y_{\textrm{inf}} \setminus Y_1$,
which vanish rapidly with respect to 
powers of a single scale.

\vskip\baselineskip \noindent
We now prove a proposition and two lemmas, 
required in the proof of Corollary 7.7.

\vskip\baselineskip
\noindent{\bf Proposition 7.4.  Direct Sums of Infinite Type Power Series.}
{\it 
If a nuclear Fr\'echet space $A$ is a direct sum of 
power series spaces $A_1$ and $A_2$ of infinite type, then
$A$ is also a power series space of infinite type.  
The scale obtained for $A$ is 
dominated by the scales defining
$A_1$ and $A_2$.
}

\noindent{Proof:} 
Assume that $A_1$ and $A_2$ are
power series spaces of infinite type,
which means they are isomorphic to
$\ell^1$-norm $\sigma_i$-rapidly vanishing functions
for single scales
$\sigma_i$, $i=1,2$, respectively 
(Remark 2.8).
By nuclearity, each scale is summable, 
and $A_i \cong {\cal S}_{\sigma_i}^1(X_i)
 \cong {\cal S}_{\sigma_i}^\infty(X_i)$ for
countably infinite sets $X_i$ (Theorem 2.9).
Let $\sigma$ be the scale on 
the disjoint union $X = X_1 \cup X_2$ defined by
\begin{equation}
\sigma(x) = \begin{cases} \sigma_1(x) & x \in X_1 \\
\sigma_2(x) & x \in X_2, \end{cases}
\label{eq:sigDef}
\end{equation}
for $x \in X$.
For $\varphi \in {\cal S}^\infty_{\sigma}(X)$,
\begin{eqnarray}
\bigl\| \varphi \bigr\|^\infty_n &=&
\sup_{x \in X} \sigma(x)^n \bigl| \varphi(x) \bigr|
\nonumber \\
&=&
\max \biggl\{
\sup_{x \in X_1} \sigma(x)^n \bigl| \varphi(x) \bigr|, \quad
\sup_{x \in X_2} \sigma(x)^n \bigl| \varphi(x) \bigr|
\biggr\}
\qquad \mbox{since $X=X_1 \cup X_2$}
\nonumber \\
&=&
\max \biggl\{
\sup_{x \in X_1} \sigma_1(x)^n 
\bigl| \varphi \restriction_{X_1} (x) \bigr|, \quad
\sup_{x \in X_2} \sigma_2(x)^n 
\bigl| \varphi \restriction_{X_2} (x) \bigr|
\biggr\}
\qquad \mbox{by (\ref{eq:sigDef})}
\nonumber \\
&=&
\max
\biggl\{ 
\bigl\| \varphi \restriction_{X_1} \bigr\|_n^\infty, \quad
\bigl\| \varphi \restriction_{X_2} \bigr\|_n^\infty
\biggr\}.
\end{eqnarray}
Hence ${\cal S}^\infty_{\sigma}(X) \cong
{\cal S}^\infty_{\sigma_1}(X_1) \oplus 
{\cal S}^\infty_{\sigma_2}(X_2)$, and so
$A \cong {\cal S}^\infty_{\sigma}(X)$.
By nuclearity
$A \cong {\cal S}^1_{\sigma}(X)$ as well, so $A$ is a power series
space of infinite type.

\vskip\baselineskip
\noindent{\bf Lemma 7.5.  Injective Scales With Partitioned Range.}
{\it 
Let $\sigma$ be a summable scale on a countably infinite set $X$. 
Then there exists an equivalent scale $\tau$ on $X$ which
is one-to-one and integral valued.
Furthermore, if $d \in \N^+$ and $i \in \N$ with $i < d$,
the range of $\tau$ can be made to lie inside
the set $d\N^+ + i$.
}

\noindent{Proof of Lemma 7.5:} 
We saw in Corollary 2.13 that $\sigma$ dominates an enumeration
$\gamma$ of $X$ which orders $\sigma$.
Define a new scale on $X$ by $\tau_1 = \gamma \sigma$.
Since $\sigma$ dominates $\gamma$, 
it is easy to show that $\sigma \thicksim \tau_1$.
If $x_k$, $x_{k+1}$ are two consecutive elements of
$X$, $\tau_1(x_k) = \tau_1(x_{k+1})$ would imply
$\sigma(x_k) = \sigma(x_{k+1}){k+1 \over k}$,
contradicting $\sigma(x_k) \leq \sigma(x_{k+1})$.
So $\tau_1 \colon X \rightarrow [1, \infty)$ is one-to-one.
Moreover 
\begin{equation}
\tau_1(x_{k+1}) - \tau_1(x_{k}) 
= (k+1)\sigma(x_{k+1})  - k\sigma(x_k)  
\geq (k+1-k) \sigma(x_k)  = \sigma(x_k) \geq 1.
\label{eq:oneStretch}
\end{equation}
Define a new scale on $X$ by 
$\tau_2 = \lceil \tau_1 \rceil$,
 the least integer greater than $\tau_1$.
By (\ref{eq:oneStretch}), $\tau_2$ is strictly increasing and
one-to-one.  Since $\tau_1 \leq \tau_2$ and $\tau_2 \leq \tau_1 + 1$,
we have $\tau_1 \thicksim \tau_2$.
Finally define 
a third equivalent scale 
$\tau_3 = d \tau_2 + i$  to see the last statement
of the Lemma.
\qed

\vskip\baselineskip
\noindent{Proof of last statement of Proposition 7.4:} 
By Lemma 7.5, 
we may assume that $\sigma_1$ and $\sigma_2$ are one-to-one,
and that they assume only odd and even integral values, respectively.
Define $\sigma$ as before, using (\ref{eq:sigDef}).
Then $\sigma$ is one-to-one and integral-valued,
and range$(\sigma)$=
range$(\sigma_1)\cup$
range$(\sigma_2)$.
Let $\gamma$ be the (unique) enumeration of $X$ that orders $\sigma$,
and let $\gamma_1$, $\gamma_2$ be the enumerations of $X_1,X_2$
which order $\sigma_1$, $\sigma_2$, respectively.
We want to show $\sigma(\gamma^{-1}(k)) \leq \sigma_1(\gamma_1^{-1}(k))$ 
and $\sigma(\gamma^{-1}(k)) \leq \sigma_2(\gamma_2^{-1}(k))$ for
all $k \in \N^+$, which is what is meant by
the last statement of the Proposition.

First show that
\begin{equation}
\gamma_1(x) \leq \gamma(x) \mbox{ for $x \in X_1$}
\qquad \mbox{ and } \qquad 
\gamma_2(x) \leq \gamma(x) \mbox{ for $x \in X_2$.}
\label{eq:order12}
\end{equation}
Assume we have shown this for $x_1, x_2, \dots x_k \in X$,
ordered with respect to $\gamma$.
If $x_{k+1} \in X_1$, 
then $\sigma_1(x_{k+1})$ is 
the smallest $\sigma_1(x)$ that has
not already occured in the list 
$\sigma(x_1)$, $\sigma(x_2)$, $\dots  \sigma(x_k)$.
Evidently, $\gamma_1(x_{k+1})\leq k+1 = \gamma(x_{k+1})$.
Similarly if $x_{k+1}$ had been in $X_2$,
we find $\gamma_2(x_{k+1}) \leq \gamma(x_{k+1})$.

It follows from (\ref{eq:order12})
that for each $k \in \N^+$, 
$k= \gamma_1(\gamma_1^{-1}(k)) \leq
\gamma(\gamma_1^{-1}(k))$
and similarly $k \leq \gamma(\gamma_2^{-1}(k))$.
Applying $\sigma \circ \gamma^{-1}$ to both
sides of the first inequality, we get 
$\sigma(\gamma^{-1}(k)) \leq
 \sigma(\gamma^{-1}\bigl(  \gamma(\gamma_1^{-1}(k)) \bigr))
= \sigma(\gamma_1^{-1}(k)) 
= \sigma_1(\gamma_1^{-1}(k)) $.
Similarly
$ \sigma(\gamma^{-1}(k)) \leq
 \sigma_2(\gamma_2^{-1}(k))$.
\qed

\vskip\baselineskip
\noindent{\bf Lemma 7.6.  Dominated by an Enumeration.} 
{\it If a summable 
scale on a countable set $X$ is dominated
by an enumeration of $X$, then the scale is equivalent
to an enumeration which orders it.}

\noindent
{Proof:}
Let $\sigma$ be a scale on $X$, and assume
\begin{equation}
\sigma(x) \leq C \gamma(x)^d, \qquad x \in X,
\end{equation}
for some $C>0$, $d \in \N^+$, and $\gamma$ an enumeration
of $X$.
Let $x_1, x_2 \in X$ with $\sigma(x_1) < \sigma(x_2)$
but $\gamma(x_2) < \gamma(x_1)$.
Then
\begin{eqnarray}
\sigma(x_1) &<& \sigma(x_2) \leq C\gamma(x_2)^d,
\nonumber \\
\sigma(x_2) &\leq &  C\gamma(x_2)^d < C \gamma(x_1)^d.
\end{eqnarray}
So we if $\gamma'$ is $\gamma$ with the values for $x_1$ and $x_2$
switched, we still have $\sigma \leq C\gamma'^d$. 
Since $C$ and $d$ are the same, we can make as many 
switches as we like, resulting in an enumeration $\gamma''$
which both satisfies the inequality $\sigma \leq C \gamma''^d$
and orders $\sigma$.
By Corollary 2.13, $\sigma$ dominates $\gamma''$, so
$\sigma \thicksim \gamma''$.
\qed

\vskip\baselineskip
\noindent{\bf Corollary 7.7. Maximal Fr\'echet Space
 Structure of The Direct Sum 
$A = A_{\textrm{fin}} \oplus A_{\textrm{inf}}$.}
{\it 
There exists a
${\frak p}$-summable scale on $Z$, which makes 
$A$, 
}$A_{\textrm{fin}}$, {\it and }
$A_{\textrm{inf}}$  {\it
 power series spaces of infinite type.

If } $d_{\textrm{fin}} = \infty$, {\it and } ${\frak p}_{\textrm{fin}}$
{\it satisfies the growth condition, then
using the scale }$\ell_{\textrm{fin}}$ {\it on } $Z_{\textrm{fin}}$, 
{\it both $A$ and } $A_{\textrm{fin}}$ {\it are isomorphic as Fr\'echet spaces
to standard Schwartz functions.

If } $d_{\textrm{inf}}\geq 1${\it, 
by choosing the scale on $Z_1 = \{ z \in Z \mid p_z = {\tilde p}_1 \}$
to be an enumeration, then
$A$  and } $A_{\textrm{inf}}$ {\it are both isomorphic as Fr\'echet spaces
to standard Schwartz functions.

If } 
$d_{\textrm{fin}} = \infty$,
 $d_{\textrm{inf}} = 0$,
 {\it and } ${\frak p}_{\textrm{fin}}$
{\it does not satisfy the growth condition, then
} $A=A_{\textrm{fin}}$ {\it is not isomorphic as a Fr\'echet space
to standard Schwartz functions.}

\vskip\baselineskip \noindent
See (\ref{eqn:standSchw})
for the definition of
 standard Schwartz functions, and Remark 2.8 for the definition
of power series space of infinite type.
By Corollary 2.13, 
every summable scale dominates any enumeration with the same ordering.
In this sense standard Schwartz functions are maximal
among nuclear power series spaces of infinite type.

Proposition 7.3 gives a 
minimality property of the scale $\ell_{\textrm{fin}}$,
which corresponds to the maximality of $A_{\textrm{fin}}$.
That result applies indepenently of 
the growth condition being satisfied.

\vskip\baselineskip \noindent
{Proof of Corollary 7.7:} 
As noted in Remark 6.3, the scale $\ell = \zeta {\frak p}$ 
on $Z$ is $\frak p$-summable, making
$A={\cal S}^{\infty, \textrm{op}}_{\ell}(Y)$,
$A_{\textrm{fin}}=
{\cal S}^{\infty, \textrm{op}}_{\ell\restriction_{Y_{\textrm{fin}}}}
(Y_{\textrm{fin}})$ and
$A_{\textrm{inf}}=
{\cal S}^{\infty, \textrm{op}}_{\ell\restriction_{Y_{\textrm{inf}}}}
(Y_{\textrm{inf}})$ all power series spaces of infinite
type.

By the last paragraph of Proposition 7.3, 
if the growth condition is satisfied, the scale $\ell_{\textrm{fin}}$
makes $A_{\textrm{fin}}$ isomorphic to standard Schwartz functions
on $Y_{\textrm{fin}}$.
Writing $A$ as a direct sum of $A_{\textrm{fin}}$ and 
 $A_{\textrm{inf}}$, with all three power series spaces of infinite type,
 Proposition 7.4 tells us the scale $\sigma$ on $Y$ used to define $A$
is dominated by the scales on both summands.  In particular
$\sigma$ is dominated by $\sigma_{\textrm{fin}}$, which is equivalent
to an enumeration.  
By Lemma 7.6, $\sigma$ is equivalent to an enumeration, and $A$
is standard Schwartz.

If $d_{\textrm{inf}} \geq 1$, we saw in (\ref{eq:AinfStand})
that $A_{\textrm{inf}}$ has a direct summand of functions
vanishing rapidly with respect to an enumeration.
Proposition 7.4 says the scale for 
 $A_{\textrm{inf}}$ is dominated by this enumeration.  By Lemma 7.6,
the scale is equivalent to some enumeration which orders it, and
therefore $A_{\textrm{inf}}$ is standard Schwartz.
Since $A$ has $A_{\textrm{inf}}$  as a direct summand, the
same argument shows that $A$ is standard Schwartz.

If $d_{\textrm{inf}}=0$, then $A=A_{\textrm{fin}}$.
If the growth condition is not satisfied, 
then by the last paragraph of Proposition 7.3, 
$\sigma_{\textrm{fin}}$ is not equivalent
to an enumeration of $Y_{\textrm{fin}}$.
Assume for a contradiction that
 ${\cal S}^1_{\sigma_{\textrm{fin}}}(Y_{\textrm{fin}})$
is isomorphic to standard Schwartz functions
${\cal S}^1_{\gamma}(X)$ for some enumeration $\gamma$ of $X$.
By the last paragraph of Proposition 2.14, 
$\sigma_{\textrm{fin}}$ is semi-equivalent to $\gamma$.
By Proposition 2.16, 
$\sigma_{\textrm{fin}}$ is equivalent to $\gamma$,
which gives the contradiction.
\qed

\vskip\baselineskip
\vskip\baselineskip
\section{Examples and Remarks}

\vskip\baselineskip
\noindent{\bf Remark 8.1. Reordering $Z_{\textrm{fin}}$ 
to Satisfy the Growth Condition.}  
If each $p_k$ occurs only finitely many times, 
and the growth
condition is satisfied for an enumeration of the
sequence $\{ p_k \}_{k=1}^\infty$,
we can always reorder so that $p_k \leq p_{k+1}$.
For let $C, d$ be such that $p_{k+1} \leq
C(kp_k)^d$.  
If $p_{k+1}$ is the
first out of order element of the sequence, find
the smallest $l> k+1$  for which $p_l$ belongs
at the $k+1$th spot.   
Define $\beta(k+1) = l$, $\beta(k+2)=k+1$, $\dots
\beta(l) = l-1$, and $\beta(i) = i $ otherwise.
Then $\{ p_{\beta(k)} \}_{k=1}^\infty$ is in order up to 
 $k+1$, and
\begin{eqnarray}
p_{\beta(k+1)} &=& p_l < p_{k+1} \leq C(kp_k)^d = C(k p_{\beta(k)})^d
\nonumber \\
p_{\beta(k+2)} &=& p_{k+1} \leq C(k p_k)^d = C(k p_{\beta(k)})^d
\leq C(k p_{\beta(k+1)})^d 
\leq C(k + 1 \cdot p_{\beta(k+1)})^d 
\nonumber \\
& \vdots & \nonumber \\
p_{\beta(l)} &=& p_{l-1} \leq C(l-2 \cdot p_{l-2})^d 
= C(l-2 \cdot p_{\beta(l-1)})^d
\leq C(l-1 \cdot p_{\beta(l-1)})^d 
\nonumber \\
p_{\beta(l+1)} &=& p_{l+1} \leq C(l \cdot p_l)^d 
\leq C(l \cdot p_{l-1})^d
= C(l \cdot p_{\beta(l)})^d,
\,\, \mbox{since $p_{l-1} \geq p_{k+1} > p_l$}
\end{eqnarray}
so the growth condition holds with the same constants $C, d$.
Continuing all the way up the sequence gives
a nondecreasing ordering, which still satisfies the
growth condition.

\vskip\baselineskip
\noindent{\bf Example 8.2. Polynomial Growth.}  
If $p_k \leq Ck^d$, $k \in \N^+$, 
then since $p_k \geq 1$ for each $k$,
\begin{equation}
p_{k+1}\, \leq\, C(k+1)^d\, \leq\,  2^d C k^d
\, \leq \,  2^d C (k p_k)^d, 
\end{equation}
so the growth condition
holds with constants $2^d C $, $d$.\footnote{Note that 
${k \over k+1} \geq {1\over 2}$.}
The sequence $p_k = e^k$ does not have polynomial
growth, but satisfies the growth
condition with $C=e$, $d=1$ since 
$p_{k+1} = e^{k+1} = e \cdot e^k = e \cdot p_k \leq e (k p_k)^1$.

\vskip\baselineskip
\noindent{\bf Remark 8.3. Alternative Construction
of Nuclear Fr\'echet Ideals.} 
Let $\underline \ell$ and 
$\underline \sigma$ be as in Theorem 6.2, with
the $\frak p$-summability condition (\ref{eqn:NucSum}) satisfied.
Pick a family of scales $\underline \beta$ on $Y$ so that
$\ell_n(z) \leq \beta_n(i,j,z)$.
Also insure that for each $n \in \N$ there is a
sufficiently large $m \in \N$ so that
$\beta_n(i,j,z)\leq \ell_m(z)$.
Then clearly ${\underline \sigma} \thicksim {\underline \beta}$.
So ${\cal S}^1_{\underline \beta} (Y)$
is a dense nuclear two-sided Fr\'echet ideal 
in the $C^\star$-algebra $B$,
by its isomorphism with
${\cal S}^{\infty,\textrm{op}}_{\underline \sigma} (Y)$.

\vskip\baselineskip
\noindent{\bf Example 8.4.  $C^\infty(G)$, $G$ Compact Connected Lie Group.}
\noindent
By [Sug, 1971], Theorem 4, the Fr\'echet space $C^\infty(G)$ is 
isomorphic to
standard Schwartz functions on $D$, where $D$ is 
the set of all dominant $G$-integral
forms on the Lie algebra of a maximal toral subgroup of $G$.
For example, consider the circle group $G=\T$.
The Fr\'echet space $C^\infty(\T)$
is topologized by seminorms $\| \varphi \|_i = 
\| \partial^i \varphi \|_\infty$,
$\varphi \in C^\infty(\T)$, where $\partial = {d \over d\theta}$.
The Fourier transform $\widehat \ $ changes 
$C^\infty(\T)$ into ${\cal S}(\Z)$, $C^\star(\T)$ into $c_0(\Z)$,
and convolution multiplication of functions 
into pointwise multiplication.
The transformed seminorms are 
$\| {\widehat \varphi } \|_i = \sup_{k \in \Z} 
|k^i {\widehat \varphi}(k)|$.

It is well-known that $C^\infty (G)$ is a dense two-sided 
Fr\'echet ideal in the convolution algebras $L^1(G)$
and $C^\star(G)$.

\vskip\baselineskip
\noindent{\bf Example 8.5.  $C^\infty$ Functions on the Cantor Group.}
The Cantor group $K$ is a compact locally compact group, which is
not a Lie group.  It is totally disconnected
and can be described as the dual group of the discrete abelian
group $\widehat{K} = \Z [\frac{1}{2}]/\Z$ of dyadic rationals
from $0$ to $1$.   
Define a proper scale $\sigma(\frac{l}{2^p}) =2^p$
on $\widehat{K}$, for any $p \in \N$.
Here $l$ is any positive odd number less than $2^p$.
The additive identity is $0$, when $p=0$ and $l=0$.
Let $A ={\cal S}_\sigma^\infty (\widehat{K})$ be 
$\sigma$-rapidly vanishing functions on $\widehat{K}$, with sup-norm.
$A$ is a dense Fr\'echet ideal in $B=c_0(\widehat{K})$, the
$C^\star$-algebra of functions vanishing at
$\infty$ on $\widehat{K}$, with pointwise multiplication.  
(Note that convolution multiplication on $K$ 
becomes pointwise multiplication on $\widehat{K}$
via the Fourier transform.)
Since $\sigma$ is summable,
$A$ is also nuclear
[Sch, 1998], Lemma 1.2 with $q=2$, Theorem 1.6.

We show that $A$ is standard Schwartz functions.
Let 
$\gamma \colon \widehat{K} \rightarrow \N^+$
be the map 
$\gamma(\frac{l}{2^p}) = 2^{p-1} + \lfloor l/2\rfloor + 1$,
where $\lfloor  \cdot \rfloor$ is the largest integer
not greater than, 
and we set $\gamma(0) = 1$.
For each $p \in \N^+$, note that $\lfloor l/2 \rfloor$ goes from
$0$ to $2^{p-1}-1$, in increments of $1$.
The scale $\gamma$ is an enumeration of
$\widehat{K}$, and
 $\gamma \leq \sigma$.  
Also $\sigma \leq 2 \gamma$,
so $\gamma$ and $\sigma$ are equivalent 
scales on $\widehat{K}$.

\vskip\baselineskip
\noindent{\bf Remark 8.6.  If Multiplication is Trivial,  
Dense Nuclear Ideals Always Exist. }
Let $B$ be any separable Banach
algebra with
$b_1 b_2 = 0$ for all $b_1, b_2 \in B$.
Let $a_1, a_2, \dots a_k, \dots$ be a countable sequence of
elements of norm 1 of $B$, with dense linear
span.
Let $A_0$
be the Fr\'echet algebra with 
standard Schwartz 
functions ${\cal S}_\gamma(X)$ as underlying Fr\'echet space,
and zero multiplication.

Define
a linear map $\theta \colon A_0 \rightarrow B$ by
$$\theta(\varphi) = \sum_{x \in X} \varphi(x) a_{\gamma(x)}, $$
$\varphi \in A_0$.  Since $\varphi$ is a Schwartz function,
and the $a_k$'s have norm 1,
this series converges absolutely to a well-defined element
of $B$.
We have $\| \theta(\varphi) \|_B
\leq \sum_{x \in X} |\varphi(x)| \| a_{\gamma(x)} \|_B
=  \| \varphi \|_1 $, where $\| \cdot \|_1 $ is the
$\ell^1$-norm, so $\theta$ is continuous,
and trivially an algebra homomorphism.
Usually
$\theta$ is not 1-1, unless for example $\{ a_k \}_{k=1}^\infty$
were a basis for $B$.  Let $A$ be the image $\theta(A_0)$ in $B$.
Then
$A$ is a dense Fr\'echet ideal in $B$, and is nuclear
since quotients by closed linear subspaces preserve 
nuclearity [Pietsch, 1972], Proposition 5.1.3
or [Treves, 1967], Proposition 50.1 (50.4).

\appendix

\vskip\baselineskip
\vskip\baselineskip
\section{Appendix. Refining the Ideal Condition, and $m$-Convexity}

We show that the constants $C_n$ can always be taken 
equal to $1$ in the dense ideal inequality (\ref{eqn:idealcond}),
and also $m_n = n$, by passing to an equivalent family
of seminorms.  We also note that any Fr\'echet algebra 
satisfying the ideal inequality 
is $m$-convex.\footnote{A Fr\'echet algebra $A$ is
{\it $m$-convex} if it can be topologized by a family
of {\it submultiplicative} seminorms, i.e. ones that
satisfy $\| a_1 a_2 \|_n\leq \| a_1\|_n \| a_2 \|_n$
for $a_1, a_2 \in A$.}

Let $A$ be a right Fr\'echet ideal in a Banach algebra $B$,
and let $\{ \| \cdot \|_n \}_{n=0}^\infty$ be an increasing
family of norms giving the topology for $A$, with
$\| \cdot \|_0 = \| \cdot \|_B$, as in Definition 3.1.
Define a new family of seminorms by
\begin{equation}
\| a \|^*_n = \sup \{\, \| a b \|_n \, | \, \| b \|_0 \leq 1,\, b \in B\, \},
\label{eqn:newsnorms}
\end{equation}
for $a \in A$.
Using the right ideal inequality (\ref{eqn:idealcond}) we see that
$\| a \|^*_n  \leq C_n\| a \|_m$, so the
topology given by $\{ \| \cdot \|^*_n\}_{n=0}^\infty$
is dominated by the original topology on $A$.
These new seminorms satisfy the right ideal inequality
with $m_n=n$ and $C_n=1$:
\begin{eqnarray}
\| a b \|^*_n & = & 
\sup \{ \, \| a b b_1 \|_n \, | \, \| b_1 \|_0 \leq 1, b_1 \in B\,  \}
\qquad \qquad \text{by definition (\ref{eqn:newsnorms})}\nonumber \\
& \leq & 
\sup \{ \, \| a b_2 \|_n \, | \, \| b_2 \|_0 \leq \| b \|_0, b_2 \in B\,  \}
\qquad \quad b_2 = b b_1, \text{ so } \|b_2\|_0 \leq \| b \|_0\nonumber \\
& = & 
\| b \|_0 \sup \{ \, \| a b_3 \|_n \, | \, \| b_3 \|_0 \leq 1, b_3 \in B\,  \}
\quad \quad\,  b_3 = b_2/\|b\|_0\nonumber \\
& = & 
\| a \|^*_n \| b \|_0  \qquad \qquad \qquad 
\qquad \qquad \qquad \qquad
\text{by definition (\ref{eqn:newsnorms}).}
\label{eqn:rightStar}
\end{eqnarray}
Note we used the submultiplicativity of the norm $\| \cdot \|_0$
on $B$ in the second step.
The inequality 
\begin{equation}
\| a b \|_n \leq \| a \|^*_n \| b\|_0
\label{eqn:rightHalfStar}
\end{equation}
can also be verified from the definition (\ref{eqn:newsnorms}).

Since $B$ may not be unital, the 
seminorms $\{\| \cdot \|^*_n\}_{n=0}^\infty$
are not necessarily equivalent to our original family.
(For example, $B$ could be a radical Banach algebra with
$a b = 0$ for all $a, b \in B$, in which 
case $\| a \|^*_n = 0$ for all $n$.)  
So define new norms by 
$\| a \|^{*+}_n = \max \{ \| a \|^*_n, \| a \|_n\}$,
$a \in A$.
Note that $\| \cdot \|^{*+}_0 = \| \cdot\|_0 = \| \cdot \|_B$, and
$\| a b \|^{*+}_n \leq \| a \|^*_n \| b \|_0$ by our estimates above.
The new family $\{ \| \cdot \|^{*+}_n\}_{n=0}^\infty$
gives a topology equivalent to the 
original family $\{ \| \cdot \|_n \}_{n=0}^\infty$,
and satisfies the right ideal inequality $\| a b \|^{*+}_n \leq
\| a \|^{*+}_n \| b \|_0$.
{\it Hence for any Banach algebra $B$ with right Fr\'echet ideal $A$, 
it is always possible to find an equivalent family
of norms giving the topology on $A$, such that the new norms 
satisfy the right Fr\'echet ideal inequality (\ref{eqn:idealcond}) 
with $C_n=1$ and $m_n= n$ for every $n \in \N$, and
the zeroth norm 
$\| \cdot \|_0$, which is the norm 
$\| \cdot \|_B$ on $B$, remains unchanged.}

The new norms  are submultiplicative for every $n$:
$\| a_1 a_2 \|^{*+}_n
 \leq \|a_1 \|^{*+}_n \| a_2\|^{*+}_0
 \leq \|a_1 \|^{*+}_n \| a_2\|^{*+}_n$,
$a_1, a_2 \in A$.
{\it Hence dense Fr\'echet ideals are always
$m$-convex Fr\'echet algebras.}

The same results hold for left ideals, by switching the order
in the above arguments.  In the left case, we define
$\| a \|^\dagger_n = 
\sup \{\, \|  b a \|_n \, | \, \| b \|_0 \leq 1,\, b \in B\, \}$,
and $\| a \|^{\dagger +}_n = \max \{ \| a \|^\dagger_n, \| a \|_n \}$.

Finally, assume that $A$ is a two-sided Fr\'echet ideal in 
the Banach algebra $B$.  
As before, $A$ is topologized by increasing norms 
$\{ \| \cdot \|_n \}_{n=0}^\infty$, with 
$\| \cdot \|_0 = \| \cdot \|_B$, and left and right
ideal inequalities are satisfied: 
$\| a b \|_n \leq C_n \| a \|_m \| b \|_0$,
$\| a b \|_m \leq C_m \| a \|_k \| b \|_0$,
and 
$\| b a \|_n \leq C_n \| b \|_0 \| a \|_m $,
$\| b a \|_m \leq C_m \| b \|_0 \| a \|_k$, 
for all $a \in A$ and $b \in B$.
Define a new family of seminorms by
\begin{equation}
\| a \|^{two}_n = \sup \{\, \| c a b \|_n \, | \, \| c \|_0, 
\| b \|_0 \leq 1,\,  c, b \in B\, \},
\label{eqn:new2snorms}
\end{equation}
for $a \in A$.
Using the right and left ideal inequalities we see that
$\| a \|^{two}_n  \leq C_n C_m \| a \|_p$, so the
topology given by $\{ \| \cdot \|^{two}_n\}_{n=0}^\infty$
is dominated by the original topology on $A$.
These new seminorms satisfy the right ideal inequalities
with $m_n=n$ and $C_n=1$:
\begin{eqnarray}
\| a b \|^{two}_n & = & 
\sup \{ \, \| c a b b_1 \|_n \, | \, \| c \|_0 , 
\| b_1 \|_0 \leq 1, \, c, b_1 \in B\,  \}
\qquad \qquad \text{by definition (\ref{eqn:new2snorms})}\nonumber \\
& \leq & 
\sup \{ \, \| c a b_2 \|_n \, | \,\| c \|_0 \leq 1, 
\| b_2 \|_0 \leq \| b \|_0,\,  c, b_2 \in B\,  \}
\quad b_2 = b b_1, \text{ so } \|b_2\|_0 \leq \| b \|_0\nonumber \\
& = & 
\| b \|_0 \sup \{ \, \| c a b_3 \|_n \, | \,\| c \|_0, 
\| b_3 \|_0 \leq 1,\,  c, b_3 \in B\,  \}
\quad \quad\,  b_3 = b_2/\|b\|_0\nonumber \\
& = & 
\| a \|^{two}_n \| b \|_0  \qquad \qquad \qquad 
\qquad \qquad \qquad \qquad \qquad \quad\,
\text{by definition (\ref{eqn:new2snorms}).}
\label{eqn:tworight}
\end{eqnarray}
Again we used the submultiplicativity of the norm $\| \cdot \|_0$
on $B$ in the second step.
Similarly, the left ideal inequality holds:
$\| b a \|^{two}_n \leq \| b \|_0 \| a \|^{two}_n$, $a \in A$,
$b \in B$.

Since the seminorms $\{ \| \cdot \|^{two}_n \}_{n=0}^\infty$
could give a weaker topology on $A$ (as we discussed in the
right ideal case above), we define our final set of norms on $A$
by 
$\| a \|^{two+}_n = \max \{ \|a \|^{two+}_n, \| a \|^*_n, 
\| a \|^\dagger_n, \| a \|_n\}$.
Note that $\| \cdot \|^{two+}_0 = \| \cdot\|_0 = \| \cdot \|_B$.
We have the four inequalities:
\begin{eqnarray}
\| a b \|^{two}_n & \leq & \| a \|^{two}_n \| b \|_0 \qquad \qquad \quad \qquad\qquad \text{by (\ref{eqn:tworight})} \nonumber \\
\| a b \|^{*}_n & \leq & \| a \|^*_n \| b \|_0 \qquad \qquad \qquad \qquad \qquad\text{by (\ref{eqn:rightStar})}\nonumber \\
\| a b \|^{\dagger}_n 
& = & \sup\{ \| c a b \|_n \, | \, \| c \|_0 \leq 1, c \in B\}\quad 
\text{by definition of $\| \cdot \|^\dagger_n$ above}\nonumber \\
& \leq  & \| a \|^{two}_n \| b \|_0 \qquad  \qquad \qquad \qquad \quad
\text{by definition (\ref{eqn:new2snorms})
of $\| \cdot \|^{two}_n$}\nonumber \\
\| a b \|_n & \leq  & \| a \|^*_n \| b \|_0, \qquad \qquad \qquad \qquad \quad \text{by (\ref{eqn:rightHalfStar})} \nonumber 
\end{eqnarray}
for $a \in A$ and $b \in B$.  Putting these together shows that 
the norms
$\{\| \cdot \|^{two+}_n \}_{n=0}^\infty$ satisfy the right ideal
inequality with $C_n=1$, $m_n = n$ for every $n\in \N$.  
Similarly, they satisfy the left ideal inequalities 
with the same constraints.

{\it We have shown that for any Banach 
algebra $B$ with two-sided Fr\'echet ideal $A$, 
it is always possible to find an equivalent family
of norms giving the topology on $A$, such that the new norms 
satisfy the right and left Fr\'echet ideal inequalities
with $C_n=1$ and $m_n= n$ for every $n \in \N$, and
the zeroth norm 
$\| \cdot \|_0$, which is the norm 
$\| \cdot \|_B$ on $B$, remains unchanged.}

\vskip\baselineskip
\vskip\baselineskip
\section{Appendix. Counterexamples}

\vskip\baselineskip
\noindent{\bf Example B.1.
Non-Summable Family of Scales That Dominates Every Power of an Enumeration,
with Regularity Satisfied.}
If $\gamma$ is an enumeration of $X$, define a family of 
scales by $\sigma_0 \equiv 1$ and $\sigma_n(x) = e^{\gamma(x)}$ for
$n \geq 1$.
Note that $\sigma_1=\sigma_2 = \cdots = \sigma_n = \cdots$,
so the summability condition
(\ref{eqn:NucSumSpace}) cannot be satisfied.
Also $\sigma_n/\sigma_m = 1$,
except at $n=0$ we have $\sigma_0/\sigma_m = 1/\sigma_m$ 
for all $m$.
Regularity is satisfied 
since $e^1 \leq e^2 \leq \cdots \leq e^{k} \leq \cdots$.
The family $\underline \sigma$ dominates every power of $\gamma$
since
$e^{\gamma(x)} \geq {\gamma(x)^k \over k!}$,
using the series expansion for $e$.

\vskip\baselineskip
\noindent{\bf Example B.2.
Enumerations are Not Minimal Scales.}
If $\gamma_1$ is an enumeration of $X$, one can find another 
enumeration $\gamma_2$ such that 
$\gamma_2 \lesssim \gamma_1$
(in fact $\gamma_2 \leq \gamma_1+1$)
but 
$\gamma_1 \not \lesssim \gamma_2$.
Let $\gamma_1$ be the identity enumeration of $X=\N^+$.
Let $\gamma_2$ be defined by $\gamma_2(k) = k + 1$ for
every $k$ except on
the set $\{\, \lceil e^{i^i} \rceil  \, | \, i \in \N \}$.\footnote{Here 
set $0^0=0$, 
so $e^{0^0}=1$, and $\lceil \cdot \rceil$
is the least integer greater than.}
Set
$\gamma_2(1) = 1$  and
 $\gamma_2(\lceil e^{(i+1)^{(i+1)}} \rceil) = 
\lceil e^{i^i} \rceil + 1$, $i \in \N$.
To be clear, $\gamma_2(1) = 1$, $\gamma_2(2) = 3$,  
$\gamma_2(3 = \lceil e^{1^1} \rceil )= 
\lceil e^{0^0} \rceil + 1 =2$, 
 $\gamma_2(4) = 5$,  $\gamma_2(5) = 6$, \dots 
$\gamma_2(54) = 55$, 
$\gamma_2(55 = \lceil e^{2^2} \rceil) = \lceil e^{1^1}\rceil + 1 = 4$,  
$\gamma_2(56) = 57$,  \dots $\gamma_2(\lceil e^{27} \rceil 
= \lceil e^{3^3} \rceil)
= \lceil e^{2^2} \rceil +1 
= 56$, $\gamma_2(\lceil e^{27} \rceil + 1) =\lceil e^{27} \rceil + 2$, 
\dots.    The function $e^{(i+1)^{(i+1)}}$
is not bounded by
$C e^{di^i}$ for $i \in \N^+$, for any fixed power $d$,
so 
$\gamma_1 \not \lesssim \gamma_2$.

\vskip\baselineskip
\noindent{\bf Example B.3. When Standard
Schwartz Functions are Not an Ideal.}   
Assume that 
$d_{\textrm{fin}}=\infty$
and 
$d_{\textrm{inf}}=0$ 
in \S 7.  We show explicitly that when the growth condition
(see Definition 7.2) fails, then
standard Schwartz functions on $Y=Y_{\textrm{fin}}$
is not an ideal in $B$.
Corollary 7.7 says this must be the case.
For an example when the growth condition fails,
take $p_k = e^{\bigl(e^{(e^k)}\bigr)}$, $k \in \N^+$.\footnote{
$e^{\bigl(e^{(e^{k+1})}\bigr)}
=e^{\bigl(e^{(e \cdot e^k)}\bigr)}
=e^{\bigl(e^{(e^k)}\bigr)^e}
> e^{\bigl(e^{(e^k)}\bigr)^2}$
so $\ln(p_{k+1}) > \ln(p_k)^2$, which contradicts 
$p_{k+1} \leq C(kp_k)^d$; $\ln(p_{k+1}) \leq C_1 + C_2 \cdot
\ln(p_k)$.}
Arrange that the sequence of dimensions 
 $\frak p = \{ p_k \}_{k=1}^\infty$ is ordered, so
$p_k \leq p_{k+1}$.
Let $\gamma(k, i, j) = 
p^2_1 + \cdots p^2_{k-1} + (i-1) + (j-1)p_k +1 $, 
$\{ k, i, j \} \in Y$ be as in Definition 7.2. 
%(where it was called $\gamma_{\textrm{fin}}$).
Then $\gamma$ is an enumeration of $Y$,
and ${\cal S}_{\gamma}(Y)$ is the Fr\'echet
space of standard Schwartz functions on $Y$.

%We show directly that 
%${\cal S}_{\gamma}(Y)$
% is not an ideal in $B$.
Define a scale $\beta$ on $Y$ by
$\beta(k, i, j) = ik p_{k-1} + (j-1)p_k$.
Then $\beta \thicksim \gamma$, and 
our calculations will simplify 
using $\beta$ in place of $\gamma$.
For each $k \in \N^+$, let 
$c_1(k)$  be
$1$'s in the first column:
$$ c_1(k) = \begin{pmatrix}
1 & 0 & \cdots & 0\\
\vdots & \vdots & \ddots & \vdots \\
1 & 0 & \cdots & 0 \end{pmatrix}
\in M_{p_k}(\C).  $$
Then $\| c_1(k) \|_{\text{op}} = \sqrt{p_k}$,
so 
$$ S_k =  {c_1(k) \over \sqrt{p_k}} \in c_f(X)$$
satisfies $\| S_k \|_B=1$.
The $n$th norm in ${\cal S}^1_\beta (X)$ is
\begin{eqnarray}
\| S_k \|_n  & = & 
 \sum_{i=1}^{p_k}
{\beta(k, i, 1)^n\over \sqrt{p_k}} \nonumber \\
&= & { k^n p^n_{k-1}  \over \sqrt{p_k}}
\sum_{i=1}^{p_k} i^n \nonumber  \\
& \geq &
{ k^n p^n_{k-1}  \over \sqrt{p_k}}
p_k^n \geq p_k^{n-1/2}. \nonumber 
\end{eqnarray}
And for each $k \in \N^+$, 
$$ T_k =  e_{k, 11} \in  c_f(X)$$
is an element of ${\cal S}_\beta (X)$ with
$m$th norm equal to 
\begin{equation}
\| T_k \|_m \quad = \quad   \beta(k, 1, 1)^m 
\quad = \quad   k^m p_{k-1}^m.
\end{equation}
Since $c_1(k) e_{k, 11} = c_1(k)$, $S_k T_k = S_k$.
The left ideal condition implies
that for each $n \in \N$, there is some $m \in \N$ for which
\begin{equation}
{\| S_k T_k \|_n \over \| S_k \|_B \| T_k \|_m} \quad =  \quad
{\| S_k  \|_n \over  \| T_k \|_m}  
\quad \geq \quad {p_k^{n-1/2} \over k^m p_{k-1}^m}
\label{eqn:blowideal}
\end{equation}
is bounded in $k$.
Taking $n = 2$, we get
\begin{equation}
p_k  \leq p_k^{2-1/2} 
\leq C k^m p_{k-1}^m 
\leq 2^m \cdot C \bigl(k-1 \cdot p_{k-1}\bigr)^m,
\end{equation}
which is the growth condition.

\vskip\baselineskip
\noindent{\bf Example B.4. Banach Algebra Not an Ideal in $c_0(X)$.}   
When the standard basis $\{ \delta_k \}_{k=0}^\infty$  
is an absolute basis for a dense Fr\'echet subspace $F$ of 
the commutative $C^\star$-algebra $c_0(\N)$, then $F$ is 
isomorphic to ${\cal S}^1_\sigma(\N)$ (Theorem 2.9), 
and is easily seen to be an ideal in $c_0(\N)$ by
the same argument as Example 3.3 (a).
We exhibit a dense Banach subalgebra $A$ of $c_0(\N)$,
for which $\{ \delta_k \}_{k=0}^\infty$ is not an absolute basis, 
and such that $A$ is not an ideal in $c_0(\N)$.  

Let $\sigma$ be a proper scale on $\N$.
Define a norm by
\begin{equation}
\| f \|_{\sigma,1} = \sup_{k \in \N}\,\biggl(\, \sigma(k) \, 
\max\, \bigl\{\,|f_+(k)|,
\, \,
\sigma(k) | f_{-}(k)|\, \bigr\}\, \biggr), 
\end{equation}
where $f_+(k) = f(2k) + f(2k+1)$,
$f_{-}(k) = f(2k) - f(2k+1)$, and $f \in c_f(\N)$.  
Define a new scale $\beta$, which is \lq\lq half of $\sigma$\rq\rq\ 
by $\beta(2k) = \beta(2k+1)=\sigma(k)$.
Define two related norms 
$\| f \|_\beta  =  \| \beta f \|_\infty$
and
$\| f \|_{\beta^2}  =  \| \beta^2 f \|_\infty$,
which topologize the Banach algebras
$c_0(\N, \beta)$ and $c_0(\N, \beta^2)$,
respectively.
Let $A$ be the completion of $c_f(\N)$ in the norm $\| \cdot \|_{\sigma, 1}$.
Then $c_f(\N) \subseteq c_0(\N,\beta^2) \hookrightarrow A \hookrightarrow
c_0(\N, \beta) \hookrightarrow c_0(\N)$, where the inclusions
$\hookrightarrow$  are continuous.
The Banach space $A$ is a Banach algebra since
\begin{equation}
\| f * g \|_{\sigma,1}\quad \leq \quad
\| f * g \|_{\beta^2}  \quad \leq \quad
\| f  \|_\beta \| g \|_\beta  \quad \leq \quad
\| f \|_{\sigma,1}
\| g \|_{\sigma,1},
\end{equation}
for $f,g \in c_f(\N)$.

Define $\delta_{+,k}=\delta_{2k} + \delta_{2k+1}$ 
and $\delta_{-,k}=\delta_{2k} - \delta_{2k+1}$.  Then
$\|\delta_{+,k} \|_{\sigma, 1} =  2\sigma(k)$,
$\|\delta_{-,k} \|_{\sigma, 1} =  2\sigma(k)^2$, 
$\|\delta_{\pm,k} \|_{\infty} = 1$, and $\delta_{+,1}*\delta_{-,1}
=\delta_{-,1}$.  So we have
\begin{equation}
{\| \delta_{+,k} * \delta_{-,k} \|_{\sigma,1} \over 
\| \delta_{+,k} \|_{\sigma,1} \| \delta_{-,k} \|_\infty}
=
{\|  \delta_{-,k} \|_{\sigma,1} \over 
\| \delta_{+,k} \|_{\sigma,1} \| \delta_{-,k} \|_\infty}
=
{2 \sigma(k)^2 \over 
2\sigma(k) \cdot 1 } = \sigma(k),
\end{equation}
which tends to $\infty$ as $k \rightarrow \infty$, since $\sigma$
was assumed proper.  So $A$ is not an ideal in $c_0(\N)$.

\vskip\baselineskip
\noindent{\bf Example B.5. Dense Nuclear Fr\'echet
Subalgebra of $c_0(X)$ Not an Ideal.}   
Let ${\cal A}$ be the Fr\'echet algebra with underlying
space ${\cal S}(\N)$\footnote{${\cal S}(\N)$ denotes standard Schwartz
functions ${\cal S}_\gamma(X)$
with $X=\N$ and $\gamma(k)=k+1$.} and multiplication
 $f * g (r) = \sum_{s=0}^{r} f(s) g(r-s)$, $r \in \N$,
for $f, g \in {\cal A}$.
Then $\| f* g \|_d \leq C_d \|f \|_d \| g \|_d$,
and ${\cal A}$ is a
commutative $m$-convex Fr\'echet algebra, with unit
$\delta_0$.
For any $t \in \N$, define the closed ideal
${\cal A}_t = \{ f \in {\cal A}| f(0) = f(1) = \cdots f(t) = 0 \}$
of ${\cal A}$.

Let $X$ be a countably infinite set, and
$\chi  \in c_0(X)$ have range in $(0,1)$.
Define a linear map
$\theta_\chi$ from ${\cal A}_0$
 to $c_0(X)$ by
\begin{equation}
\theta_\chi(f)(x) = \sum_{r=1}^\infty 
f(r) \chi(x)^r,
\end{equation}
$f \in {\cal A}_0$, $x \in X$.
Let $\epsilon > 0$ and find finite 
$S \subset X$ large enough so that
$| \chi(x)| < \epsilon$ if $x \in X - S$.
Then $| \theta_\chi(f)(x)|
\leq \epsilon \|f \|_1$,
for $x \in X - S$.
This proves
$\theta_\chi(f) \in c_0(X)$.
It is easy to see that $\| \theta_\chi(f) \|_\infty 
\leq \| f \|_1$, so $\theta_\chi$ is continuous.
For $f, g \in {\cal A}_0$ we have
\begin{eqnarray}
\bigl( \theta_\chi(f)* 
\theta_\chi(g) \bigr) (x) & = &
\sum_{r, s = 1}^\infty f(r) g(s) \chi(x)^{r + s}
\nonumber \\
&=& \sum_{r = 0}^\infty 
\sum_{s=0}^{r} f(s) g(r-s) \chi(x)^{r}
\nonumber \\
&=& \sum_{r = 1}^\infty 
 f * g(r) \chi(x)^{r}
\quad = \quad  \theta_\chi(f * g) (x),
\end{eqnarray}
$x \in X$.
So $\theta_\chi$ is an algebra homomorphism.

\noindent
{\bf Assume that $\chi$ has infinitely many values in it's range.}
We show that $\theta_\chi$ is injective.
Let $f \in {\cal A}_0$, and
note the power series $p(z) = \sum_{r=1}^\infty f(r) z^r$
converges for $|z| < 1$.
If $\theta_\chi(f) = 0$, then $p(z)$ is zero on the range of $\chi$.
Since $\chi$ is in $c_0(X)$ and has values in $(0,1)$, 
zero is an accumulation point.
By analyticity, $p(z)$ is identically zero, so $f \equiv 0$.

Moreover, the assumption implies 
$\theta_\chi({\cal A}_0) \cap c_f(X) = 0$.
If $\theta_\chi(f)$ takes on only finitely many values 
$w_1, \dots w_k$, then the power series
$(p(z) - w_1) \cdots (p(z) - w_k)$ is zero 
for $z$ in the range of $\chi$, and therefore identically zero,
by the argument of the previous paragraph.
Since $p(z)$ has no constant term, the product $w_1 \dots w_k$
is zero, so $w_i = 0 $ for some $i$.  If $p(z)$ is not
identically zero, we may divide by it to get
$(p(z) - w_1) \cdots (p(z) - w_{i-1})(p(z) - w_{i+1}) 
\cdots (p(z) - w_k) \equiv 0$, and we are in the same
situation as before.  Repeating the argument $k$ times
shows $p(z) - w_j \equiv 0$ for some $j \in \{ 1, \dots k \}$.
Finally $w_j=0$, contradicting our hypothesis that
$p(z) \not\equiv 0$.  Hence every non-zero
$\theta_\chi(f)$ takes on infinitely many values,
and cannot lie in $c_f(X)$.

\noindent
{\bf Assume in addition that $\chi(x) \not= \chi(y)$  
if $x$ and $y$ are 
distinct elements of $X$.}
Then $\theta_\chi({\cal A}_t)$ 
is dense in $c_0(X)$ for every $t \in \N$.
Let $\xi$ be
an element of the dual  $\ell^1(X) = c_0(X)'$, 
and assume $\xi$ vanishes on
$\theta_\chi({\cal A}_t)$.
Let $s \in \N$ be greater than $t$.
Taking $f = \delta_s \in {\cal A}_t$, 
we see that $\xi(\chi^s)=0$, 
so 
\begin{equation}
0 = \sum_{x \in X} \xi(x) \chi(x)^s.
\label{eqn:ChiSum}
\end{equation}
Find an ordering $\{x_i \}_{i=0}^\infty$
of $X$ for which
$\chi(x_0) > \chi(x_1) > \chi(x_2) > \cdots$.
Divide (\ref{eqn:ChiSum}) by $\chi(x_0)^s$
and let $s \rightarrow \infty$ to see
that $\xi(x_0) = 0$.  Then divide 
(\ref{eqn:ChiSum}) by
$\chi(x_1)^s$ to see $\xi(x_1)=0$, 
and continue this way through all the $x_i$'s
to see that $\xi(x_i) = 0$ for every $i$.
We have
proved that no non-zero 
element of $c_0(X)'$ 
can vanish on
$\theta_\chi({\cal A}_t)$,
which proves the density.

We have shown that $A_\chi=\theta_\chi({\cal A}_0)$
is a dense nuclear subalgebra of $c_0(X)$.
Since  $A_\chi$ does not contain $c_f(X)$, 
it cannot be an ideal in $c_0(X)$.  
Note that $A_\chi$ is not spectral invariant
in $c_0(X)$.  For example $\| \chi \|_\infty < 1$,
so the geometric series for 
the quasi-inverse $-\chi/(1+ \chi)$ converges in $c_0(X)$,
but not in $A_\chi$.  Alternately, for any $z_0 \in \C$, 
$| z_0| < 1$, $f \mapsto \sum_{r=0}^\infty f(r) z_0^r$
defines a 1-dimensional representation of ${\cal A}_0$,
and through $\theta_\chi$ a 1-dimensional
representation of $A_\chi$.  But only those
with $z_0$ in the range of $\chi$ can extend
to a representation of $c_0(X)$.

\vskip\baselineskip
\noindent{\bf Example B.6. Dense Nuclear Fr\'echet
Subalgebra of $c_0(X)$ Not an Ideal and Containing
$c_f(X)$.}   
Let $\underline \sigma$ be a family of scales on $X$ satisfying 
summability (\ref{eqn:NucSumSpace}).
Then $A_{\underline \sigma} = 
{\cal S}_{\underline \sigma} (X)$ is a dense nuclear
two-sided Fr\'echet ideal in 
the pointwise multiplication 
algebra $c_0(X)$, by Theorem 6.2.

Let $\chi \colon X \rightarrow (0, 1)$ 
satisfy the criteria of the previous Example B.5.
Arrange that no power of $\chi$ is in 
${\cal S}_{\underline \sigma}(X)$.
For example, we could take $\chi(x) = 1/\sigma_d(x)$ for 
some $d\in \N$, taking care to make sure
$\chi(x) \not= \chi(y)$ for distinct $x, y \in X$.
\footnote{$\chi(x) = 1/\gamma(x)$  works for the 
family $\sigma_n(x) = \gamma(x)^n$, where $\gamma$ 
is any enumeration of $X$.}
Let $A_\chi$ be the associated 
Fr\'echet subalgebra of $c_0(X)$ from
the previous Example B.5.    

We show $A_\chi \cap {\cal S}_{\underline \sigma}(X) =  0$.
Let $f \in {\cal A}_0$ and let $p\in \N^+$
be smallest such that $f(p)\not=0$.
Let $d \in \N$
be such that $\sigma_d(x)\chi(x)^p$
is unbounded for $x \in X$.
Let $S \subset X$ be a large enough finite set so 
that $\chi(x)\| f \|_1 <  
|f(p)|/2$ for $x\in X - S$.
Then 
\begin{eqnarray}
\sigma_d(x) \bigl|\theta_\chi(f)(x)\bigr| &=&
\sigma_d(x) \chi(x)^p \biggl| f(p) + 
\sum_{q=1}^\infty f(q+p) \chi(x)^q \biggr|
\nonumber \\
&\geq&
\sigma_d(x) \chi(x)^p
\biggl( |f(p)| - 
\biggl| \sum_{q=1}^\infty f(q+p) \chi(x)^q \biggr|
\biggr)
\nonumber \\
&\geq&
\sigma_d(x) \chi(x)^p
\biggl( |f(p)| - 
\|f \|_1 \chi(x) 
\biggr)
\nonumber \\
&\geq&
\sigma_d(x)\chi(x)^p
| f(p) | /2,
\end{eqnarray}
for $x \in X - S$.
Hence by the unboundedness of
$\sigma_d \chi^p$, 
$\theta_\chi(f)$ is not in ${\cal S}_{\underline \sigma}(X)$.
%, so $A_\chi \cap A_\sigma = 0$, 

The direct
sum $A_{\text{sum}} = A_\chi \oplus A_{\underline \sigma}$ is
naturally a Fr\'echet algebra, nuclear and
dense in $c_0(X)$.  It contains $c_f(X)$
as an ideal since $A_{\underline \sigma}$ 
does, but not densely.
In Example B.5, we noted 
the quasi-inverse $\psi = -\chi/(1- \chi)$  exists in
$c_0(X)$, but not in $A_\chi$.  
If $\psi \in A_{\underline \sigma}$, 
then $\chi \psi \in A_{\underline \sigma}$ and
$\chi \in A_{\underline \sigma}$, a contradiction.  So
$\psi \notin A_{\text{sum}}$, and
$A_{\text{sum}}$ is not 
spectral invariant, and not an ideal, in $c_0(X)$.

\vskip\baselineskip
\noindent{\bf Example B.7. Nilpotent Banach Algebra
with No Dense Nuclear Ideal.}
We exhibit a separable nilpotent
Banach algebra $B$ of order 2 ($B^3=\{0\}$) with no dense
nuclear Fr\'echet ideal.
Let $B$ be the Hilbert space direct sum of two
Hilbert spaces $B = {\cal H}_1 \oplus {\cal H}_2$, with
respective bases $\{ \alpha_k \}_{k=0}^\infty$,
$\{ \beta_k \}_{k=0}^\infty$.
Define multiplication
by $\alpha_0 \alpha_i = \beta_i$, and all
other products zero: $\alpha_{i+1} \alpha_j=0$,
$\beta_i \beta_j=0$, $\alpha_i \beta_j = \beta_j \alpha_i = 0$.
Let $\pi \colon {\cal H}_1 \cong {\cal H}_2$ be the isomorphism
of Hilbert spaces defined by
$\pi(\alpha_i) = \beta_i$.
The Hilbert space norm on $B$ is submultiplicative since
for
 $b_1=\xi_1+\eta_1$ and $b_2=\xi_2+ \eta_2$, 
%arbitrary elements of $B$,
%, where the $\xi$'s
%and $\eta$'s are in the first and second copies
%of $\ell^2(\N)$, respectively.
$$
\|  b_1 b_2\|_B
= \| \xi_1(0) \pi( \xi_2 ) \|_{{\cal H}_2}
= | \xi_1(0) |
\| \xi_2 \|_{{\cal H}_1}
\leq \| b_1\|_B
 \| b_2\|_B,
$$
so $B$ is a Banach algebra.
Clearly, $B$ is nilpotent of order 2.
Let $A$ be a
dense right Fr\'echet ideal in $B$.  By density,
we can't have $A \subseteq < \alpha_0 >^\perp$.  So let $a_0 \in A$ have
some component of $\alpha_0$ in the Hilbert space $B$.
Then ${\cal H}_2 \subset A$, since
$a_0 B = {\cal H}_2 $.
Rescale $a_0$ to arrange
that
$<a_0, \alpha_0> = 1$.
If $\eta \in {\cal H}_1$, then $\pi(\eta) = a_0 \eta$ 
and $\| \pi(\eta) \|_n = \| a_0 \eta \|_n
\leq C_n \| a_0 \|_m \| \eta \|_B = 
D_n \| \eta \|_{B}$.  So the topology on ${\cal H}_2$
inherited from $A$ is precisely the Hilbert space topology,
and $A$ cannot be nuclear (Proposition 2.1 (c)).

\vskip\baselineskip
\vskip\baselineskip
\section{Index}
\smallskip
\smallskip

\noindent
$A_{\textrm{fin}}$, 
$A_{\textrm{inf}}$ 
\dotfill  discussion after Definition 7.1

\noindent
Basis for a Fr\'echet space, coordinate functional
\dotfill Definition 2.5

\noindent
$\delta_k, \varphi_k, (\varphi_1, \varphi_2, \varphi_3, \dots)$
\dotfill Definition 2.11 

\noindent
Dominate $\lesssim$ and equivalence $\thicksim$ 
of scales, and of families of scales \dotfill Definition 2.6

\noindent
Enumeration, enumeration of $X$
\dotfill Definition 2.11

\noindent
Family of scales given by a single scale $\sigma$,
${\cal S}_{\sigma}^1(X)$, ${\cal S}_{\sigma}^\infty(X)$
\dotfill Definitions 2.6 and 2.7

\noindent
Finite support functions $c_f(X)$\dotfill Definition 2.7

\noindent
$\gamma$-regular family of scales or $\gamma$-regular basis
\dotfill Definition 2.11

\noindent
${\underline \gamma}$-family of scales
\dotfill Definition 2.11

\noindent
Growth condition
\dotfill Definition 7.2

\noindent
Ideal, Fr\'echet ideal
\dotfill Definition 3.1, Proposition 3.8

\noindent
$\ell_{\textrm{fin}}$, 
$\sigma_{\textrm{fin}}$,
$\zeta_{\textrm{inf}}$, 
${\frak p}_{\textrm{inf}}$ 
\dotfill  Definition 7.2

\noindent
$\N$, $\N^+$\dotfill  Introduction

\noindent
Nilpotent Banach algebra
\dotfill Example 8.6, Example B.7

\noindent
Not an ideal
\dotfill Examples B.3, B,4, B.5, B.6

\noindent
$\frak p$, $Z$, $Y$, $y=\{z, i, j\}$
\dotfill Introduction to \S 6

\noindent
$\frak p$-summable
\dotfill Theorem 6.2, Equation (\ref{eqn:NucSum})

\noindent
Power series space of infinite type\dotfill Remark 2.8

\noindent
Scale, family of scales \dotfill Definition 2.6

\noindent
Schauder, unconditional, equicontinuous, absolute basis\dotfill
Definition 2.5

\noindent
Schwartz Spaces
${\cal S}_{\underline \sigma}^1(X)$, 
${\cal S}_{\underline \sigma}^\infty(X)$, 
$\sup$-norm and $\ell^1$-norm\dotfill Definition 2.7

\noindent
standard Schwartz functions, Schwartz functions
\dotfill Equation (\ref{eqn:standSchw}), Definition 2.7

\noindent
Summability Condition, $\underline \sigma$ is summable,
$\sigma$ is summable
\dotfill Theorem 2.9, Equation (\ref{eqn:NucSumSpace})

\noindent
$X$
\dotfill Introduction and Definition 2.6

\noindent
$Y_{\textrm{fin}}$, 
$Z_{\textrm{fin}}$,
$Y_{\textrm{inf}}$, 
$Z_{\textrm{inf}}$ 
\dotfill  discussion after Definition 7.1

\vskip\baselineskip
\vskip\baselineskip
\section{References}
\smallskip
\smallskip
\footnotesize

\noindent [{\bf Crone Rob, 1975}]
\, L. Crone and W. B. Robinson, 
{\it Every nuclear Fr\'echet space with a regular basis has the 
quasi-equivalence property}, 
Studia Math. {\bf 52}, 
(1975), 203-207.

\noindent [{\bf Dixmier, 1982}]
\, J. Dixmier, {\it $C^\star$-algebras},
North-Holland Publishing Co., Amsterdam/New York/Oxford, 1982.

\noindent [{\bf Dubinsky, 1979}]
\, E. Dubinsky, {\it The Structure of Nuclear Fr\'echet Spaces}, 
Lect. Notes. Math. {\bf 720}, 
Berlin/Heidelberg/New York, Springer-Verlag, 1979.

\noindent [{\bf Fell Dor, 1988}]
\, J. M. G. Fell and R. S. Doran, 
{\it Representations of $^\star$-Algebras, Locally Compact Groups,
and Banach $^\star$-Algebraic Bundles, 
Volume I Basic Representation Theory of Groups and Algebras}, 
Pure and Applied Mathematics {\bf 125}, 
Academic Press, Boston, MA, 1988.

\noindent [{\bf Husain, 1991}]
\, T. Husain, 
{\it Orthogonal Schauder Bases},
Pure and Applied Mathematics, Volume 143,
Marcel Dekker, Inc., New York, 1991.

\noindent [{\bf Kad Ring II, 1997}]
\, R. V. Kadison and J. R. Ringrose, 
{\it Fundamentals of the Theory of Operator 
Algebras, Volume II: Advanced Theory},
Graduate Studies in Mathematics, Volume 16, 
American Mathematical Society, 1997.

\noindent [{\bf Kaplansky, 1948}]
\, I. Kaplansky, {\it Dual rings},
Ann. of Math. {\bf 49(3)} (1948), 689-701.

\noindent [{\bf Kaplansky, 1949}]
\, I. Kaplansky, {\it Normed algebras},
Duke Math. J. {\bf 16(3)} (1949), 399-418.

\noindent [{\bf Mars Hoff, 1999}]
\, J.E. Marsden and M.J. Hoffman, {\it Basic Complex Analysis},
W.H. Freeman, New York, 1999.

\noindent [{\bf Palmer, 1994}]
\, T. W. Palmer, {\it Banach Algebras and the General 
Theory of $\star$-algebras, 
Volume I: Algebras and Banach Algebras},
Encyclopedia of Mathematics and its Applications, Volume 49, 
Cambridge University Press, New York, 1994.

\noindent [{\bf Paterson, 1988}]
\, A.L.T. Paterson, {\it Amenability},
Mathematical Surveys and Monographs, Number 29, 
American Mathematical Society, Providence, Rhode Island, 1988.

\noindent [{\bf Pietsch, 1972}]
\, A. Pietsch, {\it Nuclear Locally Convex Spaces},
Ergebnisse Der Mathematick und Ihrer Grensgebiete, Volume 66, 
Springer-Verlag, New York/Heidelberg/Berlin, 1972.

\noindent [{\bf Rudin, 1973}]
\, W. Rudin, {\it Functional Analysis},
Series in Higher Mathematics,
McGraw-Hill, Inc., New York, 1973.

\noindent [{\bf Sch, 1998}]
\, L. B. Schweitzer, {\it $C^\infty$ functions on the Cantor
set, and a smooth $m$-convex Fr\'echet subalgebra of $O_2$},
Pac. J. Math. {\bf 184(2)} (1998), 349-365.

\noindent [{\bf Smyth, 1980}]
\, M. R. F. Smyth, {\it On problems of Olubummo and Alexander},
Proc. R. Ir. Acad. {\bf 80{\small A}(1)} (1980), 69-74.

\noindent [{\bf Sug, 1971}]
\, M. Sugiura, {\it Fourier series of smooth functions on compact Lie groups},
Osaka J. Math. {\bf 8} (1971), 33-47.

\noindent [{\bf Treves, 1967}]
\, F. Treves, {\it Topological Vector Spaces, 
Distributions and Kernels},
Academic Press, Inc, San Diego, California, 1967.

\vskip\baselineskip

\noindent{Web Page: \url{http://www.svpal.org/~lsch/Math/indexMath.html}.}

\end{document}